\documentclass[11pt]{amsart}
\usepackage{amsmath,amssymb,amsfonts}
\usepackage{times}
\usepackage[T1]{fontenc}
\usepackage[utf8]{inputenc}
\usepackage{lmodern}
\usepackage{color}
\usepackage[english]{babel}
\usepackage{geometry}
\usepackage{relsize}
\usepackage{booktabs}
\usepackage{mathtools}
\usepackage{amsthm}
\usepackage{amsmath,amssymb,graphicx,tikz,caption}
\usetikzlibrary{arrows.meta,calc}
\usepackage{emptypage}
\usepackage{xfrac}    
\usepackage{faktor}
\usepackage{listings}
\usepackage{lipsum}
\usepackage{setspace}
\AtBeginDocument{\singlespacing}
\usepackage{lipsum}
\usepackage{newlfont}
\usepackage[mathscr]{eucal}	
\usepackage{physics}
\usepackage{url}
\usepackage{hyperref}
\usetikzlibrary{positioning}
\usepackage{xcolor}
\usepackage[new]{old-arrows}
\usepackage{multirow}
\usepackage{enumitem}
\usepackage[toc]{appendix}
\numberwithin{equation}{section}
\def\eqref#1{(\ref{#1})}
\newlength{\alignedproofindent}
\NewDocumentEnvironment{alignedproof}{O{\proofname}m}
 {%
  \proof[#1]\mbox{}\par\nopagebreak
  \settowidth{\alignedproofindent}{#2\quad}%
  \everypar{\hangindent=\alignedproofindent\hangafter=0 }%
 }
 {\endproof}

\theoremstyle{plain}
\newtheorem{theorem}{Theorem}[section]
\newtheorem{lemma}[theorem]{Lemma}
\newtheorem{proposition}[theorem]{Proposition}
\newtheorem{corollary}[theorem]{Corollary}

\theoremstyle{definition}
\newtheorem{definition}[theorem]{Definition}

\newtheorem{example}[theorem]{Example}
\newtheorem{prop}[theorem]{Proposition}

\newtheorem{remark}{Remark}[section]

\usepackage{array}
\usepackage{booktabs}
\usepackage{tabularx}
\usepackage[signatures,write]{frontespizio}

\usepackage{fancyhdr}

\newcommand{\R}{\mathbb{R}} 
\newcommand{\C}{\mathbb{C}} 

\AtBeginDocument
{
\setlength\abovedisplayskip{.6\baselineskip}
\setlength\abovedisplayshortskip{.6\baselineskip}
\setlength\belowdisplayskip{.5\baselineskip}
\setlength\belowdisplayshortskip{.5\baselineskip}
}
\newcommand{\ee}{\mathrm{e}}
\newcommand{\ii}{\mathrm{i}}

\usepackage{slashed}
\usepackage{tikz-cd}
\usetikzlibrary{arrows.meta}

\newcommand{\BHEO}{BHE$_0$}
\newcommand{\BHER}{BHE$_{\mathbb{R}}$ }
\newcommand{\BHEI}{BHE$_{\ii\mathbb{R}}$ }
\newcommand{\Ht}{\mathbf{H}}

\newcommand{\AdS}{\mathrm{AdS}}
\newcommand{\SU}{\mathrm{SU}}

\newcommand{\vol}{\mathrm{vol}}

\newcommand{\epsp}{\epsilon_+}
\newcommand{\epsm}{\epsilon_-}
\newcommand{\epspm}{\epsilon_\pm}
\newcommand{\epsmp}{\epsilon_\mp}
\newcommand{\eps}{\epsilon}
\newcommand{\Jpm}{J^\pm}
\newcommand{\Jp}{J^+}

\newcommand{\fundp}{\omega^+}
\newcommand{\fundm}{\omega^-}
\newcommand{\fundpm}{\omega^\pm}

\newcommand{\holop}{\Omega^+}

\newcommand{\holopm}{\Omega^\pm}

\newcommand{\fundt}{\omega^T}

\newcommand{\holot}{\Omega^T}

\newcommand{\vvre}{V_1}
\newcommand{\vvim}{V_2}
\newcommand{\vform}{\eta}
\newcommand{\vformbar}{\bar{\eta}}
\newcommand{\vformre}{\eta_1}
\newcommand{\vformim}{\eta_2}
\newcommand{\leep}{\theta^+}
\newcommand{\leem}{\theta^-}
\newcommand{\leepm}{\theta^\pm}
\newcommand{\leet}{\theta^T}

\newcommand{\rhok}{\rho_K}
\newcommand{\Rk}{R_K}

\newcommand{\Rep}{\mathrm{Re}}
\newcommand{\Imp}{\mathrm{Im}}
\newcommand{\cf}{\text{cf. }}

\usepackage{xcolor}
\begin{document}

\title{Spinorial geometry of Bismut--Hermitian--Einstein manifolds}

\author{Beatrice Brienza}
\address[Beatrice Brienza]{Instituto Nacional de Matem\'atica Pura e Aplicada (IMPA) \\ Estrada Dona Castorina, 110, CEP 22460-320 \\
Rio de Janeiro, RJ (Brasil)}
\email{brienza.beatrice@impa.br}

\author{Matteo Kevin Crisafio}
\address[Kevin Crisafio]{Dipartimento di Matematica ``G. Peano'', Universit\`{a} degli studi di Torino \\
Via Carlo Alberto 10\\
10123 Torino, Italy}
\email{matteokevin.crisafio@unito.it}

\author{Anna Fino}
\address[Anna Fino]{Dipartimento di Matematica ``G. Peano'', Universit\`{a} degli studi di Torino \\
Via Carlo Alberto 10\\
10123 Torino, Italy\\
\& Department of Mathematics and Statistics, Florida International University\\
Miami, FL 33199, United States}
\email{annamaria.fino@unito.it, afino@fiu.edu}

\author{Dario Martelli}
\address[Dario Martelli]{Dipartimento di Matematica ``G. Peano'', Universit\`{a} degli studi di Torino \\
Via Carlo Alberto 10\\
10123 Torino, Italy}
\email{dario.martelli@unito.it}



\maketitle

\begin{abstract}
\noindent 
We present a pair of spinorial equations, 
generalizing the internal Killing
spinor equations of Type I supergravity, and show that  they characterize several classes 
of Bismut–Hermitian–Einstein structures on six-dimensional Riemannian spin manifolds. We show that in this framework the local reduction to a four-dimensional transverse Hermitian geometry  is provided by two canonical, commuting, Bismut-parallel  Killing vector fields arising as bilinears in the Killing spinors, without assuming compactness. In one class, the transverse geometry is conformally K\"ahler 
and both compact
and non-compact solutions exist.
We also discuss
variants that do not admit compact solutions, one of which 
 has hyper-Hermitian transverse geometry. 
We discuss the relation of these structures to the Hull–Strominger system with closed three-form and to five-dimensional Gauntlett–Kim geometries and provide explicit embeddings in Type
I, Type IIB and massive Type IIA supergravity, as well as in generalized supergravity.
 We illustrate our results with a number of examples, both compact and non-compact.
\end{abstract}

\tableofcontents  

\section{Introduction}

It has been known for almost half a century that the study of the solutions of supergravity theories yields very rich geometric structures on Riemannian, as well as pseudo-Riemannian (in Lorentzian signature) manifolds, that are of independent interest to the mathematics community. Supergravity combines General Relativity, that is the most geometric of the theories of Physics, with the idea of supersymmetry, extending the setting of Riemannian geometry to that of  \emph{spinorial geometry}. The story kicks off with the realization that, in the simplest set-up, supersymmetry implies the existence of a parallel  
spinor on a six-dimensional manifold, which in turn is equivalent to this being a Calabi--Yau \cite{Candelas:1985en}. More generally, parallel spinors imply a reduction of the holonomy of the Levi-Civita connection to a subgroup of the Spin group in any given dimension. This  motivates the relevance of special holonomy manifolds for physics,  either compact \cite{Shatashvili:1994zw,Papadopoulos:1995da} or non-compact \cite{Eguchi:1978xp,Gibbons:1989er,Atiyah:2001qf}.

Generically supergravity theories comprise a number of fields in addition to the metric, corresponding to the existence of extra geometric data on the manifolds. One instance is the so-called Heterotic supergravity, that includes a three-form $H$, a function $\phi$ (called the dilaton field) as well as a connection on a certain principal bundle (a Yang--Mills field). This theory is defined on a ten-dimensional  
Lorentzian manifold, and assuming this is the direct product of four-dimensional Minkowski space-time with a six-dimensional Riemannian manifold, the supersymmetry equations are equivalent to a system of differential equations originally formulated by Strominger \cite{str} and Hull \cite{Hull} in the mid '80s. These references showed  that one can view the 
 geometry as describing a Hermitian, but non-K\"ahler, manifold endowed with a 
 metric connection (Bismut) with totally antisymmetric torsion. 
 Using the spinorial formulation it is straightforward  to show that the Ricci form of this connection must vanish; in other words that the geometry describes Hermitian manifolds with vanishing Bismut Ricci form. 
 
 This set up  was revisited at  the beginning of 2000s in the light of the framework of $G$-structures and showed  \cite{LopesCardoso:2002vpf,Gauntlett:2003cy} to correspond to a six-dimensional manifold endowed with an $\SU(3)$ structure  with intrinsic torsion lying in one of the classes of \cite{chiossi2002intrinsictorsionsu3g2}.
In addition to presenting an  equivalent reformulation of the  Hull--Strominger set  up\footnote{For example, the conformally balanced condition was first highlighted in this reference.},
   reference  \cite{Gauntlett:2003cy}  discussed also different dimensions and structure groups, thus for example formulating  $\mathrm{G}_2$ and  $\mathrm{Spin}(7)$ analogs.
However, the six-dimensional case remained the most studied one, attracting considerable attention from the differential geometry community, starting from the key results of 
\cite{li2006existencesupersymmetricstringtheory} and \cite{fu2006theorysuperstringfluxnonkahler}.

A simplified system can be obtained by setting to zero the Yang--Mills field and requiring that the three-form $H$ is closed: this is known in the physics literature as Type I supergravity. 
The closedness of $H$ is equivalent to the pluriclosed condition on the Hermitian metric, thus leading to a
Bismut--Hermitian--Einstein structure. It is however important to remark that the conformally balanced condition is not affected by restricting to the simplified Type I system, and as a consequence 
one can show that no compact complete solutions can exist  with non-zero $H$ \cite{Gauntlett:2003cy}. See also \cite{IP3} for an earlier no-go theorem.  
Some  simple classes of complete examples were constructed in \cite{Gauntlett:2003cy}, generalizing the solutions in \cite{Callan:1991dj}. In  
\cite{Maldacena:2009mw} (building on  \cite{Butti:2004pk}) it was discussed a construction, with underlying topology $T^* S^3$, which was argued to realize explicitly 
the idea  of geometric transition \cite{Vafa_2001}, related to Reid's fantasy \cite{Reid1987}. Very recently, the existence of the solution discussed in
\cite{Maldacena:2009mw} has been proved  in \cite{podesta2026completepluriclosedmetricsvanishing} and \cite{foscolo2026cohomogeneitysolutionsiib}.

Lifting the conformally balanced condition, compact solutions become allowed. However, this leads one outside Type I supergravity equations, losing the spinorial geometry description. 
Notwithstanding this, in the mathematics literature, Bismut--Hermitian--Einstein structures, in any dimensions $2n$,  have been studied  \emph{assuming} that the underlying manifold is 
 compact, an assumption mostly motivated by the pluriclosed flow \cite{StreetsTian}, whose fixed points are precisely Bismut--Hermitian--Einstein metrics \cite{Ye}. 
 In particular, in \cite{ALS} (building on  \cite{GFJS}) it was shown that on compact smooth BHE manifolds there must exist two commuting Bismut-parallel Killing vector fields
 yielding a canonical $\mathbb{T}^2$-reduction to a transverse geometry in two dimensions less. Correspondingly, the Bismut connection must have (restricted) holonomy contained in $\SU(n-1)$, instead of $\SU(n)$. A natural generalization of the Bismut--Hermitian--Einstein condition to the hypercomplex realm is given by strong HKT manifolds. For further references, see \cite{Brienza:2025zjc,Papadopoulos:2025ugn}.  
 Remarkably, in six dimensions, the transverse four-dimensional geometry must be (locally) conformally K\"ahler, with a curvature condition that leads to a PDE that appeared 
 in the physics literature almost twenty years earlier in \cite{Kim:2005ez,Kim:2006qu,Gauntlett:2007ts}.  This arose in a completely different set up,  in the context of characterizing systematically the geometry 
 of supersymmetric solutions including an anti de Sitter (AdS) factor, motivated by the AdS/CFT correspondence. See \emph{e.g.} \cite{Gauntlett:2025vnh}.
 In this context the conformally K\"ahler structure describes the transverse geometry of $(2n-1)$-dimensional 
 total space, that  is a circle bundle over this\footnote{In the quasi-regular case.}.  The physically relevant dimensions are seven ($n=4$) \cite{Kim:2005ez} and nine ($n=5$) \cite{Kim:2006qu},
although   \cite{Gauntlett:2007ts} described an extension valid in arbitrary dimension\footnote{Strictly speaking the dimension five 
($n=3$) was excluded in  \cite{Gauntlett:2007ts}, while our results will cover this case. See Section \ref{Section7}.}, and following this these geometries have been named ``GK geometries''.

Despite the transverse curvature condition is a complicated fourth order PDE, explicit local solutions of co-homogeneity one and two were presented in 
 \cite{Donos:2008ug} and  \cite{Donos:2008hd}, respectively   (see also \cite{Gauntlett:2006ns}). In \cite{Couzens:2018wnk} it was shown that these local metrics can be extended 
 globally on compact five-dimensional manifolds  with topology of $S^2\times S^3$, 
named\footnote{Because of their close relationship with the five-dimensional Sasaki--Einstein manifolds of \cite{Gauntlett:2004yd} and \cite{Martelli:2005wy}.}
$\mathscr{Y}^{p,q}$ and $\mathscr{L}^{\mathrm{a,b,c}}$.   Using the results  of \cite{ALS}, the above examples lead to explicit compact examples of six-dimensional 
BHE manifolds with topology $S^3\times S^2 \times S^1$. More generally, the ``transgressive'' solutions of \cite{Donos:2008ug} can be uplifted to 
BHE manifolds with $S^3\times S^3$ topology, as discussed in \cite{ALL-26}.  Very recently \cite{apostolov2026gluingconstructionzextremalkahler} employed a  gluing construction 
to produce  genuinely new compact examples of BHE structures on manifolds  with various topologies. 



From a mathematical standpoint, in this paper we propose a spinorial system that extends the BHE set up studied \cite{KS-26,ALS}. This system is a deformation of the Type I supersymmetry equations
that lifts the conformally balanced condition and allows for compact manifolds.  We show that our system  yields the symmetry (and holonomy)  reduction discussed in \cite{ALS},
\emph{without assuming compactness}, thus opening up the study of natural non-compact BHE  manifolds. As an illustration, exploiting the cohomogeneity one local solutions presented in 
\cite{Donos:2008ug}, we construct  non-compact complete examples, with topology $S^3\times S^1 \times \mathbb{R}^2$. We also describe additional generalizations of the spinorial equations that lead us to BHE structures of distinct type with respect to those considered in \cite{KS-26,ALS,ALL-26}, where the transverse four-dimensional geometry is not conformally K\"ahler. 
In these broader contexts compact manifolds are ruled out (see Theorem \ref{prop:imaginary-kappa-noncompact}), however  we have not been able to construct complete non-compact examples (unless when the underlying hyper-Hermitian structure is locally conformally hyper-K\"ahler), which therefore remain an open problem.

From a physical standpoint,  in addition to the already quoted GK geometries giving rise to AdS$_3$ solutions of Type IIB supergravity   \cite{Kim:2005ez}  and AdS$_2$ solutions of eleven dimensional supergravity \cite{Kim:2006qu}, we discuss 
other embeddings, in Type I and massive Type IIA supergravity, as well as  generalized supergravity \cite{Arutyunov:2015mqj,Wulff:2016tju}. 
In particular, we note  that in massive Type IIA supergravity there exist supersymmetric solutions of the form AdS$_3\times I \times $BHE$_6$, where $I$ is an interval and  BHE$_6$ are (compact or non compact) six dimensional BHE manifolds of the type discussed in \cite{KS-26,ALS}.
 Indeed, these appeared in \cite{Couzens:2022agr}, although the BHE structure was not emphasised there. Moreover, we show that in generalized supergravity there exist supersymmetric solutions of the form $\mathbb{R}^{1,3} \times $ BHE$_6$, where BHE$_6$ are a new type of non-compact BHE manifolds, characterized by a transverse four-dimensional hyper-Hermitian structure. The generalized supergravity  proposed in \cite{Arutyunov:2015mqj,Wulff:2016tju}  arises from requiring that the two-dimensional string world-sheet theory is only scale invariant, as opposed to the more stringent requirement of conformal invariance, that leads to the standard supergravity theories. Although some solutions of this generalized supergravity have been studied in the physics literature  (see e.g. \cite{Orlando:2016qqu,Arutyunov:2015mqj,Sakamoto:2018krs,Seibold:2019dvf,Hong:2018tlp}), somewhat surprisingly a systematic study of the solutions of this theory, especially in the presence of supersymmetry, has not been undertaken so far. Our results thus represent a first  step in this direction, starting from the ``classic'' ten-dimensional direct product Ansatz, of four dimensional Minkowski space-time, with a six-dimensional Riemannian spin manifold. 

We defer a more detailed discussion of the different embeddings in supergravity 
and their mutual relationships to the final section of the paper and Appendix \ref{sugraembeddings}.

\vskip 3mm

\noindent 
{\bf Note added:} after the main results of this work were completed, we became aware of
the non-compact examples of \cite{podesta2026newinhomogeneousbismuthermitianeinstein}, overlapping with the examples we discuss in Section \ref{section:examples}. 

\smallskip

\section{Bismut--Hermitian--Einstein geometry}\label{subsec:BHE-structures}

In this Section we recall the preliminaries on Bismut--Hermitian--Einstein geometry. We recall here the notations and conventions used throughout the paper. Each manifold is assumed to be connected, unless otherwise specified.

Let $(M,J,g)$ be an Hermitian manifold. The fundamental form $\omega$ is defined by $\omega(X,Y)=g(X,JY)$. The Bismut (or Strominger) connection $\nabla^B$ is the unique Hermitian connection with totally skew torsion \cite{Bi,str}. Its expression is given by 
\begin{equation}
    \nabla^B=\nabla +\tfrac12 H,
\end{equation}
where $\nabla$ is the Levi-Civita connection of $g$ and $H=J \dd\omega$. Note that the complex structure $J$ acts on forms as follows: for any $\alpha \in \Omega^k(M)$, 
\begin{equation}
    J\alpha(X_1,\dots,X_k)=\alpha(-JX_1, \dots, -JX_k).
\end{equation}
A Hermitian structure is said to be pluriclosed, whenever $\dd H=0$.
\par
We adopt the following convention on the curvature: $R^B=[\nabla^B,\nabla^B]-\nabla^B_{[\cdot,\cdot]}$. The associated Bismut Ricci form is given by 
\begin{equation}
    \rho^B=\tfrac12 \sum_{a=1}^{2n} R^B(\cdot,\cdot, e_a, Je_a),
\end{equation}
where $\{e_i\}$ is any orthonormal basis. \par
We will adopt the convention that for a $k$-form $\alpha$ its norm will be given by summation over all the indexes.  

Bismut--Hermitian--Einstein structures, introduced in
\cite{GFJS,GFS}, provide a natural non-K\"ahler counterpart of
Ricci-flat K\"ahler geometry. On compact manifolds, the theory is
particularly rigid in low dimensions: the surface case was
classified in \cite{GI}, the structure theory in real dimension six
was developed in \cite{ALS}, and further topological restrictions
were established in \cite{Streets}.
We recall the following 
\begin{definition}\label{dfn:grs}
Let $(M,g)$ be a Riemannian manifold and let $H\in\Omega^3(M)$
be a closed three-form. The pair $(g,H)$ is called a
\emph{gradient generalized Ricci soliton} if there exists
a smooth function $f$ such that
\begin{equation}\label{eq:gradient-generalized-Ricci-soliton}
 \operatorname{Rc}-\frac14H^2+\operatorname{Hess}f=0,
 \qquad
 \dd^*H+\iota_{\operatorname{grad}f}H=0.
\end{equation}
Here $\operatorname{Rc}$ and $\operatorname{Hess}$ are computed
using the Levi-Civita connection.  When \(M\) is compact, the soliton potential \(f\) is unique up to an additive constant. We henceforth fix this ambiguity by requiring $\int_M \ee^{-f}\,\dd V_g=1,$ and call the resulting function the \emph{normalized soliton potential}.
\end{definition}
\begin{definition} \label{defn:GPS}
Let $(M^{2n},J,g)$ be a pluriclosed Hermitian manifold.
    The Hermitian structure $(J,g)$ is called
    \emph{Bismut--Hermitian--Einstein}, or \emph{BHE}, if
    $\rho^{B}=0$, where $\rho^{B}$ is the Bismut Ricci form. 
\end{definition}

Let $(M^{2n},J,g)$ be a compact BHE manifold. By
\cite[Proposition~2.6]{GFJS}, it is a  gradient generalized Ricci soliton. Denote by $f$ its normalized soliton potential and set
\begin{equation}\label{eq:canonical-BHE-vector-field}
    V_1:=\tfrac12\bigl(\theta^\sharp-\mathrm{grad} f\bigr),
\end{equation}
where $\theta^\sharp$ is the vector field dual to the Lee form $\theta$, defined by the relation $\dd \omega^{n-1}=\theta \wedge \omega^{n-1}$.

The structure results of \cite{ALS,GFJS} give
\begin{equation}\label{eq:BHE-parallel-fields}
    \nabla^{B} V_1=0,
    \qquad
    \mathcal{L}_{V_1} J=0.
\end{equation}
In particular, $V_1$ and $J V_1=V_2$ are commuting real holomorphic
Killing vector fields of constant norm, which, up to rescaling, is assumed to be $1$. Moreover, the soliton
potential is invariant under both vector fields, that is,
$V_1(f)= V_2(f)=0$ \cite[Lemma~2.9]{ALS}. Finally, $V_1$ vanishes
if and only if $(J,g)$ is K\"ahler. This was first proved in
\cite[Proposition~2.6]{GFJS}; a shorter proof was subsequently given
in \cite[Proposition~2.6]{ALS}. In this case, the BHE condition
reduces to the Ricci-flat K\"ahler condition.

Suppose henceforth that $(J,g)$ is non-K\"ahler. We use the
normalization
\begin{equation}\label{eq:BHE-unit-normalization}
    |V_1|=|V_2|=1.
\end{equation}
Then $\mathcal V_J:=\operatorname{span}\{V_1, V_2\}$ defines a
rank-two holomorphic Riemannian foliation. Unless otherwise stated,
\emph{basic} will mean basic with respect to the foliation
$\mathcal{V}_J$. Introduce $\eta_i:=V_i^\flat$. The pair $(\eta_1,\eta_2)$ defines an invariant connection for the foliation. Its
curvature forms $\dd\eta_i$ are basic and $J$-invariant
\cite[Lemma~2.7]{ALS}. Thus every compact non-K\"ahler BHE manifold
is locally a $\R^2$-bundle over a Hermitian manifold of real
dimension $2n-2$.

Let $\mathcal H_J:=\ker \eta_1\cap\ker\eta_2=\mathcal V_J^\perp$ be the
horizontal distribution. The metric and the fundamental form decompose
as
\begin{equation}\label{eq:BHE-transverse-decomposition}
    g
    =\eta_1\otimes\eta_1
     +\eta_2\otimes\eta_2
     +g^T,
    \qquad
    \omega
    =-\eta_1\wedge\eta_2+\omega^T.
\end{equation}
The tensors $g^T$, $J^T:=J|_{\mathcal H_J}$, and $\omega^T$ are
basic and define a transverse Hermitian structure. The torsion $H$ has the following expression
\begin{equation} \label{eqn:torsionH}
   H= \eta_1 \wedge \dd\eta_1 + \eta_2 \wedge \dd\eta_2+ H^T 
\end{equation}
with $H^T$ being the Bismut torsion of the transverse Hermitian structure. 
\begin{prop}[cf. {\cite[Proposition~2.10]{ALS}}]
\label{prop:compact-BHE-transverse-equations}
Given $(M,J,g)$ a compact (non K\"ahler) BHE manifold, the
transverse geometry satisfies
\begin{align}
    &\theta^T=\dd f ,
    \label{eq:BHE-transverse-Lee}\\
    &\operatorname{tr}_{\omega^T}\dd \eta_1=0,
    \\
    &\operatorname{tr}_{\omega^T}\dd \eta_2=4,
    \label{eq:BHE-curvature-traces}\\
    &\dd \dd_{J^T}^{\mathrm{c}}\omega^T
    =-\dd \eta_1 \wedge \dd \eta_1
      -\dd \eta_2 \wedge \dd \eta_2 .
    \label{eq:BHE-anomaly-equation}
\end{align}
Moreover,
\begin{equation}\label{eq:BHE-transverse-Ricci}
    \dd\eta_2=\tfrac12\rho^{B,T},
\end{equation}
where $\rho^{B,T}$ is the Bismut Ricci form of the transverse Hermitian structure. 
\end{prop}

Equation \eqref{eq:BHE-transverse-Lee} implies that the transverse
Hermitian structure is conformally balanced. In real dimension six,
this implies that $\omega_{K}^T:=\ee^{-f}\omega^T$ is a transverse
K\"ahler form. In complex dimension two the Bismut Ricci form is
conformally invariant. Hence \eqref{eq:BHE-curvature-traces} and
\eqref{eq:BHE-transverse-Ricci} yield the following dimensional
reduction and converse construction. Notice that the numerical
constants differ from those in \cite[Proposition~2.18]{ALS}, since
our unit vector field is one half of the vector field used there.

\begin{theorem}[cf. {\cite[Proposition~2.18]{ALS}}]
\label{thm:BHE-six-to-four}
Let $(M^6,J,g)$ be a compact non-K\"ahler BHE manifold, and let
$(M^T,g^T,J^T)$ be any local leaf space of $\mathcal V_J$, endowed
with the induced transverse Hermitian structure. Then $g_K:=\ee^{-f}g^T$
is K\"ahler. Denote its K\"ahler form, Ricci form, and scalar
curvature by $\omega_K$, $\rho_K$, and $R_K$, respectively.
Then:
\begin{enumerate}
    \item $\dd\eta_2=\tfrac12 \rho_K$ and 
    \begin{equation}\label{eq:BHE-transverse-scalar-curvature}
        R_K=8\ee^f>0.
    \end{equation}

    \item
    The form $\alpha:=\dd\eta_1$ is a closed primitive real $(1,1)$-form
    satisfying
    \begin{equation}\label{eq:BHE-four-dimensional-equation}
        \tfrac12\bigl(\dd \dd_{J^T}^{\mathrm{c}}R_K\bigr)\wedge\omega_K
        =-4\alpha \wedge\alpha
         -\rho_K \wedge\rho_K.
    \end{equation}
\end{enumerate}
Conversely, let $(S,J^T,g_K,\omega_K)$ be a K\"ahler
orbifold surface with strictly positive scalar curvature,
and let $\alpha$ be a closed primitive real $(1,1)$-form
satisfying \eqref{eq:BHE-four-dimensional-equation}.
Assume that
\begin{equation}
\operatorname{span}_{\mathbb R}
\bigl([\alpha],c_1^{\mathrm{orb}}(S,J^T)\bigr)
\subset H^2_{\mathrm{orb}}(S,\mathbb R)
\end{equation}
admits a basis consisting of integral orbifold cohomology
classes. Let $p:M\to S$ be a principal
$\mathbb T^2$-orbibundle with connection
one-forms $\eta_1,\eta_2$ satisfying
\begin{equation}
\dd\eta_1=\alpha,
\qquad
\dd\eta_2=\tfrac12 \rho_K.
\end{equation}
Define
\begin{equation}\label{eq:BHE-four-to-six-data}
g^T:=\frac{R_K}{8}g_K,
\qquad
f:=\log\!\left(\frac{R_K}{8}\right).
\end{equation}
These data define a BHE structure on $M$ by
\begin{equation}
g:=\eta_1^2+\eta_2^2+\ee^f g_K,
\qquad
\omega_J:=-\eta_1\wedge\eta_2+\ee^f \omega_K.
\end{equation}
 Here $V_1$ and $V_2$ are the vertical vector
fields dual to $\eta_1$ and $\eta_2$, respectively, and
$J$ is induced by $J^T$ on the horizontal distribution
and satisfies $J(V_1)=V_2$ on the vertical distribution.
If $S$ is compact and the total space $M$ is smooth,
the resulting structure is a compact non-K\"ahler
BHE threefold.
\end{theorem}
Note that taking the Hodge star of Equation \eqref{eq:BHE-four-dimensional-equation}, gives the following scalar equation 
\begin{equation}\label{eq:real-kappa-scalar-master-equation}
-\Delta_KR_K+\frac12R_K^2
 -|\operatorname{Rc}_K|^2
 -4| \dd \eta_1|^2=0,
 \qquad
 \Delta_K=-\dd_K^*\dd,
\end{equation}
commonly referred as the \emph{Box equation}.

\section{Spinors and spinorial bilinears}\label{sec:spinor-preliminaries}

We collect here the spin-geometric conventions and algebraic identities used
in the sequel.  Throughout, Hermitian products are conjugate-linear in the
first entry and complex-linear in the second.  Standard references for the
material in this section include \cite{LM,Fr,BFGK,Tomasiello_2022}.

\subsection{The Clifford module and complex conjugation}

Let $(M^6,g)$ be an oriented Riemannian spin manifold.  Its complex spinor
bundle $S\to M$ is the rank-eight complex vector bundle associated with the Dirac
representation
\begin{equation}
 \rho_\Delta\colon\operatorname{Spin}(6)\longrightarrow
 \operatorname{GL}_{\mathbb C}(\Delta_6),
 \qquad
 \Delta_6=\Delta_6^+\oplus\Delta_6^-,
 \qquad
 \dim_{\mathbb C}\Delta_6^\pm=4.
\end{equation}
Fix the isomorphism $\operatorname{Spin}(6)\cong\operatorname{SU}(4)$ for
which
\begin{equation}
 \Delta_6^-\cong\mathbb C^4,
 \qquad
 \Delta_6^+\cong\overline{\mathbb C^4}.
\end{equation}
With this choice, the identification
$\Delta_6\cong\overline{\mathbb C^4}\oplus\mathbb C^4$ turns the Dirac
action of $A\in\operatorname{SU}(4)$ into
\begin{equation}
 \rho_\Delta(A)(\overline v,u)=(\overline{Av},Au).
\end{equation}
Sections of $S$ will be called \emph{Dirac spinors}.
\medskip \par 
The Clifford multiplication by one-forms is the bundle map
\begin{equation}
 T^*M\otimes S\longrightarrow S,
 \qquad
 \alpha\otimes\psi\longmapsto\alpha\cdot\psi.
\end{equation}
We adopt the positive Clifford convention: for all real one-forms
$\alpha,\beta$ and all spinors $\psi$,
\begin{equation}\label{eq:spin-clifford-global}
 \alpha\cdot\beta\cdot\psi+
 \beta\cdot\alpha\cdot\psi
 =2g^{-1}(\alpha,\beta)\psi.
\end{equation}
Choose a local oriented orthonormal frame $\{e_a\}_{a=1}^6$, with dual coframe $\{e^a\}_{a=1}^6$ and a corresponding local spin lift. In the induced trivialization of $S$, the Clifford multiplication by $e^a$ is represented by the \emph{gamma matrices}
\begin{equation}
 e^a\cdot\psi=\gamma^a\psi.
\end{equation}
In particular, equation \eqref{eq:spin-clifford-global} becomes
\begin{equation}\label{eq:spin-clifford-local}
 \gamma^a\gamma^b+\gamma^b\gamma^a
 =2\delta^{ab}\operatorname{Id}_{S}.
\end{equation}

The Clifford action of a wedge product is obtained by antisymmetrization. For one-forms $\alpha_1,\ldots,\alpha_k$, we set
\begin{equation}
 (\alpha_1\wedge\cdots\wedge\alpha_k)\cdot\psi
 :=\frac{1}{k!}\sum_{\tau\in\mathfrak S_k}
 \operatorname{sgn}(\tau)\,
 \alpha_{\tau(1)}\cdot\ldots\cdot\alpha_{\tau(k)}\cdot\psi.
\end{equation}
In a local orthonormal coframe, 
\begin{equation}\label{eq:spin-antisymmetrized-gamma}
 \gamma^{a_1\cdots a_k}
 :=\gamma^{[a_1}\cdots\gamma^{a_k]}
 =\frac{1}{k!}\sum_{\tau\in\mathfrak S_k}
 \operatorname{sgn}(\tau)\,
 \gamma^{a_{\tau(1)}}\cdots\gamma^{a_{\tau(k)}},
\end{equation}
hence, if $ \beta=\frac{1}{k!}\beta_{a_1\cdots a_k}
 e^{a_1}\wedge\cdots\wedge e^{a_k} $
is a differential $k$-form, then
\begin{equation}\label{eq:spin-form-action-local}
 \beta\cdot\psi
 =\frac{1}{k!}\beta_{a_1\cdots a_k}
 \gamma^{a_1\cdots a_k}\psi.
\end{equation}

The metric identification between vectors and covectors also lets vector fields act on spinors: $X\cdot\psi:=X^\flat\cdot\psi.$
Thus, if $X=X^ae_a$, then
\begin{equation}\label{eq:spin-vector-action-local}
 X\cdot\psi=X^a\gamma_a\psi.
\end{equation}
\medskip 
We identify \(\Delta_6\) with \(\mathbb C^8\) and equip it with the standard Hermitian product $\langle\cdot,\cdot\rangle$ defined as $\langle\varphi,\psi\rangle=\varphi^\dagger \psi$. 
We shall use a \emph{Majorana realization} of the Clifford action on \(\Delta_6\): namely, we fix an orthonormal basis of \(\Delta_6\) in which the matrices \(\gamma^a\) are purely imaginary and skew-symmetric under transpose, that is
\begin{equation}\label{eq:spin-majorana-matrices}
 \overline{\gamma^a}=-\gamma^a,
 \qquad
 (\gamma^a)^T=-\gamma^a,
 \qquad
 (\gamma^a)^\dagger=\gamma^a.
\end{equation}
The third property follows from the first two, since
$(\gamma^a)^\dagger=(\overline{\gamma^a})^T=\gamma^a$. This implies also that the Clifford
multiplication by a real one-form is self-adjoint. \par
\medskip 
We also fix the globally defined $\operatorname{Spin}(6)$-equivariant antilinear map
\begin{equation}
 S\longrightarrow S,
 \qquad
 \psi\longmapsto\psi^*,
\end{equation}
which plays the role of complex conjugation on spinors. Locally, in the basis of \(\Delta_6\simeq\mathbb C^8\) in which the gamma matrices satisfy \eqref{eq:spin-majorana-matrices}, the map \(^{*}\) is component-wise complex conjugation. For every real one-form, it satisfies the following:
\begin{equation}\label{eq:spin-charge-conjugation}
 (\psi^*)^*=\psi,
 \qquad \langle\psi^*,\phi^*\rangle=\langle\phi,\psi\rangle, \qquad (\alpha\cdot\psi)^*=-\alpha\cdot\psi^*.
\end{equation}

Combining the Hermitian product with $^*$ gives the
$\operatorname{Spin}(6)$-invariant complex symmetric bilinear form
\begin{equation}\label{eq:spin-complex-bilinear-form}
 \mathsf B(\psi,\phi):=\langle\psi^*,\phi\rangle.
\end{equation}
In these local trivializations, $ \mathsf B(\psi,\phi)=\psi^T \phi$. \par

It will be useful to record the signs produced by adjoint and transpose.  Set
\begin{equation}
s_k:=(-1)^{\frac{k(k-1)}{2}},
\qquad
t_k:=(-1)^{\frac{k(k+1)}{2}}.
\end{equation}
For every increasing multi-index $(a_1,\ldots,a_k)$, one has
\begin{equation}\label{eq:spin-adjoint-transpose-kform}
(\gamma^{a_1\cdots a_k})^\dagger
=s_k\gamma^{a_1\cdots a_k},
\qquad
(\gamma^{a_1\cdots a_k})^T
=t_k\gamma^{a_1\cdots a_k}.
\end{equation}
Taking either the adjoint or the transpose reverses the order of the $k$ matrices, which contributes $k(k-1)/2$ transpositions.  Transposition also contributes one minus sign for each gamma matrix.  Thus, for every real $k$-form $\beta$,
\begin{equation}\label{eq:spin-B-symmetry-kform}
\mathsf B(\phi,\beta\cdot\psi)
=t_k\mathsf B(\psi,\beta\cdot\phi).
\end{equation}

\subsection{Chirality}

The oriented volume form defines the chirality operator
\begin{equation}\label{eq:spin-chirality-global}
 \gamma_*\psi:=\ii^3\operatorname{vol}_g\cdot\psi.
\end{equation}
In a local oriented orthonormal frame this becomes
\begin{equation}\label{eq:spin-chirality-local}
 \gamma_*=\ii^3\gamma^{123456}.
\end{equation}
The basic identities satisfied by $\gamma_*$ are
\begin{equation}\label{eq:spin-chirality-identities}
 \gamma_*^2=\operatorname{Id}_{S},
 \qquad
 \gamma_*^\dagger=\gamma_*,
 \qquad
 \gamma_*(X\cdot\psi)=-X\cdot\gamma_*\psi.
\end{equation}

Thus
\begin{equation}\label{eq:spin-chiral-splitting}
 S=S^+\oplus S^-,
 \qquad
 \psi=\psi_++\psi_-,
 \qquad
 \psi_\pm=\frac12(\operatorname{Id}_{S}\pm\gamma_*)\psi,
\end{equation}
where $S^\pm$ are the eigenbundles of $\pm 1$, respectively.
Self-adjointness of $\gamma_*$ makes this splitting orthogonal.  A section of
either summand is called a \emph{Weyl spinor}. Pointwise we have the (non-canonical) \emph{unitary} identification $S_p^\pm \cong \Delta^\pm_6$.

Complex conjugation exchanges the two chiralities.  Applying the third identity in
\eqref{eq:spin-charge-conjugation} six times gives
\begin{equation}
 (\gamma^{123456}\psi)^*
 =(-1)^6\gamma^{123456}\psi^*
 =\gamma^{123456}\psi^*.
\end{equation}
Since $^*$ is antilinear, $\overline{\ii^3}=-\ii^3$, and therefore
\begin{equation}
 (\gamma_*\psi)^*
 =(\ii^3\gamma^{123456}\psi)^*
 =-\ii^3\gamma^{123456}\psi^*
 =-\gamma_*\psi^*.
\end{equation}
For $\psi_\pm\in S^\pm$, conjugating
$\gamma_*\psi_\pm=\pm\psi_\pm$ now gives
\begin{equation}
 \pm\psi_\pm^*
 =(\gamma_*\psi_\pm)^*
 =-\gamma_*\psi_\pm^*.
\end{equation}
Hence
\begin{equation}\label{eq:spin-C-exchanges-chirality}
 \gamma_*\psi_\pm^*=\mp\psi_\pm^*,
 \qquad
 (S^\pm)^*=S^\mp.
\end{equation}
In particular, no non-zero Weyl spinor is fixed by $^*$.

The same parity argument shows that, for a $k$-form $\beta$,
\begin{equation}\label{eq:spin-form-chirality-parity}
 \gamma_*(\beta\cdot\psi)
 =(-1)^k\beta\cdot\gamma_*\psi.
\end{equation}
Thus even forms preserve chirality, while odd forms reverse it.

\subsection{Metric connections and spinorial bilinears}

Let $D$ be a metric connection on $TM$, possibly with torsion.  Its connection
form on the oriented orthonormal frame bundle lifts through the fixed spin
structure to a connection $D^S$ on $S$.  We use the same symbol $D$ for the
lifted connection.

We shall use the following standard compatibility relations.  For vector
fields $X,Y$ and spinors $\psi,\phi$,
\begin{align}
 D_X(Y\cdot\psi)
 &=(D_XY)\cdot\psi+Y\cdot D_X\psi,
 \label{eq:spin-connection-clifford-compatibility}\\
 X\langle\psi,\phi\rangle
 &=\langle D_X\psi,\phi\rangle
 +\langle\psi,D_X\phi\rangle.
 \label{eq:spin-connection-metric-compatibility}
\end{align}
The maps $^*$ and $\gamma_*$ are parallel as well:
\begin{equation}\label{eq:spin-star-chirality-parallel}
 D_X(\psi^*)=(D_X\psi)^*,
 \qquad
 D\gamma_*=0.
\end{equation}
In particular, every metric spin connection preserves the splitting
$S=S^+\oplus S^-$. The
bilinear form $\mathsf B$ is parallel in the corresponding sense:
\begin{equation} \label{eqn:metric_B}
    X\mathsf B(\psi,\phi)
 =\mathsf B(D_X\psi,\phi)+\mathsf B(\psi,D_X\phi).
\end{equation}

For $\psi,\phi\in\Gamma(S)$, define the associated differential forms by
\begin{align}
 \beta_k^{\mathrm H}(\psi,\phi)(X_1,\ldots,X_k)
 &:=\left\langle\psi,
 (X_1^\flat\wedge\cdots\wedge X_k^\flat)\cdot\phi
 \right\rangle,\label{eq:spin-H-bilinear-global}\\
 \beta_k^{\mathrm B}(\psi,\phi)(X_1,\ldots,X_k)
 &:=\mathsf B\left(\psi,
 (X_1^\flat\wedge\cdots\wedge X_k^\flat)\cdot\phi
 \right).\label{eq:spin-B-bilinear-global}
\end{align}
Equivariance of Clifford multiplication and of both pairings makes these
definitions tensorial.  In an oriented orthonormal frame they read
\begin{align}
 \beta_k^{\mathrm H}(\psi,\phi)
 &=\frac{1}{k!}
 \langle\psi,\gamma_{a_1\cdots a_k}\phi\rangle
 e^{a_1\cdots a_k},\label{eq:spin-H-bilinear-local}\\
 \beta_k^{\mathrm B}(\psi,\phi)
 &=\frac{1}{k!}
 \mathsf B(\psi,\gamma_{a_1\cdots a_k}\phi)
 e^{a_1\cdots a_k}.
 \label{eq:spin-B-bilinear-local}
\end{align}

Differentiating the two bilinears and using
\eqref{eq:spin-connection-clifford-compatibility},
\eqref{eq:spin-connection-metric-compatibility}, and
\eqref{eq:spin-star-chirality-parallel} gives
\begin{align}
 D_X\bigl(\beta_k^{\mathrm H}(\psi,\phi)\bigr)
 &=\beta_k^{\mathrm H}(D_X\psi,\phi)
 +\beta_k^{\mathrm H}(\psi,D_X\phi),
 \label{eq:spin-derivative-H-bilinear}\\
 D_X\bigl(\beta_k^{\mathrm B}(\psi,\phi)\bigr)
 &=\beta_k^{\mathrm B}(D_X\psi,\phi)
 +\beta_k^{\mathrm B}(\psi,D_X\phi).
 \label{eq:spin-derivative-B-bilinear}
\end{align}

\subsection{The $\operatorname{SU}(3)$-structure defined by a Weyl spinor}

Fix $\sigma\in\{+,-\}$ and a nowhere-vanishing unit Weyl spinor
$\epsilon_\sigma\in\Gamma(S^\sigma)$. The unit-norm assumption is only a normalization: for spinors
of constant positive norm, the same identities hold after the
corresponding constant rescaling of the spinorial bilinears.  It determines a two-form
$\omega^\sigma$ and a complex three-form $\Omega^\sigma$ by
\begin{align}
 \omega^\sigma(X,Y)
 &:=-\ii\left\langle\epsilon_\sigma,
 (X^\flat\wedge Y^\flat)\cdot\epsilon_\sigma
 \right\rangle,
 \label{eq:spin-SU3-omega-global}\\
 \Omega^\sigma(X,Y,Z)
 &:=-\mathsf B\left(\epsilon_\sigma,
 (X^\flat\wedge Y^\flat\wedge Z^\flat)\cdot\epsilon_\sigma\right).
 \label{eq:spin-SU3-Omega-global}
\end{align}
Locally, their components are
\begin{align}
 \omega^\sigma
 &=\frac12\omega^\sigma_{ab}e^{ab},
 &
 \omega^\sigma_{ab}
 &=-\ii\epsilon_\sigma^\dagger\gamma_{ab}\epsilon_\sigma,
 \label{eq:spin-SU3-omega-local}\\
 \Omega^\sigma
 &=\frac{1}{3!}\Omega^\sigma_{abc}e^{abc},
 &
 \Omega^\sigma_{abc}
 &=-\epsilon_\sigma^T\gamma_{abc}\epsilon_\sigma.
 \label{eq:spin-SU3-Omega-local}
\end{align}
Let us check explicitly that $\omega^\sigma$ is real.  Setting $k=2$ in
\eqref{eq:spin-adjoint-transpose-kform} gives
\begin{equation}
 \gamma_{ab}^\dagger=-\gamma_{ab}.
\end{equation}
Conjugate symmetry and the definition of the adjoint then give
\begin{equation}
 \overline{\left\langle\epsilon_\sigma,
 \gamma_{ab}\epsilon_\sigma\right\rangle}
 =\left\langle\gamma_{ab}\epsilon_\sigma,
 \epsilon_\sigma\right\rangle
 =\left\langle\epsilon_\sigma,
 \gamma_{ab}^\dagger\epsilon_\sigma\right\rangle
 =-\left\langle\epsilon_\sigma,
 \gamma_{ab}\epsilon_\sigma\right\rangle.
\end{equation}
Thus $\langle\epsilon_\sigma,
\gamma_{ab}\epsilon_\sigma\rangle$ is purely imaginary, and consequently
\begin{equation}
 \overline{\omega^\sigma_{ab}}
 =\ii\,\overline{\left\langle\epsilon_\sigma,
 \gamma_{ab}\epsilon_\sigma\right\rangle}
 =-\ii\left\langle\epsilon_\sigma,
 \gamma_{ab}\epsilon_\sigma\right\rangle
 =\omega^\sigma_{ab}.
\end{equation}

We use $\omega^\sigma$ to define $J^\sigma\in\operatorname{End}(TM)$ through
\begin{equation}\label{eq:spin-J-convention-global}
 \omega^\sigma(X,Y)=g(X,J^\sigma Y).
\end{equation}
Our convention for the action of $J^\sigma$ on complex one-forms is
\begin{equation}\label{eq:spin-J-action-one-forms}
 (J^\sigma\alpha)(X):=-\alpha(J^\sigma X).
\end{equation}
Accordingly, if $T^{1,0}_{J^\sigma}M$ denotes the $\ii$-eigenbundle on
$TM\otimes\mathbb C$, then
\begin{equation}\label{eq:spin-J-types-one-forms}
 \Lambda^{1,0}_{J^\sigma}
 =\{\alpha:J^\sigma\alpha=-\ii\alpha\},
 \qquad
 \Lambda^{0,1}_{J^\sigma}
 =\{\alpha:J^\sigma\alpha=\ii\alpha\}.
\end{equation}

\begin{remark}
\label{rem:spin-SU3-model}
We record an explicit model that fixes the constants and signs used below. Let $E$ be an oriented six-dimensional Euclidean vector space. Fix an oriented orthonormal basis $\{e_1,\ldots,e_6\}$ of $E$, and denote by $\{e^1,\ldots,e^6\}$ its dual basis.
Let
\begin{equation}
 I=\begin{pmatrix}1&0\\0&1\end{pmatrix},
  \qquad
 T_1=\sigma_1=\begin{pmatrix}0&1\\1&0\end{pmatrix},
 \qquad
 T_2=-\ii\sigma_2=\begin{pmatrix}0&-1\\1&0\end{pmatrix},
 \qquad
 T_3=\sigma_3=\begin{pmatrix}1&0\\0&-1\end{pmatrix}.
\end{equation}
Let $\{f_1,f_2\}$ be the standard basis of $\mathbb C^2$, and identify
$\Delta_6\cong\mathbb C^8
 \cong
 \mathbb C^2\otimes\mathbb C^2\otimes\mathbb C^2$ by ordering its standard
basis as
\begin{equation}
 u_{4(a-1)+2(b-1)+c}
 :=f_a\otimes f_b\otimes f_c,
 \qquad a,b,c\in\{1,2\}.
\end{equation}
For $P,Q,R\in\operatorname{End}(\mathbb C^2)$, the tensor-product
endomorphism $P\otimes Q\otimes R$ acts on decomposable tensors by
\begin{equation}
 (P\otimes Q\otimes R)(v_1\otimes v_2\otimes v_3)
 :=Pv_1\otimes Qv_2\otimes Rv_3,
\end{equation}
and its action is then extended linearly.  For example,
$(I\otimes T_2\otimes T_1)u_1=u_4$. 

With this convention, define
\begin{align*}
 \gamma_1&=\ii(I\otimes T_2\otimes T_1),
 &
 \gamma_2&=\ii(T_1\otimes T_2\otimes T_3),
 &
 \gamma_3&=\ii(T_3\otimes T_2\otimes T_3),\\
 \gamma_4&=\ii(T_2\otimes T_1\otimes T_1),
 &
 \gamma_5&=\ii(T_2\otimes I\otimes T_3),
 &
 \gamma_6&=\ii(I\otimes I\otimes T_2).
\end{align*}
These matrices are purely imaginary, skew-symmetric under transpose,
self-adjoint, and satisfy
\begin{equation}
 \gamma_a\gamma_b+\gamma_b\gamma_a=2\delta_{ab}\operatorname{Id}.
\end{equation}
We define the Clifford multiplication by declaring
\begin{equation}
 e^a\cdot\psi:=\gamma_a\psi,
 \qquad
 \psi\in\Delta_6.
\end{equation}
Thus $\gamma_a$ is the matrix of the Clifford multiplication by $e^a$
with respect to the spinor basis $\{u_1,\ldots,u_8\}$.
Consider the unit spinors
\begin{equation}
 \epsilon_+^0=\frac{1}{\sqrt2}(u_2+\ii u_5),
 \qquad
 \epsilon_-^0=\frac{1}{\sqrt2}(u_1-\ii u_6).
\end{equation}
A direct calculation gives
\begin{equation}
 \gamma_*\epsilon_+^0=\epsilon_+^0,
 \qquad
 \gamma_*\epsilon_-^0=-\epsilon_-^0.
\end{equation}
Let $(\omega_0^\sigma,\Omega_0^\sigma)$ be the bilinears defined by
\eqref{eq:spin-SU3-omega-global}--\eqref{eq:spin-SU3-Omega-global} from
$\epsilon_\sigma^0$.  In the coframe $\{e^a\}$,
they are
\begin{align*}
 \omega_0^\sigma
 &=e^{12}+e^{34}-\sigma e^{56},\\
 \Omega_0^\sigma
 &=-\sigma(e^5-\sigma \ii e^6)
 \wedge(e^1+\ii e^2)\wedge(e^3+\ii e^4).
\end{align*}
\end{remark}

\begin{proposition}\label{prop:spin-Weyl-SU3}
The tensors $(g,J^\sigma,\omega^\sigma,\Omega^\sigma)$ define an
$\operatorname{SU}(3)$-structure on $M$.  More precisely, $J^\sigma$ is an
orthogonal almost-complex structure, $\Omega^\sigma$ is a
nowhere-vanishing form of type $(0,3)$ with respect to $J^\sigma$, and
\begin{equation}\label{eq:spin-SU3-algebraic-identities}
 \omega^\sigma\wedge\Omega^\sigma=0,
 \qquad
 \frac{1}{3!}(\omega^\sigma)^3
 =\frac{\ii}{8}\Omega^\sigma\wedge\overline{\Omega^\sigma}
 =-\sigma\operatorname{vol}_g.
\end{equation}
In particular, $J^+$ induces the chosen orientation of $M$, while $J^-$
induces the opposite orientation.
\end{proposition}

\begin{proof}
For any $p \in M$ the group $\operatorname{Spin}(6)\cong\operatorname{SU}(4)$ acts transitively
on the unit sphere in $S^\sigma_p \cong \C^4$, with stabilizer
$\operatorname{SU}(3)$.  Fix the model spinor $\epsilon_\sigma^0$ of Remark
\ref{rem:spin-SU3-model}: under the pointwise identification introduced above, there
exists an element of $\operatorname{Spin}(6)$ sending
$\epsilon_\sigma(p)$ to the model spinor $\epsilon_\sigma^0$. By equivariance of the bilinear constructions,
there is an oriented orthonormal coframe of $T_p^*M$ in which
\begin{align}
 \omega^\sigma
 &=e^{12}+e^{34}-\sigma e^{56},
 \label{eq:spin-SU3-normal-form-omega}\\
 \Omega^\sigma
 &=-\sigma(e^5-\sigma \ii e^6)
 \wedge(e^1+\ii e^2)\wedge(e^3+\ii e^4).
 \label{eq:spin-SU3-normal-form-Omega}
\end{align}

Reading $J^\sigma$ from
\eqref{eq:spin-SU3-normal-form-omega} using
\eqref{eq:spin-J-convention-global}, we find
\begin{equation}
 \begin{aligned}
 J^\sigma e_1&=-e_2, & J^\sigma e_2&=e_1,
 &J^\sigma e_3&=-e_4, & J^\sigma e_4&=e_3,\\
 J^\sigma e_5&=\sigma e_6,
 &J^\sigma e_6&=-\sigma e_5.
 \end{aligned}
\end{equation}
This already proves that $(J^\sigma)^2=-\operatorname{Id}$ and that
$J^\sigma$ is orthogonal.  It also shows that
$e_1+\ii e_2\in T^{1,0}_{J^\sigma}M$, whereas
$e_1-\ii e_2\in T^{0,1}_{J^\sigma}M$.  Since
$(e^1+\ii e^2)(e_1+\ii e_2)=0$ and
$(e^1+\ii e^2)(e_1-\ii e_2)=2$, the form $e^1+\ii e^2$ is of type $(0,1)$.
The convention \eqref{eq:spin-J-action-one-forms} similarly gives
$J^\sigma(e^3+\ii e^4)=\ii(e^3+\ii e^4)$ and
$J^\sigma(e^5-\sigma \ii e^6)=\ii(e^5-\sigma \ii e^6)$.
By \eqref{eq:spin-J-types-one-forms}, each factor in
\eqref{eq:spin-SU3-normal-form-Omega} has type $(0,1)$; hence
$\Omega^\sigma$ has type $(0,3)$.  Expanding the normal forms gives the wedge
identities and the normalization in
\eqref{eq:spin-SU3-algebraic-identities}.  They also show that $J^+$ induces
the chosen orientation of $M$, whereas $J^-$ induces the opposite one.  We
have therefore recovered the tensors of the $\operatorname{SU}(3)$-reduction
determined by $\epsilon_\sigma$.
\end{proof}

The group-theoretic point behind the proposition is worth making
explicit. For every $p\in M$, the stabilizer of a unit vector in
either $S_p^+$ or $S_p^-$ fixes its span pointwise and acts as
$\operatorname{SU}(3)$ on its orthogonal complement.  See
\cite{ACFJ,GMSW,Gauntlett:2003cy} for further details on six-dimensional
$\operatorname{SU}(3)$-structures defined by spinors.

\subsection{Generic Dirac spinors and the induced $\operatorname{SU}(2)$ structure} \label{sub:static}

Now let
\begin{equation}
 \epsilon=\epsilon_++\epsilon_-,
 \qquad
 \epsilon_\pm\in\Gamma(S^\pm).
\end{equation}
Pointwise, an element of $\operatorname{Spin}(6)$ fixes $\epsilon$ precisely
when it fixes both chiral components.  Since $^*$ is equivariant, it fixes
$\epsilon_+$ if and only if it fixes $\epsilon_+^*$.  The two spinors
\begin{equation}
 u:=\epsilon_+^*,
 \qquad
 v:=\epsilon_-
\end{equation}
lie in the same four-dimensional half-spin module, and
\begin{equation}\label{eq:spin-Dirac-stabilizer}
 \operatorname{Stab}(\epsilon)
 =\operatorname{Stab}(u)\cap\operatorname{Stab}(v).
\end{equation}
This description also makes the possible stabilizers transparent.  If $u$
and $v$ are linearly dependent and not both zero, their common stabilizer is
$\operatorname{SU}(3)$; if both vanish, it is all of
$\operatorname{Spin}(6)$.  For a linearly independent pair, identify $S_p^-$ with $\mathbb C^4$.  An element fixing both vectors is the identity on their span and acts as $\operatorname{SU}(2)$ on its orthogonal
complement.  We call $\epsilon$ \emph{generic} when $\epsilon_+^*$ and
$\epsilon_-$ are pointwise linearly independent; equivalently, its pointwise
stabilizer is $\operatorname{SU}(2)$.

We now restrict to the following particular generic case. We assume
that $\epsilon_+$ and $\epsilon_-$ are nowhere vanishing and satisfy
\begin{equation}\label{eq:spin-static-SU2-condition}
 \langle\epsilon_-,\epsilon_+^*\rangle=0,
\end{equation}
which is equivalent to $\mathsf{B}(\epsilon_+,\epsilon_-)=0$.
 This is usually referred to as the \emph{static}
$\operatorname{SU}(2)$ case. We will further assume that $\epsilon_\pm$ have unit norm, so that $(\epsilon_+^*,\epsilon_-)$ is an ordered orthonormal pair in $S^-$.

We define the relevant tensors before choosing a frame.  Let $\eta$ be the
complex one-form
\begin{equation}\label{eq:spin-eta-global}
 \eta(X):=\langle\epsilon_-,X^\flat\cdot\epsilon_+\rangle,
\end{equation}
and let $V\in\Gamma(TM\otimes\mathbb C)$ be its metric dual with respect to
the complex-bilinear extension of $g$:
\begin{equation}\label{eq:spin-V-global}
 g(V,X)=\eta(X).
\end{equation}
Define
\begin{align}
 \omega^T
 &:=\frac12(\omega^++\omega^-),
 \label{eq:spin-omegaT-global}\\
 \Omega^T(X,Y)
 &:=\mathsf B\left(\epsilon_+,
 (X^\flat\wedge Y^\flat)\cdot\epsilon_-\right).
 \label{eq:spin-OmegaT-global}
\end{align}
Here $(\omega^+,\Omega^+)$ and $(\omega^-,\Omega^-)$ denote the
$\operatorname{SU}(3)$-structures defined by $\epsilon_+$ and $\epsilon_-$,
respectively, as in
\eqref{eq:spin-SU3-omega-global}--\eqref{eq:spin-SU3-Omega-global}.

Locally, these definitions read
\begin{align}
 \eta&=\eta_a e^a,
 &
 \eta_a&=\epsilon_-^\dagger\gamma_a\epsilon_+,
 &
 V^a&=\epsilon_-^\dagger\gamma^a\epsilon_+,
 \label{eq:spin-eta-V-local}\\
 \omega^T&=\frac12\omega^T_{ab}e^{ab},
 &
 \omega^T_{ab}
 &=-\frac{\ii}{2}\left(
 \epsilon_+^\dagger\gamma_{ab}\epsilon_+
 +\epsilon_-^\dagger\gamma_{ab}\epsilon_-
 \right),
 \label{eq:spin-omegaT-local}\\
 \Omega^T&=\frac12\Omega^T_{ab}e^{ab},
 &
 \Omega^T_{ab}
 &=\epsilon_+^T\gamma_{ab}\epsilon_-.
 \label{eq:spin-OmegaT-local}
\end{align}

\begin{remark}
\label{rem:spin-static-SU2-model}
Consider the model spinors $\epsilon_+^0$ and $\epsilon_-^0$ introduced in
Remark \ref{rem:spin-SU3-model}.  They satisfy
\begin{equation}
 \left\langle\epsilon_-^0,(\epsilon_+^0)^*\right\rangle=0,
\end{equation}
so that $((\epsilon_+^0)^*,\epsilon_-^0)$ is an ordered orthonormal pair in $\Delta_6^-$.  Let
$\eta_0,V_0,\omega_0^T$, and $\Omega_0^T$ be the tensors defined by
\eqref{eq:spin-eta-global}--\eqref{eq:spin-OmegaT-global} from this pair.  A
direct calculation in the model of Remark
\ref{rem:spin-SU3-model} gives
\begin{align*}
 \eta_0&=e^5-\ii e^6,
 &
 V_0&=e_5-\ii e_6,\\
 \omega_0^T&=e^{12}+e^{34},
 &
 \Omega_0^T&=(e^1+\ii e^2)\wedge(e^3+\ii e^4).
\end{align*}
\end{remark}

\begin{proposition}\label{prop:spin-static-SU2-normal-form}
The triple $(\eta,\omega^T,\Omega^T)$ defines an
$\operatorname{SU}(2)$-structure.  More precisely, every point has a
neighborhood carrying an oriented orthonormal coframe
$\{e^1,\ldots,e^6\}$ in which
\begin{align}
 \eta&=e^5-\ii e^6,
 \label{eq:spin-SU2-normal-form-eta}\\
 \omega^T&=e^{12}+e^{34},
 \label{eq:spin-SU2-normal-form-omega}\\
 \Omega^T&=(e^1+\ii e^2)\wedge(e^3+\ii e^4).
 \label{eq:spin-SU2-normal-form-Omega}
\end{align}
In particular, the following global algebraic identities hold:
\begin{align}
 g(V,V)&=0,
 &
 g(V,\overline V)&=2,
 \label{eq:spin-SU2-V-norms}\\
 \iota_V\omega^T&=0,
 &
 \iota_{\overline V}\omega^T&=0,
 \label{eq:spin-SU2-contraction-omega}\\
 \iota_V\Omega^T&=0,
 &
 \iota_{\overline V}\Omega^T&=0,
 \label{eq:spin-SU2-contraction-Omega}\\
 \omega^T\wedge\Omega^T&=0,
 &
 \Omega^T\wedge\Omega^T&=0,
 \label{eq:spin-SU2-wedge-zero}\\
 \omega^T\wedge\omega^T
 &=\frac12\Omega^T\wedge\overline{\Omega^T},
 \label{eq:spin-SU2-normalization}\\
 \operatorname{vol}_g
 &=\frac{\ii}{4}\overline\eta\wedge\eta
 \wedge\omega^T\wedge\omega^T.
 \label{eq:spin-SU2-volume}
\end{align}
The metric duality between $V$ and $\eta$ also gives
\begin{equation}\label{eq:spin-SU2-eta-evaluations}
 \eta(V)=\overline\eta(\overline V)=0,
 \qquad
 \eta(\overline V)=\overline\eta(V)=2.
\end{equation}
\end{proposition}

\begin{proof}
For any $p \in M$, since $\operatorname{SU}(4)$ acts transitively on ordered orthonormal two-frames in $S^-_p \cong \Delta_6^-$, there exists an element of $\operatorname{Spin}(6)$ that sends $(\epsilon_+^*,\epsilon_-)$ to
$((\epsilon_+^0)^*,\epsilon_-^0)$. By equivariance of the bilinear constructions,
there is an oriented orthonormal coframe of $T_p^*M$ in which 
the calculation in Remark
\ref{rem:spin-static-SU2-model} gives the stated forms at $p$.
Since $p$ is arbitrary, these normal forms hold pointwise on $M$.
The stabilizer of the model triple is $\operatorname{SU}(2)$;
hence the smooth triple $(\eta,\omega^T,\Omega^T)$ defines an
$\operatorname{SU}(2)$-reduction of the oriented orthonormal frame
bundle. The displayed algebraic identities follow by evaluating
them in an adapted coframe at each point.
\end{proof}

Write $V=V_1+\ii V_2$ with $V_1$ and $V_2$ real.  The two identities in
\eqref{eq:spin-SU2-V-norms} say precisely that
\begin{equation}
 |V_1|^2=|V_2|^2=1,
 \qquad
 g(V_1,V_2)=0.
\end{equation}
Hence $V_1$ and $V_2$ span an orthonormal two-plane, and we obtain the orthogonal splitting
\begin{equation}\label{eq:spin-SU2-tangent-splitting}
 TM=\mathcal H\oplus\mathcal V,
 \qquad
 \mathcal V=\operatorname{span}_{\mathbb R}\{V_1,
 V_2\},
 \qquad
 \mathcal H=\ker(\operatorname{Re}\eta)
 \cap\ker(\operatorname{Im}\eta).
\end{equation}
The contractions in the proposition show that $\omega^T$ and $\Omega^T$ are
horizontal.  Their restrictions to $\mathcal H$ give the usual
four-dimensional $\operatorname{SU}(2)$-structure.

The two $\operatorname{SU}(3)$-structures can be recovered from the transverse
forms by the global identities
\begin{align}
 \omega^\pm
 &=\omega^T\mp\frac{\ii}{2}\overline\eta\wedge\eta,
 \label{eq:spin-SU2-to-SU3-omega-global}\\
 \Omega^+
 &=-\eta\wedge\Omega^T,
 &
 \Omega^-
 &=\overline\eta\wedge\Omega^T.
 \label{eq:spin-SU2-to-SU3-Omega-global}
\end{align}
In the coframe of Proposition
\ref{prop:spin-static-SU2-normal-form}, they become
\begin{align}
 \omega^+&=e^{12}+e^{34}-e^{56},
 &
 \omega^-&=e^{12}+e^{34}+e^{56},
 \label{eq:spin-SU2-adapted-SU3-omega}\\
 \Omega^+
 &=-(e^5-\ii e^6)\wedge(e^1+\ii e^2)\wedge(e^3+\ii e^4),
 \label{eq:spin-SU2-adapted-Omega-plus}\\
 \Omega^-
 &=(e^5+\ii e^6)\wedge(e^1+\ii e^2)\wedge(e^3+\ii e^4).
 \label{eq:spin-SU2-adapted-Omega-minus}
\end{align}

We next compare the complex types determined by $J^+$ and $J^-$.  Reading the
two almost-complex structures from
\eqref{eq:spin-SU2-adapted-SU3-omega} gives
\begin{equation}\label{eq:spin-SU2-J-on-V}
 J^+V=\ii V,
 \qquad
 J^-V=-\ii V.
\end{equation}
On the splitting \eqref{eq:spin-SU2-tangent-splitting},
\begin{equation}\label{eq:spin-SU2-J-splitting}
 J^+|_{\mathcal H}=J^-|_{\mathcal H},
 \qquad
 J^+|_{\mathcal V}=-J^-|_{\mathcal V}.
\end{equation}
Thus $V$ is a $(1,0)$-vector for $J^+$ and a $(0,1)$-vector for $J^-$.  To
determine the type of its metric dual, recall that $g$ is extended complex
bilinearly.  For any complex vector field $W$, orthogonality and
\eqref{eq:spin-J-action-one-forms} give
\begin{equation}
 \bigl(J^\sigma(W^\flat)\bigr)(X)
 =-g(W,J^\sigma X)
 =g(J^\sigma W,X)
 =(J^\sigma W)^\flat(X).
\end{equation}
With $W=V$, \eqref{eq:spin-SU2-J-on-V} yields
\begin{equation}
 J^+\eta=\ii\eta,
 \qquad
 J^-\eta=-\ii\eta.
\end{equation}
Thus $\eta$ has type $(0,1)$ and $\overline\eta$ has type $(1,0)$ for
$J^+$; for $J^-$ the types are reversed.  This is also reflected in
$\eta(V)=g(V,V)=0$: the complex-bilinear metric vanishes on
$T^{1,0}\times T^{1,0}$ and pairs $T^{1,0}$ non-degenerately with
$T^{0,1}$.  Finally, \eqref{eq:spin-SU2-normal-form-Omega} shows that
$\Omega^T$ has transverse type $(0,2)$ for both structures.  The two
factorizations in \eqref{eq:spin-SU2-to-SU3-Omega-global} therefore have type
$(0,3)$, as expected.

\section{Supergravity equations}\label{sec:supergravity}

Below we recall the main features of ten-dimensional supergravity
relevant to our paper. While these are standard in the physics
literature, a mathematical audience may benefit from a brief account
of the physical origin of some familiar geometric structures.

\subsection{Heterotic and Type I supergravities}

Heterotic supergravity arises as the low-energy limit of Heterotic
string theory. As such, it is defined on a ten-dimensional spin 
manifold $M_{10}$ equipped with a Lorentzian metric $g^{(10)}$.
In addition, the theory comprises a real three-form $H$, a real
function $\phi$, called the \emph{dilaton}, a Yang--Mills connection
$A$ on a principal $G$-bundle $P\to M_{10}$ with curvature $F$,
and spinorial fields called the gravitino $\psi$, the dilatino
$\lambda$, and the gaugino $\chi$. These are sections
\begin{equation}
 \psi\in\Gamma(T^*M_{10}\otimes S_{10}^+),
 \qquad
 \lambda\in\Gamma(S_{10}^-),
 \qquad
 \chi\in\Gamma(S_{10}^+\otimes\operatorname{ad}P),
\end{equation}
where $S_{10}^\pm$ are the real chiral (Majorana--Weyl) spinor
bundles and $\operatorname{ad}P=P\times_{\operatorname{Ad}}\mathfrak g$
is the adjoint bundle. The gauge algebra $\mathfrak g$ is
$\mathfrak{so}(32)$ or $\mathfrak e_8\oplus\mathfrak e_8$.

We are interested in \emph{bosonic backgrounds}, obtained by setting
the spinorial fields to zero. The corresponding 
 action functional reduces to 
\begin{equation}\label{eq:action}
 S=\frac{1}{2\kappa_{10}^{\,2}}
 \int_{M_{10}}\ee^{-2\phi}
 \left(
 \operatorname{Scal}_{g^{(10)}}
 +4\langle \dd\phi,\dd\phi\rangle
 -\frac1{12}|H|^2
 -\alpha'\operatorname{tr}_{g^{(10)}}Q_F
 \right)\operatorname{vol}_{g^{(10)}}.
\end{equation}

Here $\kappa_{10}$ is the gravitational coupling and $\alpha'$ is
the string parameter. We also denote by $g^{(10)}$ the standard pairings induced by the Lorentzian metric on differential forms. In particular,
\begin{equation}\label{eq:sugra-contractions}
 |H|^2=g^{(10)}( H,H),
 \qquad
 H^2(X,Y)=g^{(10)}(\iota_XH,\iota_YH).
\end{equation}

Let $B$ denote the invariant inner product on the gauge algebra,
with the normalization of \cite{Gauntlett:2003cy}.
Choose a local pseudo-orthonormal frame $\{e_a\}_{a=1}^{10}$,
so that
\begin{equation}
 g^{(10)}(e_a,e_b)=s_a\delta_{ab},
 \qquad s_1=-1 \ \ \text{and} \ \ s_a=1, \ \ a=2,\dots,10.
\end{equation}
The symmetric tensor $Q_F$ appearing in the action is defined by
\begin{equation}
 Q_F(X,Y)
 =\sum_{a=1}^{10}s_a B\bigl(F(X,e_a),F(Y,e_a)\bigr).
\end{equation}
Its metric trace is
\begin{equation}
 \operatorname{tr}_{g^{(10)}}Q_F
 =\sum_{b=1}^{10}s_b Q_F(e_b,e_b).
\end{equation}
These definitions are independent of the chosen
pseudo-orthonormal frame. We retain the notation $\operatorname{Tr}(F\wedge F)$ for the
four-form obtained by contracting the gauge-algebra values of
$F\wedge F$ with $B$.

The Euler--Lagrange equations following from
\eqref{eq:action} can be expressed as
\begin{subequations}\label{eq:eom}
\begin{align}
 \operatorname{Rc}_{g^{(10)}}-\frac14H^2
 +2\operatorname{Hess}_{g^{(10)}}\phi-2\alpha'Q_F
 &=0,\label{eq:eom-a}\\
 \operatorname{tr}_{g^{(10)}}\operatorname{Hess}_{g^{(10)}}(\ee^{-2\phi})
 -\frac16 \ee^{-2\phi}|H|^2-\alpha'\ee^{-2\phi}\operatorname{tr}_{g^{(10)}}Q_F
 &=0,\label{eq:eom-b}\\
 \dd(\ee^{-2\phi}\star_{10}H)&=0,\label{eq:eom-c}\\
 \dd_A(\ee^{-2\phi}\star_{10}F)-\ee^{-2\phi}F\wedge\star_{10}H
 &=0.\label{eq:eom-d}
\end{align}
\end{subequations}
Here $\operatorname{Rc}$ and $\operatorname{Hess}$ are computed using
the Levi-Civita connection, $\star_{10}$ is the Lorentzian Hodge
operator, and $\dd_A$ is the exterior covariant derivative of $A$. Moreover, the three-form $H$ must satisfy the modified Bianchi identity
\begin{equation}
  \dd H = 2\alpha'
 \left(
 \operatorname{tr}(R^-\wedge R^-)
 -\operatorname{Tr}(F\wedge F)
 \right),
  \label{eq:bianchi-corrected}
\end{equation}
where  $R^-$  is the
curvature of $\nabla^{-H}$, the connection with opposite sign of the torsion.

Using the Clifford conventions fixed in the preliminaries, the
supersymmetry variations of the gravitino, dilatino, and gaugino
take the form
\begin{subequations}\label{eq:susy}
\begin{align}
 \delta\psi(X)
 &=\nabla_X\varepsilon-\frac14\,\iota_XH\cdot\varepsilon,
 \label{eq:susy-a}\\
 \delta\lambda
 &=\left(\dd\phi-\frac12H\right)\cdot\varepsilon,
 \label{eq:susy-b}\\
 \delta\chi&=2F\cdot\varepsilon.
 \label{eq:susy-c}
\end{align}
\end{subequations}
Here $\varepsilon$ is a Majorana--Weyl spinor, and the last equation
is understood componentwise in the gauge algebra. The gravitino
variation is the spinorial covariant derivative of the metric
connection with skew-symmetric torsion $H$, defined by
\begin{equation}\label{eq:sugra-torsion-connection}
 g^{(10)}(\nabla_X^H Y,Z)
 =g^{(10)}(\nabla_XY,Z)+\frac12H(X,Y,Z).
\end{equation}
With our positive Clifford convention, its lift to spinors is
\begin{equation}\label{eq:sugra-torsion-spin-connection}
 \nabla_X^H\varepsilon
 =\nabla_X\varepsilon-\frac14\,\iota_XH\cdot\varepsilon,
\end{equation}
so $\delta\psi=\nabla^H\varepsilon$. Our choice of $H$ is opposite
to that in the spinorial equations of \cite{Gauntlett:2003cy}; this
also accounts for the sign in \eqref{eq:bianchi-corrected}.
A bosonic background is supersymmetric if the variations
\eqref{eq:susy} vanish for some $\varepsilon\ne0$.

Setting $F = 0$ and $\alpha' = 0$ one obtains \emph{Type I supergravity}. The action functional, equations of motion, supersymmetry variations, and Bianchi identity for $H$, follow from those written above, upon setting  $F = 0$ and $\alpha'= 0$.

\subsection{Hull--Strominger system and conformally balanced BHE}

We now assume that the ten-dimensional manifold is a direct product\footnote{It can be proved that the supersymmetry equations imply that a warp factor $g^{(10)}=\ee^{2\Delta}(g_{\mathrm{Mink}}+g)$ must be $\Delta=$constant \cite{str}.}
\begin{equation}
   M_{10}=\mathbb R^{1,3}\times M_6,
 \qquad
 g^{(10)}=g_{\mathrm{Mink}}+g,  
\end{equation}
and that $\mathrm{Cliff}(1,9)$ decomposes as
\begin{equation}
\Gamma_A=\left\{
\begin{aligned}
        &\hat\gamma_\alpha\otimes \mathbf{1}_{8\times 8}\,, \qquad &&A =\alpha=0,\ldots, 3\,, \\
    &\hat\gamma_*\otimes \gamma_a\,, \qquad &&A=a+3\,,\, a=1,\ldots, 6\,.
\end{aligned}\right.
\end{equation}
We have introduced $\{\hat\gamma_\alpha\}_{\alpha=0,\ldots,3}$ as the generators of $\mathrm{Cliff(1,3)}$\footnote{For simplicity, we take the Majorana basis in $\mathrm{Cliff}(1,3)$, where the generators satisfy $\hat\gamma_\alpha^*=\hat\gamma_\alpha$.}, while $\{\gamma_a\}_{a=1,\ldots,6}$ are the generators of $\mathrm{Cliff(6)}$ defined in~\eqref{eq:spin-clifford-global} with the choice of basis in~\eqref{eq:spin-majorana-matrices}.
We also assume that the compactification preserves at least $\mathcal{N}=1$ supersymmetry, which gives the decomposition for the spinor $\varepsilon$ 
\begin{equation}
    \varepsilon=\zeta_+\otimes\epsilon+(\zeta_+)^c\otimes\epsilon^*\,
\end{equation}
where $\zeta_+$ is a constant Weyl spinor with positive chirality on $\mathbb{R}^{1,3}$, $(\zeta_+)^c$ is its conjugate, and $\epsilon$ is a non-zero complex Weyl spinor on $M_6$, with $\epsilon\in\Gamma(S^+)$.  
With $H$, $\phi$, and $F$ pulled back from the connected Riemannian
manifold $(M_6,g)$, the supersymmetry equations reduce to
\begin{subequations}\label{eq:KSE}
\begin{align}
 \nabla^H\epsilon&=0,\label{eq:KSE-2}\\
 \left(\dd\phi-\frac12H\right)\cdot\epsilon&=0,
 \label{eq:KSE-1}\\
 F\cdot\epsilon&=0.\label{eq:KSE-3}
\end{align}
\end{subequations}
The norm of $\epsilon$ is constant by \eqref{eq:KSE-2}, so we normalize it to one.
Let $(g,J,\omega,\Omega)$ be the associated
$\operatorname{SU}(3)$-structure. In the conventions of the
preliminaries,
\begin{equation}
 \omega(X,Y)=g(X,JY),
 \qquad
 \Omega\in\Omega_J^{0,3}(M_6),
 \qquad
 \operatorname{vol}_g=-\frac1{3!}\omega^3.
\end{equation}
We use the Hodge operator $\star$ determined by this orientation.

It was shown in \cite{str,Hull} that the gravitino and dilatino
equations \eqref{eq:KSE-2} and \eqref{eq:KSE-1} are equivalent to a system of equations on a complex
manifold with a Hermitian metric and a nowhere-vanishing holomorphic
three-form $\Psi$. In our conventions, this holomorphic form is
related to the normalized spinorial form by
\begin{equation}\label{eq:sugra-holomorphic-volume}
 \Psi=\ee^{-2\phi}\overline\Omega.
\end{equation}
Writing $\dd^{\mathrm{c}}=\ii(\bar\partial-\partial)$, the Hermitian form satisfies
\begin{subequations}\label{eq:strominger-original}
\begin{equation}\label{eq:so-dilatino}
 \dd^*\omega=-\dd^{\mathrm{c}}\log\|\Psi\|_g,
\end{equation}
while the three-form and dilaton are given by
\begin{equation}\label{eq:so-flux-dilaton}
 H=\dd^{\mathrm{c}}\omega,
 \qquad
 \phi=-\frac12\log\|\Psi\|_g+\phi_0,
\end{equation}
where $\phi_0$ is constant. Moreover, \eqref{eq:KSE-3} is equivalent
to the Hermitian--Yang--Mills condition
\begin{equation}\label{eq:so-hym}
 F\wedge\omega^2=0,
 \qquad F^{2,0}=F^{0,2}=0,
\end{equation}
and the Bianchi identity \eqref{eq:bianchi-corrected} becomes
\begin{equation}\label{eq:so-anomaly}
 \ii\partial\bar\partial\omega
 =\alpha'
 \left(
 \operatorname{tr}(R^-\wedge R^-)
 -\operatorname{Tr}(F\wedge F)
 \right).
\end{equation}
\end{subequations}

In \cite{Gauntlett:2003cy}, the same supersymmetry system was
expressed in terms of the normalized structure forms. Translated
into our conventions, it reads
\begin{subequations}\label{eq:gmw}
\begin{align}
 \dd(\ee^{-2\phi}\Omega)&=0,\label{eq:gmw-a}\\
 \dd(\ee^{-2\phi}\omega\wedge\omega)&=0,\label{eq:gmw-b}\\
 \star H&=\ee^{2\phi}\dd(\ee^{-2\phi}\omega),\label{eq:gmw-c}\\
 \star F&=\omega\wedge F.\label{eq:gmw-d}
\end{align}
\end{subequations}
In particular, the second equation shows that the metric is
conformally balanced. The signs in the last two equations use the orientation determined
by the positive-chirality spinor. This description in terms of
structure forms also has analogues in other dimensions
\cite{Gauntlett:2003cy}.

Setting $F=0$ and $\alpha'=0$ one obtains the Type I system, whose supersymmetry equations are given by (\ref{eq:KSE-2}) and (\ref{eq:KSE-1}), which are therefore equivalent 
to \eqref{eq:gmw-a}--\eqref{eq:gmw-c}. The Bianchi identity reduces
to $\dd H=0$, or equivalently to the pluriclosed condition $\ii\partial\bar\partial\omega=0.$
We make this specialization for the remainder of the section.

The equation of motion for $H$, given by \eqref{eq:eom-c}, is automatically satisfied using  \eqref{eq:gmw-c}. The equations of motion for the metric (\ref{eq:eom-a}) and the three-form (\ref{eq:eom-c}) can be written as 
\begin{equation}\label{eq:sugra-generalized-ricci-soliton}
 \operatorname{Rc}-\frac14H^2+2\operatorname{Hess}\phi=0,
 \qquad
 \dd^*H+2\iota_{\operatorname{grad}\phi}H=0.
\end{equation}
These are precisely the gradient generalized  Ricci soliton
equations, with potential $2\phi$. It was shown in
\cite[Appendix]{Gauntlett:2002sc} that the supersymmetry equations \eqref{eq:KSE} for $F=0$ and $\alpha'=0$,
together with the Bianchi identity and the equation of motion for
$H$, imply the metric and dilaton equations of motion.
Consequently, in the present setting, $\dd H=0$ and the spinorial
equations \eqref{eq:KSE-2} and \eqref{eq:KSE-1} imply
\eqref{eq:sugra-generalized-ricci-soliton}, without any compactness
assumption. 
For the convenience of the reader, we report the proof below, in the notation of this paper. In particular, we have:
\begin{theorem}\label{thm:grks-generalized-ricci-soliton}
Let $(M^6,g,\epsilon,H,\phi, \dd H=0)$ be a solution of the 
Killing spinor equations  \eqref{eq:KSE-2} and \eqref{eq:KSE-1}, with
$\epsilon\not\equiv0$. Then $(g,H,f)$, where $f=2\phi$, is a gradient generalized Ricci soliton. More precisely,
\begin{equation}\label{eq:steady-gradient-generalized-ricci-soliton}
 \operatorname{Rc}-\frac14H^2+\operatorname{Hess}f=0,
 \qquad
 \dd^*H+\iota_{\operatorname{grad}f}H=0.
\end{equation}
\end{theorem}

\begin{proof}
Since $\nabla^H$ is metric and preserves the orientation, it is compatible with Clifford multiplication. Thus  we have
\begin{equation}\label{eq:differentiated-dilatino}
\nabla_X^H\left(\left(\dd\phi-\frac12 H\right) \cdot \epsilon \right)= \left(\nabla_X^H \dd\phi-\frac12\nabla_X^HH\right)\cdot\epsilon=0.
\end{equation}

We take the following conventions:
\begin{equation}
 R^H(X,Y)Z=[\nabla_X^H,\nabla_Y^H]Z-\nabla_{[X,Y]}^HZ,
 \qquad
 \operatorname{Rc}^H(Y,Z)=\sum_a R^H(e_a,Y,Z,e_a), 
\end{equation}
and an analogous one for the Ricci tensor of the Levi-Civita connection, that we denote $ \operatorname{Rc}$. 

Let $\{e_a\}_{a=1}^6$ be a local oriented orthonormal frame of $TM$, with
dual coframe $\{e^a\}_{a=1}^6$. As in the spinorial preliminaries (Section \ref{sec:spinor-preliminaries}), we use the
fixed Majorana representation of the Clifford algebra. Thus, in the
corresponding local trivialization of $S$, Clifford multiplication by $e^a$ is
represented by the gamma matrices $\gamma^a$ introduced there.
 
The standard curvature identities for a metric connection with
skew torsion are proved in
\cite[Proposition~3.1, equations~(3.15), (3.19)--(3.21)]{IP}, \cite[Theorem~1.1]{AgrFe}, and
\cite[Propositions~3.18 and~3.20]{GFS}.
With the conventions above and $\dd H=0$, 
\begin{equation}\label{eq:contracted-bismut-curvature}
 R^H_{abcd}\gamma^{bcd}
 =\left(
 \frac13\nabla_aH_{bcd}
 -\frac12H_{ab}{}^eH_{ecd}
 \right)\gamma^{bcd}
 =\frac13(\nabla_a^HH)_{bcd}\gamma^{bcd}.
\end{equation}
Since a $3$-form acts on spinors with
the factor $1/3!$, equation \eqref{eq:contracted-bismut-curvature} gives
\begin{equation}
 -\frac12(\nabla_a^HH)\cdot\epsilon
 =-\frac14R^H_{abcd}\gamma^{bcd}\epsilon.
\end{equation}
Thus \eqref{eq:differentiated-dilatino} becomes
\begin{equation}\label{eq:dilatino-curvature}
 \left(
 (\nabla_a^H\dd\phi)_b\gamma^b
 -\frac14R^H_{abcd}\gamma^{bcd}
 \right)\epsilon=0.
\end{equation}

On the other hand, the curvature of the spin connection annihilates the
$\nabla^H$-parallel spinor. With our positive Clifford convention this is
\begin{equation}\label{eq:spinorial-integrability}
 0=-\frac14R^H_{abcd}\gamma^{cd}\epsilon.
\end{equation}
Multiplying by $\gamma^b$ and using the orthonormal-frame identity
\begin{equation}
 \gamma^b\gamma^{cd}
 =\gamma^{bcd}+\delta^{bc}\gamma^d-\delta^{bd}\gamma^c
\end{equation}
and the first--fourth Ricci contraction, we obtain
\begin{equation}\label{eq:contracted-spinorial-integrability}
 \left(
 R^H_{abcd}\gamma^{bcd}
 +2\operatorname{Rc}^H_{ad}\gamma^d
 \right)\epsilon=0.
\end{equation}
Substitution in \eqref{eq:dilatino-curvature} therefore yields
\begin{equation}\label{eq:M-clifford}
 \left(
 \operatorname{Rc}^H_{ab}
 +2(\nabla_a^H\dd\phi)_b
 \right)\gamma^b\epsilon=0.
\end{equation}

Set
\begin{equation}
 \mathcal{M}:=\operatorname{Rc}^H+2\nabla^H\dd\phi.
\end{equation}
Equation \eqref{eq:M-clifford} says intrinsically that the spinor-valued
$1$-form $X\mapsto \mathcal{M} (X,\mathord\cdot)\cdot\epsilon$ vanishes.  Taking the pointwise norm of
the preceding spinor-valued $1$-form gives
\begin{align*}
 0
 &=\sum_a\left|\mathcal{M}_{ab}\gamma^b\epsilon\right|^2\\
 &=(\mathcal{M}^T\mathcal{M})_{bc}\,\epsilon^\dagger\gamma^b\gamma^c\epsilon.
\end{align*}
Since $\mathcal{M}^T\mathcal{M}$ is symmetric, the Clifford relation implies
\begin{equation}
 (\mathcal{M}^T\mathcal{M})_{bc}\gamma^b\gamma^c
 =\frac12(\mathcal{M}^T\mathcal{M})_{bc}\{\gamma^b,\gamma^c\}
 =\operatorname{tr}(\mathcal{M}^T\mathcal{M})\operatorname{Id}.
\end{equation}
Consequently,
\begin{equation}
 0=\operatorname{tr}(\mathcal{M}^T\mathcal{M})|\epsilon|^2=|\mathcal{M}|^2|\epsilon|^2.
\end{equation}
The norm of a $\nabla^H$-parallel spinor is constant; hence the assumed
non-zero spinor is nowhere vanishing, and therefore
\begin{equation}\label{eq:bismut-soliton-tensor}
 \operatorname{Rc}^H+2\nabla^H\dd\phi=0.
\end{equation}

Finally, taking its symmetric and skew-symmetric parts and   using the identity \cite{FI}
\begin{equation}
  \operatorname{Rc}^H (X,Y) = \operatorname{Rc} (X,Y) - \frac{1}{2}\dd^* H(X,Y) - \frac{1}{4}\sum_{a} g\big(H(X,e_a), H(Y,e_a)\big) ,
\end{equation}
we get
\begin{equation}
 \operatorname{Rc}-\frac14H^2+2\operatorname{Hess}\phi=0,
 \qquad
 \dd^*H+2\iota_{\operatorname{grad}\phi}H=0.
\end{equation}
Upon setting $f=2\phi$, these are precisely the equations in Definition
\ref{dfn:grs}.
\end{proof}

With the convention $\Delta=-\dd^*\dd=\operatorname{tr}_g
\operatorname{Hess}$ on functions, the dilaton equation (\ref{eq:eom-b}) reads
\begin{equation}\label{eq:sugra-scalar-dilaton}
 \Delta(\ee^{-2\phi})-\frac16\ee^{-2\phi}|H|^2=0.
\end{equation}

We will rederive this equation in the notation of this paper in the next section in Proposition \ref{prop:scalar-dilaton-equation}, as a special case, by setting $\kappa =0$. This leads to a simple obstruction to the existence of compact solutions. In particular we have:

\begin{theorem} \label{nogotheorem} \cite{Gauntlett:2003cy} Let $(M^6,g,\epsilon,H,\phi, \dd H=0)$ be a solution of the 
Killing spinor equations  \eqref{eq:KSE-2} and \eqref{eq:KSE-1}, with
$\epsilon\not\equiv0$. If $M^6$ is compact without boundary then $H=0$, $\phi$ is constant and  $M^6$ is a Calabi--Yau manifold.
\end{theorem}

\begin{proof}
Since $H\wedge\star H=\tfrac16|H|^2\operatorname{vol}_g$, eq. (\ref{eq:sugra-scalar-dilaton}) can equivalently be written as
\begin{equation}\label{eq:dilaton-sixform}
 \dd\star \dd(\ee^{-2\phi})-\ee^{-2\phi}H\wedge\star H=0.
\end{equation}
If $M_6$ is compact without boundary, integrating this equation and
applying Stokes' theorem gives
\begin{equation}
 \int_{M_6} \ee^{-2\phi}H\wedge\star H=0,
\end{equation}
and hence $H=0$ \cite{Gauntlett:2003cy}.  The dilaton
equation (\ref{eq:sugra-scalar-dilaton}) then implies that $\phi$ is constant, and the remaining equations give a Ricci-flat K\"ahler metric.
\end{proof}

In view of the generalization that we will discuss in the next section, it is useful to recall the proof of some well-known properties. 

\begin{theorem}\label{wellknowtheorem}
Let $(M^6,g,\epsilon,H)$ be a solution to  $\nabla^H\epsilon=0$, with $\epsilon$ a chiral 
non-zero spinor. Then $\epsilon$ determines a 
$\nabla^H$-parallel $\operatorname{SU}(3)$-structure
$(g,J,\omega,\Omega)$. 
$\nabla^H$ is the
characteristic connection and
\begin{equation}
 N \in\Omega^3(M),
 \qquad
 H=\dd^{\mathrm{c}}\omega-N,
 \qquad
 \rho^H=0.
\end{equation}
Thus  the $\operatorname{SU}(3)$-structure defines an almost Calabi--Yau with torsion. 
\end{theorem}

\begin{proof}
Since $\nabla^H$ preserves the Hermitian inner product on spinors,
$\epsilon$ has constant norm. The proof of the general case is analogous, after suitably rescaling. We may therefore normalize it to
have unit norm.
Proposition \ref{prop:spin-Weyl-SU3}  gives the stated
$\operatorname{SU}(3)$-structure. Equations
\eqref{eq:spin-derivative-H-bilinear}--
\eqref{eq:spin-derivative-B-bilinear} and $\nabla^H\epsilon=0$ yield
\begin{equation}
 \nabla^H\omega=0,
 \qquad
 \nabla^H\Omega=0.
\end{equation}
Since $\nabla^Hg=0$, also $\nabla^HJ=0$. The assertion about the
Nijenhuis tensor and the torsion formula now follows from \cite[Theorem~10.1]{FI}.

It remains to consider the Ricci form. In a local orthonormal frame, \eqref{eq:spin-SU3-omega-local} gives
\begin{align*}
 \rho^H(X,Y)
 &=-\frac12R^H(X,Y,e_a,e_b) \omega^{ab}=\frac{\ii}{2}\left\langle\epsilon,
 R^H(X,Y,e_a,e_b)\gamma^{ab}\epsilon\right\rangle.
\end{align*}
The curvature of the spinorial lift is
\begin{equation}
 R^{S,H}_{X,Y}
 :=[\nabla_X^H,\nabla_Y^H]-\nabla_{[X,Y]}^H
 =-\frac14R^H(X,Y,e_a,e_b)\gamma^{ab}.
\end{equation}
Since $\nabla^H\epsilon=0$, one has
$R^{S,H}_{X,Y}\epsilon=0$, and therefore
\begin{equation}
 \rho^H(X,Y)
 =-2\ii\left\langle\epsilon,
 R^{S,H}_{X,Y}\epsilon\right\rangle=0.
\end{equation}
\end{proof}

Thus, the gravitino equation alone (without using the dilatino equation) gives a Bismut-parallel spinor 
and it implies that the Bismut--Ricci form vanishes.  Further employing the dilatino equation, one shows  $N=0$ so that the $J$ becomes integrable and 
$\dd H=0$ implies the pluriclosed condition. In conclusion, this yields
a Bismut--Hermitian--Einstein structure, which is also conformally
balanced by \eqref{eq:gmw-b}. For reasons that will become clear later, in this paper we shall sometimes denote these geometries as \BHEO~structures. However,  Theorem \ref{nogotheorem} shows that if $M^6$ is compact, the above  conformally
balanced BHE structure cannot have non-zero torsion. Below we
will propose a deformation of \eqref{eq:KSE-1} that leads to
other types of BHE structures, for which compact solutions are possible. In particular, we will recover the BHE structure  considered in
\cite{KS-26,ALS,ALL-26}, as well as further examples.

\section{Generalized Killing spinor equations}

Let $(M^6,g)$ be an oriented Riemannian spin manifold with complex spinor
bundle $S\to M$. Let $H\in\Omega^3(M)$ be a closed $3$-form,
$\phi\in C^\infty(M,\mathbb R)$, and $\kappa\in\mathbb C$.

We denote by $\nabla$ the Levi-Civita connection of $g$, as well as its lift
to the spinor bundle. The $3$-form $H$ defines a metric connection $\nabla^H$
on $TM$ by
\begin{equation}\label{eq:torsion-connection}
 g(\nabla_X^H Y,Z)
 =
 g(\nabla_XY,Z)
 +\frac12 H(X,Y,Z).
\end{equation}

With our positive Clifford convention (see Section \ref{sec:spinor-preliminaries}), the induced connection on the spinor
bundle is
\begin{equation}\label{eq:torsion-spin-connection}
 \nabla_X^H\epsilon
 =
 \nabla_X\epsilon
 -\frac14\,\iota_XH\cdot\epsilon,
 \qquad
 X\in\Gamma(TM),\quad \epsilon\in\Gamma(S).
\end{equation}
Here and below, we use the same notation $\nabla^H$ for the connection on
$TM$ and for its spinorial lift. Partly motivated by the relation with the Killing spinor equations of \emph{generalized Type I supergravity} (see Appendix \ref{gensugemb}), we will define the following spinorial system.

\begin{definition}
The \emph{generalized Killing spinor equations} (GKSE) associated with the data
$(g,H,\phi,\kappa, \dd H=0)$ are
\begin{align}
 \nabla^H\epsilon&=0,
 \label{eq:grav}\\
 \left(\dd\phi-\frac12H\right)\cdot\epsilon
 +\kappa\gamma_*\epsilon
 &=0. 
 \label{eq:dil}
 \end{align}
A spinor $\epsilon\in\Gamma(S)$ satisfying \eqref{eq:dil} and
\eqref{eq:grav} will be called a \emph{generalized Killing spinor}.
\end{definition}

Here $\gamma_*$ is the chirality operator introduced in the spinorial
preliminaries.  Equation \eqref{eq:grav} is the \emph{gravitino equation} and
states that $\epsilon$ is parallel with respect to the metric connection with
torsion $H$.  For simplicity, we will refer to equation \eqref{eq:dil} as the \emph{dilatino equation}, which provides an algebraic constraint.
When $\kappa=0$, the system reduces, up to the Clifford-sign
conventions adopted here, to the Killing spinor equations arising in the NS
sector of Heterotic supergravity (see equations \eqref{eq:KSE-1}, \eqref{eq:KSE-2}) \cite{FI, Gauntlett:2003cy,Hull,str}.

Let $\{e_a\}_{a=1}^6$ be a local oriented orthonormal frame of $TM$, with
dual coframe $\{e^a\}_{a=1}^6$. As in the spinorial preliminaries (Section \ref{sec:spinor-preliminaries}), we use the
fixed Majorana representation of the Clifford algebra. Thus, in the
corresponding local trivialization of $S$, Clifford multiplication by $e^a$ is
represented by the gamma matrices $\gamma^a$ introduced there.

In these local frames, equations \eqref{eq:dil} and \eqref{eq:grav} take the local
form
\begin{subequations}\label{eq:KSE-modified}
\begin{align}
 \left(
 \nabla_a
 -\frac18H_{abc}\gamma^{bc}
 \right)\epsilon
 &=0.
 \label{eq:KSE-modified-2}\\
 \left(
 \nabla_a\phi\,\gamma^a
 -\frac1{12}H_{abc}\gamma^{abc}
 +\kappa\gamma_*
 \right)\epsilon
 &=0 
 \label{eq:KSE-modified-1}
\end{align}
\end{subequations}

\subsection{Gradient generalized Ricci solitons} \label{section4.1}
It is simple to show that Equations \eqref{eq:grav} and \eqref{eq:dil} imply that $(g, H, \phi)$ is a gradient generalized Ricci soliton (see Definition \ref{dfn:grs}), as in the case of $\kappa=0$ of Theorem
\ref{thm:grks-generalized-ricci-soliton}. We will assume that $\epsilon$ is non-zero, and the computation will be performed in the local orthonormal frame $\{e_a\}$.

\begin{theorem}\label{thm:grks-generalized-ricci-soliton-kappa1}
Let $(M^6,g,\epsilon,H,\phi,\kappa, \dd H=0)$ be a solution of the generalized
Killing spinor equations \eqref{eq:dil}--\eqref{eq:grav}, with
$\epsilon\not\equiv0$. Then $(g,H,f)$, where $f=2\phi$, is a gradient
generalized Ricci soliton. More precisely,
\begin{equation}\label{eq:steady-gradient-generalized-ricci-soliton}
 \operatorname{Rc}-\frac14H^2+\operatorname{Hess}f=0,
 \qquad
 \dd^*H+\iota_{\operatorname{grad}f}H=0.
\end{equation}
\end{theorem}

\begin{proof}
Since $[\nabla^H,\gamma_*]=0$, the identity proved in \eqref{eq:differentiated-dilatino} is still valid. The rest of the proof of Theorem \ref{thm:grks-generalized-ricci-soliton} uses only the undeformed equation \eqref{eq:grav} and therefore it remains valid. 
\end{proof}

We conclude with the scalar equation satisfied by the dilaton. 
\begin{proposition}\label{prop:scalar-dilaton-equation}
Under the assumptions of Theorem
\ref{thm:grks-generalized-ricci-soliton-kappa1}, the dilaton satisfies
\begin{equation}\label{eq:scalar-dilaton-equation}
 2|\operatorname{grad}\phi|^2
 -\Delta\phi
 -\frac1{12}|H|^2
 +2\kappa^2
 =0.
\end{equation}
\end{proposition}

\begin{proof}
We start from the global identity
\eqref{eq:differentiated-dilatino}. Let
$\{e_a\}$ be the fixed local orthonormal frame, with dual coframe $\{e^a\}$. Taking
$X=e_a$, multiplying by $e^a$ on the left, and summing over $a$, we obtain
\begin{equation}\label{eq:clifford-trace-differentiated-dilatino}
 \left[\sum_a e^a\cdot
 \left(\nabla_{e_a}^H\dd\phi-\frac12\nabla_{e_a}^HH\right)
 \right]\cdot\epsilon=0.
\end{equation}
We compute the two Clifford traces in
\eqref{eq:clifford-trace-differentiated-dilatino} separately.

First consider the term containing $\dd\phi$. Since, for every form $\omega$,
\begin{equation}
 e^a\cdot\omega=e^a\wedge\omega+\iota_{e_a}\omega,
\end{equation}
we have
\begin{equation}\label{eq:split-trace-dphi}
 \sum_a e^a\cdot\nabla_{e_a}^H\dd\phi
 =\sum_a e^a\wedge\nabla_{e_a}^H\dd\phi
 +\sum_a\iota_{e_a}\nabla_{e_a}^H\dd\phi.
\end{equation}
The contraction term is the trace of $\nabla^H\dd\phi$. Indeed,
\begin{equation}
 \begin{aligned}
 \sum_a(\nabla_{e_a}^H\dd\phi)(e_a)
 &=\sum_a(\nabla_{e_a}\dd\phi)(e_a)
   -\frac12\sum_aH(e_a,e_a,\operatorname{grad}\phi)\\
 &=\Delta\phi.
 \end{aligned}
\end{equation}
To compute the exterior term, let $Y,Z$ be vector fields. We have
\begin{align*}
 \left(\sum_a e^a\wedge\nabla_{e_a}^H\dd\phi\right)(Y,Z)
 &=(\nabla_Y^H\dd\phi)(Z)-(\nabla_Z^H\dd\phi)(Y)\\
 &=\dd(\dd\phi)(Y,Z)-H(Y,Z,\operatorname{grad}\phi)\\
 &=-\bigl(\iota_{\operatorname{grad}\phi}H\bigr)(Y,Z).
\end{align*}
Therefore \eqref{eq:split-trace-dphi} becomes
\begin{equation}\label{eq:trace-nablaH-dphi}
 \sum_a e^a\cdot\nabla_{e_a}^H\dd\phi
 =\Delta\phi-\iota_{\operatorname{grad}\phi}H.
\end{equation}

We now compute the term containing $H$. Again,
\begin{equation}\label{eq:split-trace-H}
 \sum_a e^a\cdot\nabla_{e_a}^HH
 =\sum_a e^a\wedge\nabla_{e_a}^HH
 +\sum_a\iota_{e_a}\nabla_{e_a}^HH.
\end{equation}
It is a straightforward check using the expression of $\nabla^H$ that 
\begin{align*}
\sum_a\iota_{e_a}\nabla_{e_a}^HH=-\dd^*H.
\end{align*}

For the exterior part, introduce the canonical $4$-form
\begin{align*}
 \sigma_H(X,Y,Z,W)
 ={}&g( H(X,Y,\mathord\cdot),H(Z,W,\mathord\cdot))-g( H(X,Z,\mathord\cdot),H(Y,W,\mathord\cdot))+g(H(X,W,\mathord\cdot),H(Y,Z,\mathord\cdot)).
\end{align*}
Antisymmetrizing the formula for $\nabla^HH$ gives
\begin{equation}
 \left(\sum_a e^a\wedge\nabla_{e_a}^HH\right)_{bcde}
 =4\nabla^H_{[b}H_{cde]}
 =(\dd H)_{bcde}-2(\sigma_H)_{bcde}.
\end{equation}
Since $\dd H=0$, the exterior part in \eqref{eq:split-trace-H} is
$-2\sigma_H$. We have thus proved
\begin{equation}\label{eq:trace-nablaH-H}
 \sum_a e^a\cdot\nabla_{e_a}^HH
 =-\dd^*H-2\sigma_H.
\end{equation}

Substituting \eqref{eq:trace-nablaH-dphi} and
\eqref{eq:trace-nablaH-H} into
\eqref{eq:clifford-trace-differentiated-dilatino}, we obtain
\begin{equation}\label{eq:contracted-differentiated-dilatino}
 \left[
 \Delta\phi
 -\iota_{\operatorname{grad}\phi}H
 +\frac12\dd^*H
 +\sigma_H
 \right]\cdot\epsilon=0.
\end{equation}

It remains to eliminate the $4$-form $\sigma_H$. For this purpose, we square
the algebraic dilatino equation. The form $\dd\phi-\tfrac12H$ has odd degree,
and hence its Clifford action anticommutes with $\gamma_*$. Since
$\gamma_*^2=\operatorname{Id}_S$, equation \eqref{eq:dil} implies
\begin{equation}\label{eq:expanded-square-dilatino}
 \left[
 \dd\phi\cdot \dd\phi
 -\frac12\bigl(\dd\phi\cdot H+H\cdot \dd\phi\bigr)
 +\frac14H\cdot H
 +\kappa^2\operatorname{Id}_S
 \right]\epsilon=0.
\end{equation}
The three Clifford products in this expression are
\begin{align*}
 \dd\phi\cdot \dd\phi
 &=|\operatorname{grad}\phi|^2\operatorname{Id}_S,\\
 \dd\phi\cdot H+H\cdot \dd\phi
 &=2\iota_{\operatorname{grad}\phi}H,\\
 H\cdot H
 &=2\sigma_H-\frac16|H|^2\operatorname{Id}_S.
\end{align*}
The second identity is the positive-Clifford formula
$v\cdot\alpha+\alpha\cdot v=2\iota_v\alpha$ for a $3$-form $\alpha$.
In the third identity, the degree-six part vanishes because $H\wedge H=0$,
the degree-two part vanishes by symmetry, the degree-four part is
$2\sigma_H$, and the scalar part is
$-\tfrac16|H|^2\operatorname{Id}_S$.
Consequently, \eqref{eq:expanded-square-dilatino} is equivalent to
\begin{equation}\label{eq:squared-dilatino}
 \left[
 |\operatorname{grad}\phi|^2
 -\iota_{\operatorname{grad}\phi}H
 +\frac12\sigma_H
 -\frac1{24}|H|^2
 +\kappa^2
 \right]\cdot\epsilon=0.
\end{equation}

We now subtract twice \eqref{eq:squared-dilatino} from
\eqref{eq:contracted-differentiated-dilatino} and we therefore obtain
\begin{equation}\label{eq:scalar-dilaton-before-H-equation}
 \left[
 \left(
 \iota_{\operatorname{grad}\phi}H+\frac12\dd^*H
 \right)\!\cdot
 +\Delta\phi
 -2|\operatorname{grad}\phi|^2
 +\frac1{12}|H|^2
 -2\kappa^2
 \right]\epsilon=0.
\end{equation}
By Theorem \ref{thm:grks-generalized-ricci-soliton-kappa1}, $\dd^*H+2\iota_{\operatorname{grad}\phi}H=0,$
so the $2$-form term in
\eqref{eq:scalar-dilaton-before-H-equation} vanishes. Since $\epsilon$ is
nowhere zero, \eqref{eq:scalar-dilaton-before-H-equation} reduces to
\begin{equation}
 \Delta\phi
 -2|\operatorname{grad}\phi|^2
 +\frac1{12}|H|^2
 -2\kappa^2=0.
\end{equation}
\end{proof}
As previously remarked, equation \eqref{eq:scalar-dilaton-equation} reduces to the dilaton equation of motion in Type I supergravity (see Equation \eqref{eq:sugra-scalar-dilaton}) for $\kappa=0$. 
\begin{remark} \label{rmk:gotheorem}
The presence of the constant~$\kappa$ allows us to have compact solutions with non-trivial $H$ and~$\phi$, which are excluded when~$\kappa=0$.
Indeed, equation \eqref{eq:scalar-dilaton-equation} is equivalent to
\begin{align}
0
&=
\ee^{-2\phi}
\left(
2|\operatorname{grad}\phi|^2
-\Delta\phi
-\frac1{12}|H|^2
+2\kappa^2
\right)
\nonumber\\
&=
\frac12\Delta(\ee^{-2\phi})
-\frac1{12}\ee^{-2\phi}|H|^2
+2\ee^{-2\phi}\kappa^2.
\label{eq:weighted-scalar-dilaton}
\end{align}
Here we have used
\begin{equation}
\frac12\Delta(\ee^{-2\phi})
=
\ee^{-2\phi}
\left(
2|\operatorname{grad}\phi|^2-\Delta\phi
\right).
\end{equation}

Suppose now that $M$ is compact. Integrating
\eqref{eq:weighted-scalar-dilaton} over $M$ gives
\begin{equation}
0
=
\int_M
\left(
-\frac1{12}\ee^{-2\phi}|H|^2
+2\ee^{-2\phi}\kappa^2
\right)
\,\dd V,
\end{equation}
and hence
\begin{equation}\label{eq:integral-H-kappa}
\frac1{24}
\int_M \ee^{-2\phi}|H|^2\,\dd V
=
\int_M \ee^{-2\phi}\kappa^2\,\dd V.
\end{equation}

In particular, if $\kappa=0$, equation
\eqref{eq:integral-H-kappa} implies $H=0$ (and, in turn, $\phi=\text{const}$).  On the other hand, when $\kappa\neq0$, the integral identity no longer forces $H$ to vanish. Thus the modified Killing spinor equations avoid the \emph{no-go theorem} \cite{Gauntlett:2003cy}.
\end{remark}
\medskip

Since $H$, $\phi$ and the metric are real, by Equation \eqref{eq:scalar-dilaton-equation}, we conclude that either $\kappa \in \mathbb{R}$  or $\kappa \in \ii \mathbb{R}$. This is analogous to the standard argument regarding the ``Killing spinors'' on Einstein manifolds, where the two cases correspond to positive or negative curvature, respectively. We will return to the interpretation of these equations in the context of Bismut--Hermitian--Einstein geometry later. We stress out that in order to derive these equations we used only $\dd H=0$ and (\ref{eq:dil}),  (\ref{eq:grav}), and we did not assume that $M$ is compact.

\subsection{Almost CYT $\mathrm{SU}(3)$-structures} \label{sec:CYT}

By \eqref{eq:spin-chiral-splitting}, write
\begin{equation}
 \epsilon=\epsilon_++\epsilon_-,
 \qquad
 \epsilon_\pm\in\Gamma(S^\pm).
\end{equation}
The parallelism of $\gamma_*$ in
\eqref{eq:spin-star-chirality-parallel}, together with the parity identity
\eqref{eq:spin-form-chirality-parity}, decomposes the GKSE equations into
\begin{subequations}\label{eq:KSE-chiral-decomposition}
\begin{align}
 \nabla^H\epsilon_\pm&=0.
 \label{eq:KSE-chiral-gravitino}\\
 \left(\dd\phi-\frac12H\right)\cdot\epsilon_\mp
 \pm\kappa\epsilon_\pm
 &=0,
 \label{eq:KSE-chiral-dilatino}
\end{align}
\end{subequations}

If $\kappa\neq0$, \eqref{eq:KSE-chiral-dilatino} implies that
\(
 \epsilon_-=0 \Longleftrightarrow \epsilon_+=0.
\)
Thus a non-zero solution has two nowhere-vanishing chiral components. Each spinor  defines an $\SU(3)$ structure obeying the properties of Theorem \ref{wellknowtheorem}.

\begin{theorem}\label{thm:GRKS-induced-SU3-structures}
Let $(M^6,g,\epsilon,H,\phi,\kappa)$ be a solution of the generalized
Killing spinor equations \eqref{eq:dil}--\eqref{eq:grav}, with $\epsilon$ non-zero and $\kappa \neq 0$.  Then the chiral components of $\epsilon$ determine two distinct
$\nabla^H$-parallel $\operatorname{SU}(3)$-structures
$(g,J^\pm,\omega^\pm,\Omega^\pm)$. For each sign, $\nabla^H$ is the
characteristic connection and
\begin{equation}
 N_\pm\in\Omega^3(M),
 \qquad
 H=\dd_\pm^{\mathrm{c}}\omega^\pm-N_\pm,
 \qquad
 \rho_\pm^H=0.
\end{equation}
Thus both $\operatorname{SU}(3)$-structures are almost Calabi--Yau with
torsion (CYT). 
\end{theorem}

\begin{proof}
Since parallel, both $\epsilon_+$ and $\epsilon_-$ are nowhere vanishing.
Applying Theorem \ref{wellknowtheorem} to each chirality spinor concludes the proof.
\end{proof}

\begin{remark}
Theorem 5.4 holds verbatim by replacing the condition $\kappa \neq 0$ with the fact that $\epsilon_\pm$ are both non-vanishing.  The geometry associated to the equations with $\kappa=0$ and two generic (see below) chiral spinors $\epsilon_\pm$ was studied in \cite{Gauntlett:2003cy}, see Section 4.2, thus we will not discuss this case here.
\end{remark}

\section{$\SU(3)$ and $\SU(2)$ structure analysis} \label{sec:su2su3}
We now derive the differential constraints satisfied by the two
parallel $\operatorname{SU}(3)$-structures and by the associated
$\operatorname{SU}(2)$-structure. Throughout this section, we work
in the static case, namely
$\mathsf B(\epsilon_+,\epsilon_-)=0$.
The following proposition shows that any generic solution of
\eqref{eq:dil}--\eqref{eq:grav} with $\kappa\neq0$ can be replaced
by a static solution with the same data $(g,H,\phi,\kappa)$.
For simplicity, we will further assume that the two chiral
components have equal norms, which can then be normalized to one
by a common constant rescaling. The case of unequal norms is
discussed in Appendix~\ref{app:beyond-static}. The non generic case, that is $\epsilon_+=\alpha \epsilon_-^*$ is studied in Appendix \ref{sec:appendix}.

\begin{proposition}\label{prop:beyond-static}
Let $(M^6,g,\epsilon,H,\phi,\kappa)$ be a solution of
\eqref{eq:dil}--\eqref{eq:grav} with $\kappa\neq0$ and $\epsilon$
generic. The functions $|\epsilon_+|$, $|\epsilon_-|$ and
$\mathsf B(\epsilon_+,\epsilon_-)$ are constant on each connected
component of $M$. Moreover, there exists a generic static spinor
$\widetilde\epsilon$ such that
$(M^6,g,\widetilde\epsilon,H,\phi,\kappa)$ is a solution of the
same equations. Furthermore, if $|\epsilon_+|=|\epsilon_-|$ and
$\kappa\in\mathbb R$, then $\widetilde\epsilon$ can be
chosen with equal chiral norms.
If $\kappa\in\ii\mathbb R$ and
$|\epsilon_+|=|\epsilon_-|$, the static spinor constructed
below has equal chiral norms if and only if
$\mathsf B(\epsilon_+,\epsilon_-)=0$.
\end{proposition}
\begin{proof}
By equation \eqref{eq:KSE-chiral-gravitino}, $\epsilon_+$ and $\epsilon_-$
are $\nabla^H$-parallel. Furthermore, equation \eqref{eqn:metric_B} implies that  $|\epsilon_\pm|$ and
$\mathsf B(\epsilon_+,\epsilon_-)$ are constant.

Recall that by
Section~\ref{section4.1}, $\kappa\in\mathbb R\cup \ii\mathbb R$. Set $C\epsilon=\epsilon^*$ if
$\kappa\in\mathbb R$ and $C\epsilon=\gamma_*\epsilon^*$ if
$\kappa\in \ii\mathbb R$.
By \eqref{eq:spin-star-chirality-parallel} and by $[\nabla^H, \gamma_*]=0$, it follows that $C\epsilon$ satisfies \eqref{eq:grav}.
By \eqref{eq:spin-charge-conjugation},
$(\beta\cdot\psi)^*=(-1)^k\beta\cdot\psi^*$ for a real
$k$-form $\beta$, and
$(\gamma_*\psi)^*=-\gamma_*\psi^*$.
Conjugating \eqref{eq:dil} therefore gives
\begin{equation}
\left(\dd\phi-\frac12H\right)\cdot\epsilon^*
+\overline\kappa\gamma_*\epsilon^*=0.
\end{equation}
Since odd forms anticommute with $\gamma_*$ by
\eqref{eq:spin-form-chirality-parity}, it also follows that
\begin{equation}
\left(\dd\phi-\frac12H\right)\cdot\gamma_*\epsilon^*
-\overline\kappa\,\epsilon^*=0.
\end{equation}
Hence $\epsilon^*$ solves \eqref{eq:dil} with
$\overline\kappa$, and $\gamma_*\epsilon^*$ with
$-\overline\kappa$. In either case $C\epsilon$ is a solution
with the same $\kappa$, since $C\epsilon=\epsilon^*$ if
$\kappa\in\mathbb R$ and $C\epsilon=\gamma_*\epsilon^*$ if
$\kappa\in \ii\mathbb R$. Therefore, by complex linearity, 
$\epsilon_\nu:=\epsilon+\nu C\epsilon$ is a solution of the GKSE for every
$\nu\in\mathbb C$.

By \eqref{eq:spin-C-exchanges-chirality},
$(C\epsilon)_+=\epsilon_-^*$ and
$(C\epsilon)_-=\pm\epsilon_+^*$, with the upper sign in the
real case and the lower sign in the imaginary case.
Using
$\mathsf B(\epsilon_\pm,\epsilon_\pm^*)=|\epsilon_\pm|^2$
and
$\mathsf B(\epsilon_-^*,\epsilon_+^*)
=\overline{\mathsf B(\epsilon_+,\epsilon_-)}$, one finds
\begin{equation} \label{eqn:static}
\mathsf B\bigl((\epsilon_\nu)_+,(\epsilon_\nu)_-\bigr)
=\mathsf B(\epsilon_+,\epsilon_-)
\pm\nu\bigl(|\epsilon_+|^2\pm|\epsilon_-|^2\bigr)
\pm\nu^2\,\overline{\mathsf B(\epsilon_+,\epsilon_-)},
\end{equation}
with the same sign convention as above.
If $\mathsf B(\epsilon_+,\epsilon_-)=0$, take $\nu=0$;
otherwise choose a root of this quadratic polynomial.
Then for such a root $\mathsf B((\epsilon_\nu)_+,(\epsilon_\nu)_-)=0$.

It remains to check genericity. First, $\epsilon_\nu$ is non-zero. If not, $\epsilon$ would be proportional to $C\epsilon$, which contradicts that $\epsilon$ is generic. Since $\kappa\neq0$, the chiral dilatino equations imply
that both chiral components of $\epsilon_\nu$ are nowhere
vanishing. Moreover,
$\mathsf B((\epsilon_\nu)_+,(\epsilon_\nu)_-)=0$ means that
$(\epsilon_\nu)_+^*$ and $(\epsilon_\nu)_-$ are orthogonal.
They are therefore linearly independent, so $\epsilon_\nu$
is generic.

Suppose now that $|\epsilon_+|=|\epsilon_-|=a$. Since $
\langle\epsilon_+,\epsilon_-^*\rangle
=\langle\epsilon_-,\epsilon_+^*\rangle
=\overline{\mathsf B(\epsilon_+,\epsilon_-)},$
we obtain
\begin{align}
|(\epsilon_\nu)_+|^2
&=a^2(1+|\nu|^2)
 +2\operatorname{Re}\bigl(
 \nu\overline{\mathsf B(\epsilon_+,\epsilon_-)}\bigr),\\
|(\epsilon_\nu)_-|^2
&=a^2(1+|\nu|^2)
 \pm2\operatorname{Re}\bigl(
 \nu\overline{\mathsf B(\epsilon_+,\epsilon_-)}\bigr),
\end{align}
with the same sign convention as above.
Thus, in the real case the two chiral components have
equal norms for every $\nu$.

In the imaginary case, the static condition \eqref{eqn:static} becomes $
\mathsf B(\epsilon_+,\epsilon_-)
-\nu^2\overline{\mathsf B(\epsilon_+,\epsilon_-)}=0.$
If $\mathsf B(\epsilon_+,\epsilon_-)=0$, we take $\nu=0$,
so the norms remain equal. Otherwise, the two roots are $
\nu=\pm\frac{\mathsf B(\epsilon_+,\epsilon_-)}
{|\mathsf B(\epsilon_+,\epsilon_-)|},$
and hence
\begin{equation}
    |(\epsilon_\nu)_+|^2-|(\epsilon_\nu)_-|^2
=\pm4|\mathsf B(\epsilon_+,\epsilon_-)|\neq0.
\end{equation}
Therefore, in the imaginary case the static spinor
constructed above has equal chiral norms if and only if
$\mathsf B(\epsilon_+,\epsilon_-)=0$.
\end{proof}

\subsection{Integrability and Lee forms}
In this subsection we prove that the two almost complex structures $\Jpm$ introduced in Theorem \ref{thm:GRKS-induced-SU3-structures} are both integrable under the assumption that $\epsilon$ is static, that is $\langle \epsilon_+^*, \epsilon_- \rangle=0$, and that the two chiral spinors have unit norm. After establishing the integrability, we provide an explicit computation of the Lee forms of $(g, \Jpm)$ under the further assumption that $\epsilon_\pm$ have unit norms. \par

\medskip
Recall that the Nijenhuis tensor (as a $(0,3)$-tensor) has the following expression (see, for instance \cite[Theorem 10.1]{FI}):
\begin{equation} \label{eqn:integrability}
N(X,Y,Z)
=
-H(X,Y,Z)
+H(J^{\pm}X, J^{\pm}Y, Z)
+H(J^{\pm}X, Y, J^{\pm}Z)
+H(X, J^{\pm}Y, J^{\pm}Z).
\end{equation}

\begin{proposition}\label{prop:integrability-induced-complex-structures}
Let $(M^6,g,\epsilon,H,\phi,\kappa)$ be a solution of
\eqref{eq:dil}--\eqref{eq:grav}, with $\epsilon$ static and $\epsilon_\pm$ with unit norm. Then the almost
complex structures $J^\pm$ induced by the two chiral components of
$\epsilon$ are integrable. Consequently, the two induced
$\operatorname{SU}(3)$-structures are Bismut--Hermitian--Einstein.
\end{proposition}

\begin{proof}
By equation \eqref{eqn:integrability},
\begin{equation}
 N^\pm=0
 \quad\Longleftrightarrow\quad
 H^{3,0}_{J^\pm}=H^{0,3}_{J^\pm}=0.
\end{equation}
It therefore suffices to show that the pure-type component of $H$
vanishes.
Let $g_{\mathbb C}$ denote the complex-bilinear extension of the metric
contraction to complex-valued differential forms.

We first relate the contraction of $H$ with $\Omega^\pm$ to the
spinorial action of $H$. At a point, choose an orthonormal coframe and
write
\begin{equation}
H=\frac1{3!}H_{abc}e^{abc},
\qquad
\Omega^\pm
=
\frac1{3!}(\Omega^\pm)_{abc}e^{abc}.
\end{equation}
Since the metric coefficients in this coframe are Kronecker deltas,
Equation~\eqref{eq:spin-SU3-Omega-global} gives
\begin{equation}
(\Omega^\pm)_{abc}
=
-\mathsf B
\bigl(
\epsilon_\pm,\gamma_{abc}\epsilon_\pm
\bigr).
\end{equation}
Consequently,
\begin{equation}
\begin{aligned}
g_{\mathbb C}(H,\Omega^\pm)
=
H_{abc}(\Omega^\pm)^{abc}
=
-H_{abc}\,
\mathsf B
\bigl(
\epsilon_\pm,\gamma^{abc}\epsilon_\pm
\bigr)
=
-3!\,
\mathsf B
\bigl(
\epsilon_\pm,H\cdot\epsilon_\pm
\bigr)=
-6\,
\mathsf B
\bigl(
\epsilon_\pm,H\cdot\epsilon_\pm
\bigr),
\end{aligned}
\end{equation}
where we used the convention
\begin{equation}
H\cdot\epsilon_\pm
=
\frac1{3!}H_{abc}\gamma^{abc}\epsilon_\pm.
\end{equation}
In particular,
\begin{equation}
g_{\mathbb C}(H,\Omega^\pm)=0
\quad\Longleftrightarrow\quad
\mathsf B
\bigl(
\epsilon_\pm,H\cdot\epsilon_\pm
\bigr)=0.
\end{equation}

Since $\Omega^\pm$ has type $(0,3)$ and the metric pairing is extended complex-bilinearly,
\begin{equation}
  g_\C(H,\Omega^\pm)
 =g_\C(H^{3,0}_{J^\pm},\Omega^\pm).
\end{equation}
Moreover, $\Lambda^{3,0}_{J^\pm}$ is a complex line generated by
$\overline{\Omega^\pm}$, thus the contraction vanishes if and only if
$H^{3,0}_{J^\pm}=0$. Since $H$ is real, this is equivalent to
\begin{equation}
 H^{3,0}_{J^\pm}=H^{0,3}_{J^\pm}=0.
\end{equation}

It remains to prove that $\mathsf B
\bigl(
\epsilon_\pm,H\cdot\epsilon_\pm
\bigr)=0$. The two chiral components
of the dilatino equation give
\begin{equation}\label{eq:H-action-chiral-spinors}
 H\cdot\epsilon_\pm
 =
 2\,\dd\phi\cdot\epsilon_\pm
 \mp2\kappa\epsilon_\mp.
\end{equation}
Pairing this identity with $\epsilon_\pm$ using the complex-bilinear
form $\mathsf B$ yields
\begin{equation}
\mathsf B(\epsilon_\pm,H\cdot\epsilon_\pm)
=
2\mathsf B(\epsilon_\pm,\dd\phi\cdot\epsilon_\pm)
\mp2\kappa\mathsf B(\epsilon_\pm,\epsilon_\mp).
\end{equation}

The first term vanishes because $\mathsf B$ is symmetric and Clifford
multiplication by a one-form is skew with respect to $\mathsf B$ (see equation \eqref{eq:spin-B-symmetry-kform}):
\begin{equation}
\mathsf B(\epsilon_\pm,\dd\phi\cdot\epsilon_\pm)
=
\mathsf B(\dd\phi\cdot\epsilon_\pm,\epsilon_\pm)
=
-\mathsf B(\epsilon_\pm,\dd\phi\cdot\epsilon_\pm).
\end{equation}
The second term vanishes by the static condition. Indeed,
\begin{equation}
\mathsf B(\epsilon_+,\epsilon_-)
=
\langle\epsilon_+^*,\epsilon_-\rangle
=0,
\end{equation}
and the symmetry of $\mathsf B$ gives the same identity with the two
spinors interchanged. Consequently,
\begin{equation}
 \mathsf B(\epsilon_\pm,H\cdot\epsilon_\pm)=0.
\end{equation}
Therefore,
\begin{equation}
 H^{3,0}_{J^\pm}=H^{0,3}_{J^\pm}=0,
\end{equation}
and hence the type decomposition above gives $N^\pm=0$.

Thus both $J^+$ and $J^-$ are integrable. Using Theorem~\ref{thm:GRKS-induced-SU3-structures} the proof is completed.
\end{proof}

We compute the Lee forms associated to the two Hermitian structures $(g,\Jpm)$. Recall that the Lee form is defined as
\begin{equation}\label{Lee-form-definition}
    \leepm(X)=-\frac{1}{2}H(\Jpm X,e_i,\Jpm e_i). 
\end{equation}
With our conventions for $\omega^\pm$, this definition implies that we also have 
\begin{equation}\label{eq:Lambda-H-Lee}
 \mathrm{tr}_{\omega^\pm}H=2J^\pm\theta^\pm;
 \end{equation}
 where
\begin{equation}
 \bigl(\mathrm{tr}_{\omega^\pm}H\bigr)(X)
 :=H(X,e_a,e_b)(\omega^\pm)^{ab}.
\end{equation}
\begin{lemma}\label{lem:Lee-forms}
Let $(M^6,g,\epsilon,H,\phi,\kappa)$ be a solution of
\eqref{eq:dil}--\eqref{eq:grav}, with $\epsilon$ static and $\epsilon_\pm$ with unit norm. Let
$(g,J^\pm,\omega^\pm,\Omega^\pm)$ be the induced Hermitian
$\operatorname{SU}(3)$-structures and let $\theta^\pm$ denote their Lee
forms. Then
\begin{align}
 \theta^+
 &=2\dd\phi-\kappa\overline\eta-\overline\kappa\eta
 =2\dd\phi-2\operatorname{Re}(\kappa\overline\eta),
 \label{eq:Lee-plus}\\
 \theta^-
 &=2\dd\phi+\kappa\eta+\overline\kappa\,\overline\eta
 =2\dd\phi+2\operatorname{Re}(\kappa\eta).
 \label{eq:Lee-minus}
\end{align}
\end{lemma}

\begin{proof}
The chiral components of the dilatino equation read
\begin{equation}\label{eq:chiral-dilatino-Lee}
 \left(\dd\phi-\frac12H\right)\cdot\epsilon_+
 =\kappa\epsilon_-,
 \qquad
 \left(\dd\phi-\frac12H\right)\cdot\epsilon_-
 =-\kappa\epsilon_+.
\end{equation}

We first record a global identity that applies to either chiral component.
Let $X$ be a real vector field and denote its metric-dual one-form by
$X^\flat$. With our positive Clifford convention, Clifford multiplication
by a real one-form is Hermitian, whereas Clifford multiplication by a real
three-form is skew-Hermitian (see equation \eqref{eq:spin-adjoint-transpose-kform}). It follows that
\begin{align}
\left\langle
 X^\flat\cdot\epsilon_\pm,
 \left(\dd\phi-\frac12H\right)\cdot\epsilon_\pm
 \right\rangle
-\left\langle
 \left(\dd\phi-\frac12H\right)\cdot\epsilon_\pm,
 X^\flat\cdot\epsilon_\pm
 \right\rangle &=
\left\langle
 \epsilon_\pm,
 \left(
 [X^\flat\cdot,\dd\phi\cdot]
 -\frac12\{X^\flat\cdot,H\cdot\}
 \right)\epsilon_\pm
\right\rangle
\nonumber\\
&=
2\left\langle
 \epsilon_\pm,
 (X^\flat\wedge \dd\phi)\cdot\epsilon_\pm
\right\rangle
-\left\langle
 \epsilon_\pm,
 (\iota_XH)\cdot\epsilon_\pm
\right\rangle.
\label{eq:global-skew-dilatino-pairing}
\end{align}
Here we used the global Clifford identities
\begin{equation}
 [X^\flat\cdot,\dd\phi\cdot]
 =2(X^\flat\wedge \dd\phi)\cdot,
 \qquad
 \{X^\flat\cdot,H\cdot\}=2(\iota_XH)\cdot.
\end{equation}

The definition of $\omega^\pm$ implies that, for every real two-form
$\beta$,
\begin{equation}
 \left\langle\epsilon_\pm,
 \beta\cdot\epsilon_\pm\right\rangle
 =\frac{\ii}{2}g(\beta,\omega^\pm).
\end{equation}
Indeed, with respect to the orthonormal coframe $\{e^a\}$, $
 \beta=\frac12\beta_{ab}e^{ab}$, and since $(\omega^\pm)_{ab}=-\ii\left\langle
 \epsilon_\pm,\gamma_{ab}\epsilon_\pm\right\rangle,$ we have
\begin{equation}
 \begin{aligned}
 \left\langle
 \epsilon_\pm,\beta\cdot\epsilon_\pm
 \right\rangle
 &=
 \frac12\beta_{ab}
 \left\langle
 \epsilon_\pm,\gamma^{ab}\epsilon_\pm
 \right\rangle=
 \frac{\ii}{2}\beta_{ab}(\omega^\pm)^{ab}
 =\frac{\ii}{2}g(\beta,\omega^\pm).
 \end{aligned}
\end{equation}
Moreover,
\begin{equation}
 \begin{aligned}
 g(X^\flat\wedge \dd\phi,\omega^\pm)
 &=2\omega^\pm(X,\operatorname{grad}\phi),\\
 g(\iota_XH,\omega^\pm)
 &=\bigl(\mathrm{tr}_{\omega^\pm}H\bigr)(X)=2J^\pm\theta^\pm (X).
 \end{aligned}
\end{equation}
Replacing $X$ by
$J^\pm X$ in \eqref{eq:global-skew-dilatino-pairing} and using
\begin{equation}
 \omega^\pm(J^\pm X,\operatorname{grad}\phi)=\dd\phi(X),
 \qquad
 (J^\pm\theta^\pm)(J^\pm X)=\theta^\pm(X),
\end{equation}
the right-hand side becomes
\begin{equation}\label{eq:global-skew-dilatino-Lee}
 \ii\bigl(2\dd\phi-\theta^\pm\bigr)(X).
\end{equation}

We now apply this identity to the first equation in
\eqref{eq:chiral-dilatino-Lee}. The left hand side of
\eqref{eq:global-skew-dilatino-pairing} is
\begin{equation}
 \kappa\left\langle J^+X^\flat\cdot\epsilon_+,\epsilon_-\right\rangle
 -\overline\kappa
 \left\langle\epsilon_-,J^+X^\flat\cdot\epsilon_+\right\rangle
 =\kappa\overline\eta(J^+X)-\overline\kappa\eta(J^+X).
\end{equation}
By Proposition~\ref{prop:spin-static-SU2-normal-form}, $\overline\eta$ is
of type $(1,0)$ and $\eta$ is of type $(0,1)$ with respect to $J^+$.
Hence
\begin{equation}
 \overline\eta(J^+X)=\ii\overline\eta(X),
 \qquad
 \eta(J^+X)=-\ii\eta(X).
\end{equation}
Combining this with \eqref{eq:global-skew-dilatino-Lee} gives
\begin{equation} \label{eqn:theta_su3}
 \ii\bigl(2\dd\phi-\theta^+\bigr)(X)
 =\ii\bigl(\kappa\overline\eta+\overline\kappa\eta\bigr)(X),
\end{equation}
which proves \eqref{eq:Lee-plus}.

For the second equation in \eqref{eq:chiral-dilatino-Lee}, the corresponding
right-hand side is
\begin{equation}
 -\kappa\left\langle X^\flat\cdot\epsilon_-,\epsilon_+\right\rangle
 +\overline\kappa
 \left\langle\epsilon_+,X^\flat\cdot\epsilon_-\right\rangle
 =-\kappa\eta(X)+\overline\kappa\,\overline\eta(X).
\end{equation}
Now $\eta$ is of type $(1,0)$ and $\overline\eta$ is of type $(0,1)$ with
respect to $J^-$, and therefore
\begin{equation}
 \eta(J^-X)=\ii\eta(X),
 \qquad
 \overline\eta(J^-X)=-\ii\overline\eta(X).
\end{equation}
It follows that
\begin{equation}
 \ii\bigl(2\dd\phi-\theta^-\bigr)(X)
 =-\ii\bigl(\kappa\eta+\overline\kappa\,\overline\eta\bigr)(X),
\end{equation}
which proves \eqref{eq:Lee-minus}.
\end{proof}

\subsection{Properties of $V_1$ and $V_2$}
We determine the component of $\dd\phi$ along the parallel complex vector field
$V$. 

Recall that 
\begin{equation}
 \eta=\eta_1+\ii\eta_2,
 \qquad
 V=V_1+\ii V_2,
 \qquad
 V_i=\eta_i^\sharp.
\end{equation}
Since $\nabla^H$ is real and metric,
\begin{equation}\label{eq:static-SU2-parallel-vector-fields}
 \nabla^H\eta_i=0,
 \qquad
 \nabla^HV_i=0,
 \qquad i=1,2.
\end{equation}
In particular, the vector fields $V_i$ have constant norm. Moreover, 
\begin{equation}
 g(\nabla_XV_i,Y)
 =-\frac12H(X,V_i,Y)
 =\frac12H(V_i,X,Y).
\end{equation}
The right-hand side is skew-symmetric in $X$ and $Y$, and hence each $V_i$
is Killing. Since the Levi--Civita connection is torsion-free, the same
identity gives
\begin{align*}
 \dd\eta_i(X,Y)
 &=g(\nabla_XV_i,Y)-g(\nabla_YV_i,X)=H(V_i,X,Y).
\end{align*}
Consequently,
\begin{equation}\label{eq:static-SU2-deta}
 \dd\eta_i=\iota_{V_i}H,
 \qquad i=1,2.
\end{equation}
Finally, Cartan's formula and $\dd H=0$ imply
\begin{equation}\label{eq:static-SU2-H-invariant}
 \mathcal L_{V_i}H
 =\dd\iota_{V_i}H+\iota_{V_i}\dd H
 =\dd^2\eta_i=0.
\end{equation}
\begin{remark}
For comparison, on a compact complex manifold the implication from
the Bismut--Hermitian--Einstein condition to the gradient generalized Ricci soliton is proved in 
\cite[Proposition~2.6]{GFJS}. Furthermore, the vector fields $V_1$ and $V_2$ defined spinorially play the role of the canonical Bismut-parallel symmetries associated with a compact
non-K\"ahler BHE threefold (Section \ref{subsec:BHE-structures}); compare
\cite[Proposition~2.6 and Theorem~1.1]{ALS} (see also \cite{GFJS, SU}). Although similar in spirit, the vector fields $V_1$ and $V_2$ are not exactly the canonical vector fiels introduced in Section  \ref{subsec:BHE-structures}. Note for instance that $J^\pm V_1=\mp V_2$.
\end{remark}
\begin{lemma}\label{lem:dilaton-vertical-derivatives}
Let $(M^6,g,\epsilon,H,\phi,\kappa)$ be a solution of
\eqref{eq:dil}--\eqref{eq:grav}, with $\epsilon$ static and with
 $|\epsilon_+|=|\epsilon_-|=1$. Let $\eta$ and $V$ be defined by
\eqref{eq:spin-eta-global} and \eqref{eq:spin-V-global}, and write
\begin{equation}
 \eta=\eta_1+\ii\eta_2,
 \qquad
 V=V_1+\ii V_2,
 \qquad
 V_i=\eta_i^\sharp.
\end{equation}
Then
\begin{equation}\label{eq:dilaton-vertical-derivatives}
 \mathcal L_{V_1}\phi=0,
 \qquad
 \mathcal L_{V_2}\phi=\operatorname{Im}\kappa.
\end{equation}
\end{lemma}

\begin{proof}
The lower-sign equation in \eqref{eq:KSE-chiral-dilatino} is
\begin{equation}
 \left(\dd\phi-\frac12H\right)\cdot\epsilon_+
 =\kappa\epsilon_-.
\end{equation}
Taking its Hermitian product with $\epsilon_-$ and using
$|\epsilon_-|=1$ gives
\begin{equation}\label{eq:dilatino-pairing-minus}
 \kappa
 =V(\phi)
 -\frac12\left\langle
 \epsilon_-,H\cdot\epsilon_+
 \right\rangle,
\end{equation}
because
\begin{equation}
 \left\langle\epsilon_-,\dd\phi\cdot\epsilon_+\right\rangle
 =\eta(\operatorname{grad}\phi)
 =V(\phi).
\end{equation}

Similarly, the upper-sign equation gives
\begin{equation}
 \left\langle
 \epsilon_+,
 \left(\dd\phi-\frac12H\right)\cdot\epsilon_-
 \right\rangle
 =-\kappa.
\end{equation}
With our Clifford convention, multiplication by a real one-form is
self-adjoint, whereas multiplication by a real three-form is skew-adjoint:
$
 (\dd\phi\cdot)^\dagger=\dd\phi\cdot,
 (H\cdot)^\dagger=-H\cdot
$
(see \eqref{eq:spin-adjoint-transpose-kform}). Taking the complex conjugate of
the preceding identity therefore yields
\begin{equation}\label{eq:dilatino-pairing-plus}
 -\overline\kappa
 =V(\phi)
 +\frac12\left\langle
 \epsilon_-,H\cdot\epsilon_+
 \right\rangle.
\end{equation}
Adding \eqref{eq:dilatino-pairing-minus} and
\eqref{eq:dilatino-pairing-plus}, we obtain
\begin{equation}
 2V(\phi)
 =\kappa-\overline\kappa
 =2\ii\,\operatorname{Im}\kappa.
\end{equation}
Since $\phi$ is real and $V=V_1+\ii V_2$, this is equivalent to
\begin{equation}
 V_1(\phi)+\ii V_2(\phi)
 =\ii\,\operatorname{Im}\kappa,
\end{equation}
which proves \eqref{eq:dilaton-vertical-derivatives}.
\end{proof}

\medskip

We finally examine the infinitesimal action of the two parallel vector
fields on the spinor and on the induced structure forms.  

For any vector field $X$ the Kosmann derivative \cite{Kosmann} is 
\begin{equation}\label{eq:Kosmann-positive-convention}
 \mathcal L_X^K\psi
 =\nabla_X\psi+\frac14\dd X^\flat\cdot\psi.
\end{equation}

\begin{proposition}\label{prop:Kosmann-static-vector-fields}
Let $(M^6,g,\epsilon,H,\phi,\kappa)$ be a solution of
\eqref{eq:dil}--\eqref{eq:grav}, with $\epsilon$ static and with
 $|\epsilon_+|=|\epsilon_-|=1$. 
Then
\begin{equation}\label{eq:Kosmann-static-vector-fields}
 \mathcal L_{V_1}^K\epsilon=0,
 \qquad
 \mathcal L_{V_2}^K\epsilon
 =\ii\operatorname{Re}(\kappa)\epsilon.
\end{equation}
In particular, both vector fields preserve the spinor if and only if
$\operatorname{Re}(\kappa)=0$.
\end{proposition}

\begin{proof}
Equations \eqref{eq:static-SU2-deta} and \eqref{eq:grav} give
\begin{equation}
 \dd\eta_i=\iota_{V_i}H,
 \qquad
 \nabla_{V_i}\epsilon=\frac14(\iota_{V_i}H)\cdot\epsilon.
\end{equation}
Substitution in \eqref{eq:Kosmann-positive-convention} therefore yields the
global identity
\begin{equation}\label{eq:Kosmann-torsion-contraction}
 \mathcal L_{V_i}^K\epsilon
 =\frac12(\iota_{V_i}H)\cdot\epsilon.
\end{equation}
It remains to determine the two contractions of $H$ on the spinor.

The static normal form of Proposition
\ref{prop:spin-static-SU2-normal-form}, together with equivariance, gives
\begin{equation}\label{eq:static-vector-Clifford-action}
 V\cdot\epsilon_+=0,
 \qquad V\cdot\epsilon_-=2\epsilon_+,
 \qquad \overline V\cdot\epsilon_+=2\epsilon_-,
 \qquad \overline V\cdot\epsilon_-=0.
\end{equation}
On the other hand, the chiral dilatino equations are equivalently
\begin{equation}\label{eq:H-action-for-Kosmann}
 H\cdot\epsilon_+=2\dd\phi\cdot\epsilon_+-2\kappa\epsilon_-,
 \qquad
 H\cdot\epsilon_-=2\dd\phi\cdot\epsilon_-+2\kappa\epsilon_+.
\end{equation}
For every one-form $\alpha$, the positive Clifford relation gives
$\{\alpha\cdot,H\cdot\}=2(\iota_{\alpha^\sharp}H)\cdot$.  Using this
identity, \eqref{eq:static-vector-Clifford-action}, and
\eqref{eq:H-action-for-Kosmann}, we find
\begin{align*}
 (\iota_VH)\cdot\epsilon_+
 &=\frac12\{V\cdot,H\cdot\}\epsilon_+
 =V\cdot \dd\phi\cdot\epsilon_+-\kappa V\cdot\epsilon_-
 =2\bigl(V(\phi)-\kappa\bigr)\epsilon_+,\\
 (\iota_VH)\cdot\epsilon_-
 &=\frac12\{V\cdot,H\cdot\}\epsilon_-
 =V\cdot \dd\phi\cdot\epsilon_-+H\cdot\epsilon_+
 =2\bigl(V(\phi)-\kappa\bigr)\epsilon_-.
\end{align*}
The same calculation with $\overline V$ gives
\begin{equation}
 (\iota_{\overline V}H)\cdot\epsilon_\pm
 =2\bigl(\overline V(\phi)+\kappa\bigr)\epsilon_\pm.
\end{equation}
By Lemma~\ref{lem:dilaton-vertical-derivatives},
\begin{equation}
 V(\phi)=\ii\operatorname{Im}(\kappa),
 \qquad
 \overline V(\phi)=-\ii\operatorname{Im}(\kappa).
\end{equation}
Since $V=V_1+\ii V_2$, we conclude that
\begin{equation}\label{eq:vertical-H-action}
 (\iota_{V_1}H)\cdot\epsilon=0,
 \qquad
 (\iota_{V_2}H)\cdot\epsilon
 =2\ii\operatorname{Re}(\kappa)\epsilon.
\end{equation}
The result now follows from \eqref{eq:Kosmann-torsion-contraction}.
\end{proof}

The above proposition has the following

\begin{corollary}\label{cor:Lie-derivatives-structure-forms}
Under the hypotheses of Proposition
\ref{prop:Kosmann-static-vector-fields}, one has
\begin{equation}\label{eq:Lie-invariant-Hermitian-bilinears}
 \mathcal L_{V_i}\omega^\pm=0,
 \qquad
 \mathcal L_{V_i}\omega^T=0,
 \qquad
 \mathcal L_{V_i}\eta=0,
 \qquad i=1,2,
\end{equation}
whereas
\begin{align}
 \mathcal L_{V_1}\Omega^\pm&=0,
 &\mathcal L_{V_1}\Omega^T&=0,
 \label{eq:Lie-invariant-complex-bilinears-V1}\\
 \mathcal L_{V_2}\Omega^\pm
 &=2\ii\operatorname{Re}(\kappa)\Omega^\pm,
 &\mathcal L_{V_2}\Omega^T
 &=2\ii\operatorname{Re}(\kappa)\Omega^T.
 \label{eq:Lie-equivariant-complex-bilinears-V2}
\end{align}
In particular, both $V_1$ and $V_2$ are holomorphic with respect to
$J^+$ and $J^-$.  Here $\Omega^T$ is the globally defined two-form on
$M$ introduced in \eqref{eq:spin-OmegaT-global}.
\end{corollary}
\begin{proof}
We use that for a Killing vector field
\begin{align*}
 X\langle\psi,\chi\rangle
 &=\left\langle\mathcal L_X^K\psi,\chi\right\rangle
   +\left\langle\psi,\mathcal L_X^K\chi\right\rangle,\\
 X\mathsf B(\psi,\chi)
 &=\mathsf B\left(\mathcal L_X^K\psi,\chi\right)
   +\mathsf B\left(\psi,\mathcal L_X^K\chi\right).
\end{align*}
From the formulas above,
\begin{equation}\label{eq:Kosmann-Hermitian-bilinear-derivation}
 \mathcal L_X\bigl(\beta_k^{\mathrm H}(\psi,\chi)\bigr)
 =\beta_k^{\mathrm H}\bigl(\mathcal L_X^K\psi,\chi\bigr)
 +\beta_k^{\mathrm H}\bigl(\psi,\mathcal L_X^K\chi\bigr),
\end{equation}
and
\begin{equation}\label{eq:Kosmann-complex-bilinear-derivation}
 \mathcal L_X\bigl(\beta_k^{\mathrm B}(\psi,\chi)\bigr)
 =\beta_k^{\mathrm B}\bigl(\mathcal L_X^K\psi,\chi\bigr)
 +\beta_k^{\mathrm B}\bigl(\psi,\mathcal L_X^K\chi\bigr),
\end{equation}
are immediate. We will use such identities together with the fact that $V_1$ and $V_2$ are Killing.

For $X=V_1$, Proposition~\ref{prop:Kosmann-static-vector-fields} makes
both terms on the right-hand sides vanish.  We now discuss the case $X=V_2$. Since
$\mathcal L_{V_2}^K\epsilon_\pm=\ii\operatorname{Re}(\kappa)\epsilon_\pm$, conjugate-linearity in
the first argument of the Hermitian pairing gives
\begin{equation}
 \beta_k^{\mathrm H}(\ii
 \operatorname{Re}(\kappa)\psi,\chi)
 +\beta_k^{\mathrm H}(\psi,\ii \operatorname{Re}(\kappa)\chi)
 =(-\ii\operatorname{Re}(\kappa)+\ii\operatorname{Re}(\kappa))\beta_k^{\mathrm H}(\psi,\chi)=0.
\end{equation}
By contrast, $\mathsf B$ is complex bilinear, and therefore
\begin{equation}
 \beta_k^{\mathrm B}(\ii\operatorname{Re}(\kappa)\psi,\chi)
 +\beta_k^{\mathrm B}(\psi,\ii\operatorname{Re}(\kappa)\chi)
 =(\ii\operatorname{Re}(\kappa)+\ii\operatorname{Re}(\kappa))\beta_k^{\mathrm B}(\psi,\chi)
 =2\ii\operatorname{Re}(\kappa)\beta_k^{\mathrm B}(\psi,\chi).
\end{equation}
The forms $\omega^\pm$, $\eta$, and $\omega^T$ are Hermitian bilinears,
whereas $\Omega^\pm$ and $\Omega^T$ are complex-bilinear ones.  Equations
\eqref{eq:Lie-invariant-Hermitian-bilinears}--
\eqref{eq:Lie-equivariant-complex-bilinears-V2} now follow directly from
\eqref{eq:Kosmann-Hermitian-bilinear-derivation} and
\eqref{eq:Kosmann-complex-bilinear-derivation}.

Finally, $\mathcal L_{V_i}g=0$ because $V_i$ is Killing.  Since
$\omega^\pm(X,Y)=g(X,J^\pm Y)$ and both $g$ and $\omega^\pm$ are
preserved, $\mathcal L_{V_i}J^\pm=0$.
\end{proof}

\begin{corollary}\label{cor:static-local-transverse-geometry}
Under the hypotheses of Proposition \ref{prop:Kosmann-static-vector-fields}: 
\begin{enumerate}
    \item the vector fields $V_1$ and $V_2$ commute;
    \item $\dd\eta_i\in
\Lambda^{1,1}_{J^+}\cap\Lambda^{1,1}_{J^-}$; 
\item $\dd\eta_i, \omega^T$
and $H^T$ are basic; 
\item The torsion can be written as 
\begin{equation}\label{eq:H-transverse-decomposition}
 H=H^T+\eta_1\wedge \dd\eta_1+\eta_2\wedge \dd\eta_2.
\end{equation}
\end{enumerate}  
Moreover, locally there are coordinates $t_1,t_2$ along the two-dimensional
foliation such that
\begin{equation}\label{eq:local-static-fibration-forms}
 V_i=\partial_{t_i},
 \qquad
 \eta_i=\dd t_i+\sigma_i,
 \qquad i=1,2,
\end{equation}
where $\sigma_1$ and $\sigma_2$ are basic one-forms.  The complex form $\Omega^T$ is basic when
$\operatorname{Re}(\kappa)=0$; in general it is horizontal and satisfies
\eqref{eq:Lie-equivariant-complex-bilinears-V2}.  On each such coordinate
neighborhood, the locally defined phase-corrected form
\begin{equation}\label{eq:phase-corrected-transverse-volume}
 \widehat\Omega^T
 :=\ee^{-2\ii\operatorname{Re}(\kappa)t_2}\Omega^T
\end{equation}
is basic. 
\end{corollary}

\begin{proof}
The fact that $[V_1,V_2]=0$ follows since $V_1$ is $J^\pm$-holomorphic and $V_2=\mp J^\pm V_1$. Since $V_1$ and $V_2$ are Killing, this implies also that 
\begin{equation} \label{{eq:deta-horizontal}}
    \iota_{V_{i}} \dd\eta_j=0.
\end{equation}

Let
\begin{equation}
\mathcal V:=\operatorname{span}\{V_1,V_2\},
\qquad
\mathcal H:=\mathcal V^\perp
=\ker\eta_1\cap\ker\eta_2.
\end{equation}
Since \(V_1\) and \(V_2\) are Killing, the horizontal distribution \(\mathcal H\) is invariant under the
flows of \(V_1\) and \(V_2\). Consequently,
\begin{equation}
g^T
=
g-\eta_1\otimes\eta_1-\eta_2\otimes\eta_2
\end{equation}
is a basic transverse metric. Moreover, since
\(\mathcal L_{V_i}J^\pm=0\) and \(J^\pm\) preserve both
\(\mathcal V\) and \(\mathcal H\), their common restriction
\begin{equation}
J^T:=J^+|_{\mathcal H}=J^-|_{\mathcal H}
\end{equation}
is also basic. It follows that the transverse fundamental form
\begin{equation}
\omega^T(\mathord\cdot,\mathord\cdot)
=
g^T(\mathord\cdot,J^T\mathord\cdot)
\end{equation}
is basic as well. The horizontal component \(H^T\) of \(H\) is precisely the
Bismut torsion of the transverse Hermitian structure:
\begin{equation}\label{eq:transverse-horizontal-torsion-Lee}
H^T=\dd_T^{\mathrm{c}}\omega^T=-\star_T\theta^T,    
\end{equation}
where $\star_T$ is the Hodge star of $g^T$ with respect to the orientation $\tfrac12 \omega^T \wedge \omega^T$.
Since exterior differentiation and the action of the basic tensor \(J^T\)
preserve basic forms, this identity shows directly that \(H^T\) is basic.

We next determine the type of the curvature forms \(\dd\eta_i\). Since $V_i$ are $J^\pm$-holomorphic, 
\begin{equation}
\begin{aligned}
g((\mathcal L_{V_{i}}J^\pm )(X), Y)
&=g([V_i,J^\pm X]-J^\pm [V_i,X], Y) \\
&=-H(V_i,J^\pm X, Y)-H(V_i,X,J^\pm  Y)\\
&= - \dd\eta_i (J^\pm X, Y)-\dd\eta_i (X, J^\pm Y), 
\end{aligned}
\end{equation}
where we have used \eqref{eq:static-SU2-deta} for the last equality. This proves that $
\dd\eta_i\in
\Lambda^{1,1}_{J^+}\cap\Lambda^{1,1}_{J^-}
$. From this, equation \eqref{eq:H-transverse-decomposition} easily follows. 

The vector fields $V_1$ and $V_2$ are everywhere linearly independent and commute, so the
local coordinate description in \eqref{eq:local-static-fibration-forms}
follows from the Frobenius theorem. Finally,
\begin{align*}
 \mathcal L_{V_1}\widehat\Omega^T&=0,\\
 \mathcal L_{V_2}\widehat\Omega^T
 &=\ee^{-2\ii\operatorname{Re}(\kappa)t_2}
 \left(-2\ii\operatorname{Re}(\kappa)\Omega^T
 +\mathcal L_{V_2}\Omega^T\right)=0,
\end{align*}
and $\widehat\Omega^T$ is horizontal because $\Omega^T$ is.  This proves
that it is basic on the chosen coordinate neighborhood.
\end{proof}

We now retain only the differential identities needed for the real and
purely imaginary reductions.  Write
\begin{equation}
 a:=\operatorname{Re}(\kappa),
 \qquad
 b:=\operatorname{Im}(\kappa).
\end{equation}
Since \(J^\pm\) are integrable and \(\nabla^H\) is their Bismut
connection, the Bismut torsion identities, with our orientation and
fundamental-form conventions, are
\begin{equation} \label{eqn:symplectic}
    \dd\omega^+=\theta^+\wedge\omega^++\star H,
 \qquad
 \dd\omega^-=\theta^-\wedge\omega^--\star H. 
\end{equation}
Using Lemma~\ref{lem:Lee-forms},
\(\omega^\pm=\omega^T\pm\eta_1\wedge\eta_2\), and
\begin{equation}
 \operatorname{Re}(\kappa\overline\eta)=a\eta_1+b\eta_2,
 \qquad
 \operatorname{Re}(\kappa\eta)=a\eta_1-b\eta_2,
\end{equation}
we compute $  \dd\bigl(\ee^{-2\phi}\omega^T\bigr)=\tfrac12 \dd\bigl(\ee^{-2\phi}(\omega^++\omega^-)\bigr)$ and obtain
\begin{equation}\label{eq:weighted-transverse-fundamental-form}
 \dd\bigl(\ee^{-2\phi}\omega^T\bigr)
 =-2b\,\eta_2\wedge \ee^{-2\phi}\omega^T.
\end{equation}
Expanding \eqref{eq:weighted-transverse-fundamental-form} yields
\begin{equation}\label{eq:transverse-Lee-form}
 \dd\omega^T=\theta^T\wedge\omega^T,
 \qquad
 \theta^T:=2\dd\phi-2b\,\eta_2
 =\frac{\theta^++\theta^-}{2}.
\end{equation}
This is the Lee form of the transverse Hermitian structure. Furthermore, since the transverse Hermitian structure is invariant under \(V_1,V_2\); hence $\theta^T$ is basic and descends to every local leaf space.

The complex transverse form \(\Omega^T\) itself is global on \(M\), by
its spinorial definition \eqref{eq:spin-OmegaT-global}.  Since
\(\nabla^H\Omega^\pm=0\) and \(J^\pm\) are integrable, one has
\(\dd\Omega^\pm=\theta^\pm\wedge\Omega^\pm\) (see \cite[Theorem 34]{II}).  
Using \(\dd\Omega^\pm=\theta^\pm\wedge\Omega^\pm\) and the expressions for
the Lee forms (see Lemma \ref{lem:Lee-forms}), we obtain
\begin{subequations}\label{eq:weighted-holomorphic-forms}
\begin{align}
\dd\bigl(\ee^{-2\phi}\Omega^+\bigr)
&=
-(\kappa\bar\eta+\bar\kappa\eta)
\wedge \ee^{-2\phi}\Omega^+,
\label{eq:weighted-holomorphic-plus}\\
\dd\bigl(\ee^{-2\phi}\Omega^-\bigr)
&=
(\kappa\eta+\bar\kappa\bar\eta)
\wedge \ee^{-2\phi}\Omega^-.
\label{eq:weighted-holomorphic-minus}
\end{align}
\end{subequations}
Recall that by equation \eqref{eq:spin-SU2-to-SU3-Omega-global},
\begin{equation}
\Omega^+=-\eta\wedge\Omega^T,
\qquad
\Omega^-=\bar\eta\wedge\Omega^T.
\end{equation}
Moreover, \(\dd\eta\) and \(\dd\bar\eta\) are transverse forms of type
\((1,1)\). Since the transverse complex dimension is two,
\begin{equation}
\dd\eta\wedge\Omega^T
=
\dd\bar\eta\wedge\Omega^T
=
0.
\end{equation}
Substituting the preceding expressions for \(\Omega^\pm\) into
\eqref{eq:weighted-holomorphic-forms} therefore gives
\begin{subequations}\label{eq:transverse-holomorphic-wedge}
\begin{align}
\eta\wedge \dd\bigl(\ee^{-2\phi}\Omega^T\bigr)
&=
-\kappa\,\eta\wedge\bar\eta\wedge \ee^{-2\phi}\Omega^T,
\label{eq:transverse-holomorphic-eta}\\
\bar\eta\wedge \dd\bigl(\ee^{-2\phi}\Omega^T\bigr)
&=
-\kappa\,\eta\wedge\bar\eta\wedge \ee^{-2\phi}\Omega^T.
\label{eq:transverse-holomorphic-bareta}
\end{align}
\end{subequations}
The right-hand sides of these two identities are obtained by wedging
\begin{equation}
\kappa(\eta-\bar\eta)\wedge \ee^{-2\phi}\Omega^T
\end{equation}
with \(\eta\) and \(\bar\eta\), respectively. 
Consequently,
\begin{equation}
\dd\bigl(\ee^{-2\phi}\Omega^T\bigr)
=
\kappa(\eta-\bar\eta)\wedge \ee^{-2\phi}\Omega^T.
\end{equation}
Observe that $\dd\bigl(\ee^{-2\phi}\Omega^T\bigr)$ has no component along $\eta \wedge \bar \eta$, as it is an horizontal form and $[V,\bar V]=0$.
Since \(\eta-\bar\eta=2\ii\eta_2\), we conclude that
\begin{equation}\label{eq:transverse-holomorphic-form}
\dd\bigl(\ee^{-2\phi}\Omega^T\bigr)
=
2\ii\kappa\,\eta_2\wedge \ee^{-2\phi}\Omega^T.
\end{equation}

We next derive the relation
between the six- and four-dimensional Lee forms of the Hermitian structures $(J^+, \omega^+)$ and $(J^T,\omega^T)$.  Since
\begin{equation}
 \omega^+=\omega^T+\eta_1\wedge\eta_2,
 \qquad
 \dd\omega^T=\theta^T\wedge\omega^T,
\end{equation}
we have
\begin{equation}
 \dd\omega^+
 =\theta^T\wedge\omega^T
  +\dd\eta_1\wedge\eta_2-\eta_1\wedge \dd\eta_2.
\end{equation}

Using $\theta^\pm= \tfrac12 \operatorname{tr}_{\omega^\pm} \dd\omega^\pm$, the above equation gives
\begin{equation}\label{eq:Lee-trace-decomposition}
 \theta^+
 =\theta^T
 +\frac12\operatorname{tr}_{\omega^T}(\dd \eta_1)\,\eta_2
 -\frac12\operatorname{tr}_{\omega^T}( \dd \eta_2)\,\eta_1.
\end{equation}

Combining \eqref{eq:transverse-horizontal-torsion-Lee} with the
decomposition in Corollary~\ref{cor:static-local-transverse-geometry} and the formula above,
we obtain
\begin{equation} \label{eq:H-from-transverse-structure}
 H=-\star_T\theta^T+\eta_1\wedge \dd \eta_1+\eta_2\wedge \dd \eta_2,
\end{equation}
and the Bianchi identity \(\dd H=0\) is exactly
\begin{equation}\label{eq:universal-transverse-Bianchi}
 -\dd\star_T\theta^T+ \dd \eta_1 \wedge \dd \eta_1+\dd \eta_2 \wedge \dd \eta_2=0.
\end{equation}
This is the common starting point for the two reductions below. From now on we will assume that $\kappa \neq 0$.

\subsection{Real \(\kappa\)}

Below we assume that \(\kappa\in\mathbb R\setminus\{0\}\). We will refer to the corresponding geometries as \BHER structures. We will omit pullbacks: same letters denote both the forms on the base and the total space.

\begin{theorem}\label{thm:real-kappa-local-model}
Let \((M^6,g,\epsilon,H,\phi,\kappa)\) be a solution
of the GKSE with $\epsilon$ static and $\epsilon_\pm$ with unit norm. Assume
\(\kappa\in\mathbb R\setminus\{0\}\).  Every point has a neighborhood
which is a local principal \(\mathbb R^2\)-bundle
\begin{equation}
 \pi:U\longrightarrow B^4
\end{equation}
whose vertical fields are \(V_1,V_2\) and whose connection one-forms are
\(\eta_1,\eta_2\).  The local leaf space \(B\) carries the Kähler metric
\begin{equation}
 g_K^T=\ee^{-2\phi}g^T,
 \qquad
 \omega_K=\ee^{-2\phi}\omega^T.
\end{equation}
Let \(\rho_K\), \(R_K=\mathrm{tr}_{\omega_K}\rho_K\), and
\(\operatorname{Rc}_K\) denote its Ricci form, scalar curvature, and
Ricci tensor, with the conventions fixed above.  Then
\begin{equation}
 \dd\eta_1\in\Lambda^{1,1}_0(B),
 \qquad
 \dd\eta_2=\frac1{2\kappa}\rho_K,
 \qquad
 R_K=8\kappa^2\ee^{2\phi}>0,
\end{equation}
where $\Lambda^{1,1}_0(B)$ denotes the space of ASD forms on $B$.  The Bianchi identity is equivalent to
the scalar equation
\begin{equation}\label{eq:real-kappa-scalar-master-equation}
-\Delta_KR_K+\frac12R_K^2
 -|\operatorname{Rc}_K|^2
 -4\kappa^2| \dd \eta_1|^2=0,
 \qquad
 \Delta_K=-\dd_K^*\dd.
\end{equation}

Conversely, let \((B^4,g_K,J^T,\omega_K)\) be a Kähler surface with
\(R_K>0\), and let \(\alpha\) be a closed real primitive
\((1,1)\)-form satisfying
\eqref{eq:real-kappa-scalar-master-equation}.  On a sufficiently small
open set choose connection one-forms on
\(B\times\mathbb R^2\) such that
\begin{equation}
 \dd\eta_1=\alpha,
 \qquad
 \dd\eta_2=\frac1{2\kappa}\rho_K.
\end{equation}
Then
\begin{equation}
 g=\frac{R_K}{8\kappa^2}g_K+\eta_1^2+\eta_2^2,
 \qquad
 \omega=\frac{R_K}{8\kappa^2}\omega_K+\eta_1\wedge\eta_2
\end{equation}
define a Hermitian structure whose Bismut torsion is
\begin{equation} \label{eq:H-real}
 H=\frac1{8\kappa^2}\dd^{\mathrm{c}}R_K\wedge\omega_K
   +\eta_1\wedge\alpha
   +\frac1{2\kappa}\eta_2\wedge\rho_K.
\end{equation}
This Hermitian structure is pluriclosed and has vanishing Bismut Ricci
form; in particular, it is Bismut--Hermitian--Einstein.
\end{theorem}

\begin{proof}
Since \(\operatorname{Im}(\kappa)=0\), Lemma
\ref{lem:dilaton-vertical-derivatives} shows that \(\phi\) is basic.
Equation~\eqref{eq:weighted-transverse-fundamental-form} becomes
\(\dd(\ee^{-2\phi}\omega^T)=0\), proving that
\((g_K^T,J^T,\omega_K)\) is Kähler on every local leaf space.

In the local coordinates
\(t_1,t_2\) of Corollary
\ref{cor:static-local-transverse-geometry}, write
\(\eta_2=\dd t_2+\sigma_2\) and set
\(\widehat\Omega^T=\ee^{-2\ii\kappa t_2}\Omega^T\).
Equation~\eqref{eq:transverse-holomorphic-form} gives
\begin{equation}
 \begin{aligned}
 \dd\bigl(\ee^{-2\phi}\widehat\Omega^T\bigr)
 &=\ee^{-2\ii\kappa t_2}
   \left[
   \dd\bigl(\ee^{-2\phi}\Omega^T\bigr)
   -2\ii\kappa\,\dd t_2\wedge \ee^{-2\phi}\Omega^T
   \right]\\
 &=2\ii\kappa\,\sigma_2\wedge
   \ee^{-2\phi}\widehat\Omega^T.
 \end{aligned}
\end{equation}
Thus \(2\kappa\sigma_2\) is the local connection one-form of the
Chern connection.  With our
definition of the Ricci form,
\begin{equation}
 \rho_K=2\kappa\,\dd\sigma_2=2\kappa\,\dd\eta_2.
\end{equation}

For real \(\kappa\), Lemma~\ref{lem:Lee-forms} and
\eqref{eq:transverse-Lee-form} give
\begin{equation}
 \theta^+=2\dd\phi-2\kappa\eta_1,
 \qquad
 \theta^T=2\dd\phi.
\end{equation}
Comparing with \eqref{eq:Lee-trace-decomposition}, separately in the
\(\eta_1\)- and \(\eta_2\)-directions, yields
\begin{equation}
 \operatorname{tr}_{\omega^T}(\dd \eta_1)=0,
 \qquad
 \operatorname{tr}_{\omega^T}(\dd \eta_2)=4\kappa.
\end{equation}
The first identity, together with the \((1,1)\)-type proved in Corollary \ref{cor:static-local-transverse-geometry},
shows that \(\dd \eta_1\) is primitive.  On a Hermitian four-manifold every
real primitive \((1,1)\)-form is anti-self-dual, so
\(\star_K \dd\eta_1=- \dd \eta_1\).  Since
\(\operatorname{tr}_{\omega_K}
=\ee^{2\phi}\operatorname{tr}_{\omega^T}\), we obtain
\begin{equation}
 \begin{aligned}
 R_K
 &=\operatorname{tr}_{\omega_K}(\rho_K)
 =\ee^{2\phi}
   \operatorname{tr}_{\omega^T}(2\kappa \dd \eta_2)\\
 &=8\kappa^2\ee^{2\phi}.
 \end{aligned}
\end{equation}

It remains to rewrite \(\dd H=0\).  Since
\(\omega^T=\ee^{2\phi}\omega_K\) and
\(\ee^{2\phi}=R_K/(8\kappa^2)\), the horizontal torsion is
\begin{equation}
 H^T=\dd^{\mathrm{c}}\omega^T
 =\frac1{8\kappa^2}\dd^{\mathrm{c}}R_K\wedge\omega_K.
\end{equation}
Together with
\( \dd \eta_2=(2\kappa)^{-1}\rho_K\), this gives the asserted formula for
\(H\).  Differentiation gives
\begin{equation}
 \frac1{8\kappa^2}\dd \dd^{\mathrm{c}}R_K\wedge\omega_K
 + \dd\eta_1 \wedge \dd \eta_1
 +\frac1{4\kappa^2}\rho_K\wedge\rho_K=0.
\end{equation}
After multiplication by \(4\kappa^2\), this is
\begin{equation}
 \frac12 \dd \dd^{\mathrm{c}}R_K\wedge\omega_K
 +4\kappa^2  \dd\eta_1 \wedge \dd \eta_1
 +\rho_K\wedge\rho_K=0.
\end{equation}

The Hodge star with respect to the K\"ahler metric gives
\begin{equation}
 -\frac12\Delta_KR_K
 -2\kappa^2|\dd\eta_1|^2
 -\frac12|\operatorname{Rc}_K|^2
 +\frac14R_K^2=0,
\end{equation}
which is exactly
\eqref{eq:real-kappa-scalar-master-equation}.

For the converse, the scalar box equation is equivalent, by the same
Hodge-star calculation, to the preceding four-form equation.  Hence the
displayed three-form \(H\) is closed.  The closed primitive
\((1,1)\)-form \( \dd \eta_1\) and the Ricci form provide the two local
connection curvatures, and the standard inverse construction for
six-dimensional Bismut--Hermitian--Einstein metrics applies; see
\cite[Theorem~1.1]{ALS}.  In our conventions the
vertical complex structure is fixed by
\(J V_1=-V_2\), \(J V_2=V_1\).
A direct calculation gives \(H=\dd^{\mathrm{c}}\omega\), while for a $\R^2$ bundle, 
\begin{equation}
\rho^H=\rho_K - 2\kappa \cdot (2\kappa)^{-1} \rho_K,
\end{equation} 
where we have used \cite[Proposition 4]{GGP} together with the fact that the Bismut Ricci form is conformally invariant in dimension $4$. Since \(\dd H=0\), the metric
is pluriclosed and therefore Bismut--Hermitian--Einstein.
\end{proof}

\begin{remark}
The corresponding canonical reduction of compact non-K\"ahler BHE
threefolds to a transverse K\"ahler surface, together with its local
converse, is established in
\cite[Theorem~1.1 and Proposition~2.18]{ALS}. For the global quasi-regular reconstruction and the corresponding
integrality conditions on the two curvature classes, see
\cite[ Section 2.1 and Definition~2.2]{ALL-26}
and \cite[Theorem~1.1]{ALS}.
\end{remark}

\subsection{Purely imaginary \(\kappa\)} 

Below we assume that 
\(\kappa=\ii\lambda\), \(\lambda\in\mathbb R\setminus\{0\}\). We will refer to the corresponding geometries as \BHEI structures. We start by proving the following no-go theorem
\begin{theorem}\label{prop:imaginary-kappa-noncompact}
Let $(M^6,g,\epsilon,H,\phi,\kappa, dH=0)$ be a static solution of the
GKSE with $\epsilon_\pm$ of unit norm. If
\begin{equation}
\operatorname{Im}(\kappa)\neq0,
\end{equation}
then $M$ is non-compact. In particular, there are no compact solutions
in the branch $\kappa\in \ii\mathbb{R}\setminus\{0\}$.
\end{theorem}

\begin{proof}
Suppose, by contradiction, that $M$ is compact. By
Lemma~\ref{lem:dilaton-vertical-derivatives},
\begin{equation}
V_2(\phi)=\mathcal L_{V_2}\phi=\operatorname{Im}(\kappa).
\end{equation}
Moreover, $V_2$ is Killing and hence divergence-free:
\begin{equation}
\operatorname{div}_g(V_2)=0.
\end{equation}
The divergence theorem therefore gives
\begin{equation}
\begin{aligned}
0
&=
\int_M\operatorname{div}_g(\phi V_2)\,\dd V_g\\
&=
\int_M\left(
V_2(\phi)+\phi\,\operatorname{div}_g(V_2)
\right)\dd V_g\\
&=
\operatorname{Im}(\kappa)\operatorname{Vol}_g(M).
\end{aligned}
\end{equation}
Since $\operatorname{Vol}_g(M)>0$, it follows that
$\operatorname{Im}(\kappa)=0$, contradicting the assumption. Hence
$M$ cannot be compact.
\end{proof}

\begin{theorem}\label{thm:imaginary-kappa-local-model}
Let \((M^6,g,\epsilon,H,\phi,\kappa)\) be a  solution
of the GKSE  with $\epsilon$ static and $\epsilon_\pm$ with unit norm. Assume
\(\kappa=\ii\lambda\), \(\lambda\in\mathbb R\setminus\{0\}\).
Every point has a neighborhood
which is a local principal \(\mathbb R^2\)-bundle
\begin{equation}
 \pi:U\longrightarrow B^4
\end{equation}
whose vertical fields are \(V_1,V_2\) and whose connection one-forms are
\(\eta_1,\eta_2\).  The local leaf space \(B\) carries a hyper-Hermitian
structure whose fundamental forms are
\begin{equation}
 \omega_1=\omega^T,
 \qquad
 \omega_2=\operatorname{Re}\Omega^T,
 \qquad
 \omega_3=-\operatorname{Im}\Omega^T.
\end{equation}
Their common Lee form is the basic one-form
\begin{equation}
 \theta^T=2\dd\phi-2\lambda\eta_2.
\end{equation}
Moreover, $\dd\eta_1, \dd\eta_2$
are closed anti-self-dual transverse two-forms and
\(\dd\theta^T=-2\lambda \dd\eta_2\).  The Bianchi identity is equivalent to
\begin{equation}\label{eq:box-equation-imaginary}
 \dd_T^*\theta^T
 -\frac12|\dd \eta_1|^2
 -\frac1{8\lambda^2}|\dd\theta^T|^2=0.
\end{equation}

Conversely, let
\((B^4,g^T,I_1,I_2,I_3)\) be a hyper-Hermitian surface with Lee
form \(\theta^T\), and let \(\alpha,\beta\) be closed anti-self-dual
two-forms satisfying $
 \dd\theta^T=-2\lambda\beta$
and $$\dd_T^*\theta^T
 -\frac12|\alpha|^2
 -\frac1{8\lambda^2}|\dd\theta^T|^2=0.$$  Locally choose connection
one-forms on \(B\times\mathbb R^2\) with
\(\dd\eta_1=\alpha\), \(\dd\eta_2=\beta\).  Then
\begin{equation}
 g=g^T+\eta_1^2+\eta_2^2,
 \qquad
 \omega=\omega_1+\eta_1\wedge\eta_2,
\end{equation}
define a six-dimensional Bismut--Hermitian--Einstein structure.
\end{theorem}

\begin{proof}
When \(\kappa=\ii\lambda\), equations
\eqref{eq:weighted-transverse-fundamental-form} and
\eqref{eq:transverse-holomorphic-form} become
\begin{equation}
 \begin{aligned}
 \dd\bigl(\ee^{-2\phi}\omega^T\bigr)
 &=-2\lambda\eta_2\wedge \ee^{-2\phi}\omega^T,\\
 \dd\bigl(\ee^{-2\phi}\Omega^T\bigr)
 &=-2\lambda\eta_2\wedge \ee^{-2\phi}\Omega^T.
 \end{aligned}
\end{equation}
After expanding the factor \(\ee^{-2\phi}\), both equations are precisely
\begin{equation}
 \dd\omega^T=\theta^T\wedge\omega^T,
 \qquad
 \dd\Omega^T=\theta^T\wedge\Omega^T,
 \qquad
 \theta^T=2\dd\phi-2\lambda\eta_2.
\end{equation}
Because \(\operatorname{Re}(\kappa)=0\), Corollary
\ref{cor:Lie-derivatives-structure-forms} shows that
\(\omega^T\) and \(\Omega^T\) are global basic transverse forms.
The three real forms
\(\omega^T,\operatorname{Re}\Omega^T,-\operatorname{Im}\Omega^T\)
satisfy the quaternionic algebraic identities of Proposition
\ref{prop:spin-static-SU2-normal-form} and have the same Lee form.
They therefore define the asserted transverse hyper-Hermitian structure.

Differentiating the global identity for \(\theta^T\) gives
\begin{equation}
 \dd\theta^T=-2\lambda\,\dd\eta_2.
\end{equation}
Furthermore, Lemma~\ref{lem:Lee-forms} gives
\(\theta^+=\theta^-=\theta^T\).  Comparing this with
\eqref{eq:Lee-trace-decomposition} gives, separately in the two vertical
directions,
\begin{equation}
 \operatorname{tr}_{\omega^T}(\dd \eta_1)=0,
 \qquad
 \operatorname{tr}_{\omega^T}(\dd \eta_2)=0.
\end{equation}
Since \(\dd \eta_i\) already have type \((1,1)\) (see Corollary \ref{cor:static-local-transverse-geometry}), both are primitive
and hence anti-self-dual.

Applying \(\star_T\) to
\eqref{eq:universal-transverse-Bianchi}, and using
\(\dd_T^*\theta^T=-\star_T\dd\star_T\theta^T\), gives
\begin{equation}
 \dd_T^*\theta^T-\frac12| \dd \eta_1|^2-\frac12|\dd \eta_2|^2=0.
\end{equation}
Substitution of
\(\dd \eta_2=-(2\lambda)^{-1}\dd\theta^T\) yields the box equation.

For the converse, anti-self-dual forms on a hyper-Hermitian surface are
of type \((1,1)\) for each of its compatible complex structures.  Hence
the horizontal complex structure \(I_1\), extended by
\(J V_1=-V_2\), \(J V_2=V_1\), is integrable on the local
\(\mathbb R^2\)-bundle. The scalar equation $\dd_T^*\theta^T
 -\frac12|\alpha|^2
 -\frac1{8\lambda^2}|\dd\theta^T|^2=0$ is equivalent to
\begin{equation}
 -\dd\star_T\theta^T+\alpha\wedge\alpha+\beta\wedge\beta=0,
\end{equation}
so \(\dd H=0\).  Finally, applying \cite[Proposition 4]{GGP}, 
\begin{equation}
\rho^H=\rho^{H^{T}}=0,
\end{equation} 
where $\rho^T$ is the Bismut Ricci form of the transverse metric which vanishes since the transverse geometry is hyper-Hermitian. The resulting pluriclosed Hermitian structure is therefore Bismut--Hermitian--Einstein.
\end{proof}
\begin{remark}
By \cite[Formula 2.12]{AI}, for any hyper-Hermitian (HKT) surface
\begin{equation}
    \dd^*\leet
    =
    \frac{1}{3}R^T-\frac{1}{2}|\leet|^2.
\end{equation}
Combining this identity with~\eqref{eq:box-equation-imaginary} gives
\begin{equation}
    R^T
    =
    \frac{3}{2}|\leet|^2
    +\frac{3}{8\lambda^2}|\dd\leet|^2
    +\frac{3}{2}|\alpha|^2.
\end{equation}
In particular, the transverse scalar curvature is non-negative.

\end{remark}

The results of this section are summarized in Table \ref{tab:real-vs-imaginary}.
\begin{table}[ht] 
\centering
\renewcommand{\arraystretch}{1.8}
\setlength{\tabcolsep}{7pt}
\newcommand{\sumlabel}[1]{{\footnotesize #1}}
\newcommand{\sumcell}[1]{$\displaystyle #1$}
\begin{tabular}{| l | c | c |} 
\hline
 & \sumcell{\boldsymbol{\kappa\in\R}} & \sumcell{\boldsymbol{\kappa=\mathrm i\lambda\in\mathrm i\R}} \\
\hline
\sumlabel{Transverse $\omega^T$}
 & \sumcell{\mathrm d\big(\mathrm e^{-2\phi}\omega^T\big)=0}
 & \sumcell{\mathrm d\big(\mathrm e^{-2\phi}\omega^T\big)=-2\lambda\,\eta_2\wedge \mathrm e^{-2\phi}\omega^T} \\
\sumlabel{Transverse $\Omega^T$}
 & \sumcell{\mathrm d\big(\mathrm e^{-2\phi}\Omega^T\big)=2\mathrm i\kappa\,\eta_2\wedge \mathrm e^{-2\phi}\Omega^T}
 & \sumcell{\mathrm d\big(\mathrm e^{-2\phi}\Omega^T\big)=-2\lambda\,\eta_2\wedge \mathrm e^{-2\phi}\Omega^T} \\
\sumlabel{Transverse metric}
 & conformally K\"ahler,\quad \sumcell{g^T=\mathrm e^{2\phi}g_K}
 & hyper-Hermitian \\
\sumlabel{Transverse Lee form}
 & \sumcell{\theta^T=2\,\mathrm d\phi}
 & \sumcell{\theta^T=2\,\mathrm d\phi-2\lambda\,\eta_2} \\
\sumlabel{BHE metric}
 & \sumcell{g_6=\mathrm e^{2\phi}g_K+\eta_1^2+\eta_2^2,\quad \mathrm e^{2\phi}=\frac{R}{8\kappa^2}}
 & \sumcell{g_6=g^T+\eta_1^2+\eta_2^2} \\
\sumlabel{Fibration}
 & \sumcell{\mathrm d\eta_1=\alpha,\quad \mathrm d\eta_2=\frac{1}{2\kappa}\,\rho_K}
 & \sumcell{\mathrm d\eta_1=\alpha,\quad \mathrm d\eta_2=\beta=-\frac{1}{2\lambda}\,\mathrm d\theta^T} \\
\sumlabel{$\mathrm dH=0$}
 & \sumcell{\Box R=\tfrac12R^2-R_{mn}R^{mn}-4\kappa^2\alpha_{mn}\alpha^{mn}}
 & \sumcell{\nabla^m\theta^T_m+\tfrac12\big(\alpha_{mn}\alpha^{mn}+\beta_{mn}\beta^{mn}\big)=0} \\
\sumlabel{Compact $M_6$}
 & allowed
 & excluded \\
 \hline
\end{tabular}

\caption{Summary of the two classes of solutions of the GKSE (for a static spinor with $|\epsilon_\pm|=1$); see Theorems \ref{thm:real-kappa-local-model} and \ref{thm:imaginary-kappa-local-model}. In the left column $g_K$ is a K\"ahler metric on the four-dimensional base, with Ricci form $\rho_K$, Ricci tensor $R_{mn}$ and scalar curvature $R$; in the right column $g^T$ is a hyper-Hermitian metric. Indices in the ``Box equations'' are raised with $g_K$ in the left column and $g^T$, in the right one.}
\label{tab:real-vs-imaginary}
\end{table}

\section{Sufficiency} \label{sec:sufficiency}
We have shown that our spinorial system~\eqref{eq:KSE-chiral-decomposition} implies the existence of two BHE structures that satisfy equations ~\eqref{eqn:symplectic},~\eqref{eq:weighted-holomorphic-forms} and Lemma \ref{lem:Lee-forms}. A natural question, at this point, is whether the converse is true, as well. Now our starting point is an $\SU(2)$ structure as in Proposition ~\ref{prop:spin-static-SU2-normal-form}, whose associated $\SU(3)$ structures are defined by~\eqref{eq:spin-SU2-to-SU3-omega-global} and~\eqref{eq:spin-SU2-to-SU3-Omega-global} and satisfy the system
\begin{subequations} \label{eq:modified-Hull-Strominger}
\begin{align}
        \dd\holopm &=\leepm\wedge\holopm\,, \label{holopm-integrability}\\
    \dd(\fundpm\wedge\fundpm) &=\theta^\pm\wedge\fundpm\wedge\fundpm\,, \label{conformally-bal-generalised} \\
    \pm\star H &=\dd\fundpm-\leepm\wedge\fundpm\,, \label{fundamental-form-su3pm-intrinsic-torsion}
\end{align}
where
\begin{align} 
      \leep=2\dd\phi-2\mathrm{Re}(\kappa  \vformbar)\,, \label{eq:sufficiency-Leep-form}\\
    \leem= 2\dd\phi+2\mathrm{Re}(\kappa \vform) \, \label{eq:sufficiency-Leem-form},
\end{align}
\end{subequations}
and $dH=0$.
It is clear that the two structures $(M,\fundpm,\holopm)$ satisfy the usual Hull--Strominger system, but with the replacement $2\dd\phi\mapsto \leepm$, where $\leepm$ is in general \emph{not} closed.

We consider now the structure $(M,\fundp,\holop)$, given that the argument proceeds in a similar fashion for the other one with opposite orientation. Observe that~\eqref{holopm-integrability}, by the standard theory of $\SU(3)$ intrinsic torsion (see for example~\cite{chiossi2002intrinsictorsionsu3g2}), implies that the associated almost complex structure $\Jp$ is integrable, while~\eqref{conformally-bal-generalised} implies that the Lee form associated to $\fundp$ is exactly $\leep$. The final equation implies that $H=\Jp\dd\fundp$, so $H$ is exactly the torsion of the Bismut connection, namely the unique 3-form which preserves $\Jp$, that is $\nabla^H\fundp=0$.   

\subsection{Dilatino equation}
On a spin $6$-manifold, the existence of such a $\SU(3)$ structure implies the existence of a Weyl unit spinor  $\epsilon_+$ defining $(g,J^+,\omega^+,\Omega^+)$ and compatible with the $\SU(2)$-reduction. Given a local orthonormal coframe $\{e^a\}_{a=1}^6$, let
$\gamma_a$ denote the corresponding Clifford matrices in
the Majorana realization fixed in
Section~\ref{sec:spinor-preliminaries}.
The eight spinors $\epsilon_+^*$, $\ii\epsilon_+^*$ and
$\gamma_a\epsilon_+$, $a=1,\ldots,6$, form a basis of $S^-$
over $\mathbb R$. Since Clifford multiplication by the three-form $H$ reverses
chirality, $H\cdot\epsilon_+\in S^-$ and hence
\begin{equation}
     H\cdot\epsilon_+
 =u\epsilon_+^*+v\ii\epsilon_+^*
  +\sum_{a=1}^6w_a\gamma_a\epsilon_+,
 \qquad u,v,w_a\in\mathbb R.
\end{equation}
Setting $a=u+\ii v$ and $w=\sum_{a=1}^6 w_a e^a$, we obtain
\begin{equation}
 H\cdot\epsilon_+=a\epsilon_+^*+w\cdot\epsilon_+,
\end{equation}
where $a$ is a complex function and $w$ is a real one-form.
To determine $a$, pair the decomposition with $\epsilon_+$
using $\mathsf B$. Since $|\epsilon_+|=1$, $
 \mathsf B(\epsilon_+,\epsilon_+^*)=1,$ and $
 \mathsf B(\epsilon_+,w\cdot\epsilon_+)=0.$
Hence the definition of $\Omega^+$ gives $
 a=\mathsf B(\epsilon_+,H\cdot\epsilon_+)
   =-\frac16g_\C (H, \Omega^+)=0,$
because $H$ has type $(2,1)+(1,2)$. Thus $H\cdot\epsilon_+=w\cdot\epsilon_+$.
To determine $w$, let $X$ be a real vector field and take
the imaginary part of the Hermitian pairing with
$(J^+X)^\flat\cdot\epsilon_+$. Since $w$ is real, the definition of $\omega^+$ gives
\begin{equation}
\operatorname{Im}
 \langle (J^+X)^\flat\cdot\epsilon_+,w\cdot\epsilon_+\rangle
 =\omega^+(J^+X,w^\sharp)=w(X).
\end{equation}
On the other hand, Clifford multiplication by $H$ is
skew-Hermitian, so
\begin{equation}
\begin{split}
 2\ii\operatorname{Im}
 \langle (J^+X)^\flat\cdot\epsilon_+,H\cdot\epsilon_+\rangle&=\langle\epsilon_+,
 \{(J^+X)^\flat\cdot,H\cdot\}\epsilon_+\rangle=2\langle\epsilon_+,(\iota_{J^+X}H)\cdot\epsilon_+\rangle\\
 &=\ii(\operatorname{tr}_{\omega^+}H)(J^+X)=2\ii\theta^+(X),
 \end{split}
\end{equation}
where we used the two-form bilinear identity and
\eqref{eq:Lambda-H-Lee}.
Comparing the two expressions gives $w=\theta^+$.
Therefore $H\cdot\epsilon_+=\theta^+\cdot\epsilon_+.$
After performing the analogous computation for the other structure, we find
\begin{equation} \label{eq:gravitino-converse-lee-form}
    (H-\leepm)\cdot\epspm=0\,.
\end{equation}
From~\eqref{eq:static-vector-Clifford-action} we have
\begin{equation}
    \vformre\cdot\epspm=\epsmp\,,\quad \vformim\cdot \epspm=\pm\ii\epsmp\,,
\end{equation}
and after using
\begin{subequations}
  \begin{align}
          \mathrm{Re}(\kappa\bar\eta) &=\mathrm{Re}(\kappa)\vformre+\mathrm{Im}(\kappa)\vformim\,, \\
          \mathrm{Re}(\kappa\eta) &=\mathrm{Re}(\kappa)\vformre-\mathrm{Im}(\kappa)\vformim\,,
  \end{align}  
\end{subequations}
we find
\begin{equation}
    \mathrm{Re}(\kappa\vformbar)\cdot\epsp=\kappa \epsm\,,\quad \mathrm{Re}(\kappa\vform)\cdot\epsm=\kappa \epsp\,.
\end{equation}
Therefore,~\eqref{eq:gravitino-converse-lee-form} becomes
\begin{equation}
   \left(\dd\phi - \frac{1}{2} H\right) \cdot \epspm \mp \kappa \epsmp = 0\,. 
\end{equation}
\subsection{Gravitino equation}
We have shown that $\nabla^H$ is the Bismut connection of
$(g,J^+)$. It therefore preserves the complex line
$\Lambda^{0,3}_{J^+}$ and its Hermitian metric.
Since $\Omega^+$ has constant norm, there exists a real
one-form $\alpha$ such that $\nabla^H\Omega^+=\ii\alpha\otimes\Omega^+.$

Let $\{e_a\}$ be a local orthonormal frame, with dual
coframe $\{e^a\}$. Since the torsion of $\nabla^H$ is $H$,
exterior differentiation gives
\begin{align}
 \dd\Omega^+
 &=\sum_a e^a\wedge\nabla^H_{e_a}\Omega^+
   +\sum_a(\iota_{e_a}H)\wedge(\iota_{e_a}\Omega^+)
 \nonumber\\
 &=\ii\alpha\wedge\Omega^+
   +\sum_a(\iota_{e_a}H)\wedge(\iota_{e_a}\Omega^+).
\end{align}

We compute the torsion term using the decomposition $
 H=\frac12J^+\theta^+\wedge\omega^++H_0,$ and $
 \operatorname{tr}_{\omega^+}H_0=0,$
where $H_0$ has type $(2,1)+(1,2)$.
The algebraic identities of the $\SU(3)$ structure give
\begin{equation}
     \sum_a
 \iota_{e_a}(J^+\theta^+\wedge\omega^+)
 \wedge\iota_{e_a}\Omega^+
 =2\theta^+\wedge\Omega^+,
 \qquad
 \sum_a(\iota_{e_a}H_0)\wedge(\iota_{e_a}\Omega^+)=0.
\end{equation}
Consequently, $
 \dd\Omega^+=(\ii\alpha+\theta^+)\wedge\Omega^+.$

Comparing with \eqref{holopm-integrability}, we obtain
$\alpha\wedge\Omega^+=0$.
Since $\Omega^+$ is a nowhere-vanishing $(0,3)$-form,
wedging with $\Omega^+$ is injective on $(1,0)$-forms.
Thus $\alpha^{1,0}=0$. As $\alpha$ is real, also
$\alpha^{0,1}=0$, and hence $\alpha=0$.

Therefore $\nabla^H$ preserves $\Omega^+$, as well as
$g$ and $J^+$. It preserves the corresponding
$\SU(3)$-reduction, so its spinorial lift annihilates
the defining unit spinor $\epsilon_+$.
The same argument applies to the negative-chirality
structure, giving $\nabla^H\epsilon_\pm=0.$

\section{Examples}\label{section:examples}
This section illustrates the reverse constructions established in
Theorems~\ref{thm:real-kappa-local-model} and
\ref{thm:imaginary-kappa-local-model}. We distinguish the real and
purely imaginary cases and describe their relation to known constructions.

\subsection{Examples of \BHER geometries}
\label{sec:real-kappa-examples}

In the compact setting, the construction of Theorem~\ref{thm:real-kappa-local-model} (see also \cite{ALS}) has already
been exploited in \cite{ALL-26, apostolov2026gluingconstructionzextremalkahler}. The following non-compact examples are
related to the quadratic K\"ahler ansatz of
\cite[Appendix~C]{Donos:2008ug} and recover members of the
cohomogeneity-one family of
\cite{podesta2026newinhomogeneousbismuthermitianeinstein}.
We recall their transverse description in our conventions.

\begin{example}
\label{ex:real-complete-line-bundles}
Fix $\kappa=\tfrac12$ and $n\in\{1,2,3\}$. Let
$\Sigma_n=\mathbb{CP}^1$ carry the round K\"ahler metric of Gaussian
curvature $k_n=2/n$. We use the orientation for the bundle periods in
which
\begin{equation}\label{eq:base-area}
 \int_{\Sigma_n}\omega_{\Sigma_n}=2\pi n.
\end{equation}
Let $P_n\to\Sigma_n$ be the circle bundle with $2\pi$-periodic fibres
and connection $\vartheta$ satisfying $\dd\vartheta=\omega_{\Sigma_n}$.
Let $B_n=\mathcal O_{\mathbb{CP}^1}(-n)$ be the associated line bundle
for the weight $-1$ representation.
On the complement of the zero section, put $z=1+r^2/2$, where $r$
is the fibre radius, and set
\begin{equation}\label{eq:parameters}
 q_n=\frac1n,\qquad
 p_n=\frac{4n-n^2-1}{2n},\qquad
 c_n=k_n-p_n>0.
\end{equation}
The profile and K\"ahler structure are
\begin{align}
 F_n(z)&=(z-1)(p_nz+2-p_n),\label{eq:profile}\\
 g_{K,n}&=z\,g_{\Sigma_n}
       +\frac{z}{F_n(z)}\,\dd z^2
       +\frac{F_n(z)}{z}\,\vartheta^2,\label{eq:gK}\\
 \omega_{K,n}&=z\,\omega_{\Sigma_n}+\dd z\wedge\vartheta.
 \label{eq:omegaK}
\end{align}
The conditions $F_n(1)=0$ and $F_n'(1)=2$ give a smooth extension
across the zero section. The Ricci form, scalar curvature and closed
primitive $(1,1)$-form needed in the construction are
\begin{align}
 \rho_{K,n}
 &=\left(c_n-\frac{1-p_n}{z}\right)\omega_{\Sigma_n}
   +\frac{1-p_n}{z^2}\,\dd z\wedge\vartheta,\label{eq:rho}\\
 R_{K,n}&=\frac{2c_n}{z}>0,\label{eq:scalar}\\
 \alpha_n
 &=q_n\left(\frac1z\omega_{\Sigma_n}
              -\frac1{z^2}\dd z\wedge\vartheta\right).
 \label{eq:alpha}
\end{align}
With our tensor-norm convention, the reduced equation is
\begin{equation}\label{eq:reduced-equation}
 -\Delta_{K,n}R_{K,n}+\frac12R_{K,n}^2
 -|\operatorname{Rc}_{K,n}|^2-|\alpha_n|^2=0.
\end{equation}
It reduces to the identity
\begin{equation}\label{eq:parameter-relation}
 q_n^2=k_n(2-p_n)-1,
\end{equation}
which holds for the parameters above. Moreover,
\begin{equation}\label{eq:real-example-periods}
 \frac1{2\pi}\int_{\Sigma_n}\alpha_n=1,
 \qquad
 \frac1{2\pi}\int_{\Sigma_n}\rho_{K,n}=2-n.
\end{equation}
Thus there is a principal $\mathbb T^2$-bundle $M_n\to B_n$ with
connection forms satisfying
\begin{equation}\label{eq:eta-curvature}
 \dd\eta_1=\alpha_n,\qquad \dd\eta_2=\rho_{K,n}.
\end{equation}
Theorem~\ref{thm:real-kappa-local-model} gives the BHE metric
\begin{equation}\label{eq:real-example-metric}
 g_n=\frac{c_n}{z}g_{K,n}+\eta_1^2+\eta_2^2,
 \qquad
 \phi_n=\frac12\log\frac{c_n}{z}.
\end{equation}
The Chern vector $(1,2-n)$ is primitive. After an integral change of
basis of the torus, the bundle is the product of a trivial circle with
the pullback of the Hopf bundle to $B_n$. Since the pullback of
$\mathcal O(-n)$ to $S^3$ is topologically trivial, this gives
\begin{equation}\label{eq:topology}
 M_n\simeq S^3\times S^1\times\mathbb R^2.
\end{equation}
The transverse metric is
\begin{equation}\label{eq:gT-explicit}
 g_n^T=c_n\left(g_{\Sigma_n}+\frac{\dd z^2}{F_n(z)}
                    +\frac{F_n(z)}{z^2}\vartheta^2\right).
\end{equation}
Since $F_n(z)\sim p_nz^2$ at infinity, the radial length is infinite.
Together with smoothness at the zero section, this proves completeness
of $g_n^T$, and hence of the connection metric $g_n$ with compact torus
fibres.

More precisely, after the homothety $g_n\mapsto g_n/4$, these 
metrics agree locally with those of
\cite[Theorem~3.13]{podesta2026newinhomogeneousbismuthermitianeinstein}
with parameters
\begin{equation}
 m=1,\qquad a=\frac12,\qquad b=\frac{n-2}{2},\qquad c=n,
 \qquad D=a^2+b^2=\frac{1+(n-2)^2}{4}<1.
\end{equation}
\end{example}

\subsection{Examples of \BHEI geometries}
\label{sec:imaginary-kappa-examples}

Throughout this subsection, let $\kappa=\ii\lambda$ with
$\lambda\in\mathbb R\setminus\{0\}$. The first example uses the
classical harmonic-function construction. The second starts from a
known hyper-Hermitian surface with non-closed Lee form and gives an
explicit solution of the transverse scalar equation.

\begin{example}
\label{subsec:imaginary-explicit-HK}
Let $N=\mathbb R^4\setminus\{a,-a\}$, where $a\neq0$, and let
$(h,\omega_1^{HK},\omega_2^{HK},\omega_3^{HK})$ be its standard flat
hyper-K\"ahler structure. Set
\begin{equation}\label{eq:two-centre-harmonic-function}
 q(X)=\frac1{|X-a|^2}+\frac1{|X+a|^2}.
\end{equation}
The conformal structure
\begin{equation}\label{eq:two-centre-transverse-structure}
 g^T=qh,\qquad \omega_L^T=q\,\omega_L^{HK},\qquad L=1,2,3,
\end{equation}
is the standard construction from a positive harmonic function on a
hyper-K\"ahler four-manifold; see \cite[Section~3, equation~(3.4)]{Gauntlett:2003cy}.
Its common Lee form is $\theta^T=\dd\log q$, and
$\star_T\theta^T=\star_h \dd q$ gives $\dd_T^*\theta^T=0$.

On $M=N\times\mathbb R^2$, take
\begin{equation}\label{eq:two-centre-six-dimensional-data}
 \eta_1=\dd t_1,\qquad \eta_2=\dd t_2,\qquad
 g_6=qh+\dd t_1^2+\dd t_2^2.
\end{equation}
Since $\dd\eta_1=\dd\eta_2=0$, $\dd\theta^T=0$ and $\dd_T^*\theta^T=0$,
the hypotheses of Theorem~\ref{thm:imaginary-kappa-local-model} hold,
and the forms $\omega^\pm=q\omega_1^{HK}\pm \dd t_1\wedge \dd t_2$
define two BHE structures. The metric $g_6$ is the two-centre
Callan--Harvey--Strominger metric \cite{Callan:1991dj}.
Its Bismut holonomy is
$\operatorname{SU}(2)$ by
\cite[Proposition~3.3]{foscolo2026cohomogeneitysolutionsiib}.
In both descriptions, the transverse Lee form is
\begin{equation}
 \theta^T=\dd\log q.
\end{equation}

Near either removed point, $q=r_\pm^{-2}+O(1)$, while at infinity
$q=2r^{-2}+O(r^{-4})$. All three ends of $(N,qh)$ are therefore
asymptotically cylindrical and complete. The product metric $g_6$ is
complete as well.
\end{example}

\begin{example}
\label{ex:nonclosed-transverse-Lee}
Let $B=(0,\infty)_t\times\mathbb T^3_{x,y,z}$, with $x,y,z$ of period
$2\pi$, and put
\begin{equation}
 A=\ee^{-t}(\sin z\,\dd x+\cos z\,\dd y),\qquad
 e^1=\dd t+A,\quad e^2=\dd x,\quad e^3=\dd y,\quad e^4=\dd z.
\end{equation}
The metric $h=\sum_{a=1}^4(e^a)^2$ is a special case of the
hyper-Hermitian construction of
\cite[Section~3.2]{DunajskiGrantStrachan2007}. With orientation
$\operatorname{vol}_h=e^{1234}$, its fundamental forms are
\begin{equation}
 \omega_1=e^{12}+e^{34},\qquad
 \omega_2=e^{13}-e^{24},\qquad
 \omega_3=-e^{14}-e^{23}.
\end{equation}
They satisfy the quaternionic identities in our convention
$\omega_i(X,Y)=h(X,J_iY)$ and
$\dd\omega_i=\theta\wedge\omega_i$, with common Lee form $\theta=A$.
In particular,
\begin{equation}
 \dd\theta=-\ee^{-t}\sin z\,(e^{12}-e^{34})
          -\ee^{-t}\cos z\,(e^{13}+e^{24}),
 \qquad \star_h \dd\theta=-\dd\theta,
 \qquad |\dd\theta|_h^2=4\ee^{-2t}.
\end{equation}

Define
\begin{equation}
 q(t)=\frac1{4\lambda^2}
       \left(\frac2{\sqrt{1+\ee^{-2t}}}-1\right)>0,
 \qquad g^T=qh,\qquad \theta^T=\dd\log q+\theta.
\end{equation}
On $M=B\times\mathbb R^2_{\psi_1,\psi_2}$, set
\begin{equation}
 \eta_1=\dd\psi_1,\qquad
 \eta_2=\dd\psi_2-\frac1{2\lambda}\theta,
 \qquad \phi=\frac12\log q+\lambda\psi_2.
\end{equation}
Then
\begin{equation}
 \dd\eta_1=0,\qquad \dd\eta_2=-\frac1{2\lambda}\dd\theta\neq0,
 \qquad \theta^T=2\dd\phi-2\lambda\eta_2.
\end{equation}
Both curvature forms are anti-self-dual for $g^T$.

For a function of $t$, write
\begin{equation}
 L_hq:=\operatorname{div}_h
       \bigl(\operatorname{grad}_h q+q\theta^{\sharp_h}\bigr)
       =(1+\ee^{-2t})q''-3\ee^{-2t}q'+2\ee^{-2t}q.
\end{equation}
The chosen $q$ satisfies $L_hq=-\ee^{-2t}/(2\lambda^2)$. Consequently,
\begin{equation}
 \dd_T^*\theta^T=-q^{-2}L_hq
       =\frac{\ee^{-2t}}{2\lambda^2q^2},
 \qquad |\dd\theta^T|_{g^T}^2=\frac{4\ee^{-2t}}{q^2},
\end{equation}
and hence
\begin{equation}
 \dd_T^*\theta^T-\frac12|\dd\eta_1|_{g^T}^2
       -\frac1{8\lambda^2}|\dd\theta^T|_{g^T}^2=0.
\end{equation}
Theorem~\ref{thm:imaginary-kappa-local-model} now gives two BHE
structures with metric and fundamental forms
\begin{equation}
 g=qh+\eta_1^2+\eta_2^2,
 \qquad \omega^\pm=q\omega_1\pm\eta_1\wedge\eta_2.
\end{equation}
The metric $g$ is incomplete on the chosen domain, since the
curves with all coordinates fixed except $t$  approach $t=0$
with finite length.
\end{example}

\section{Global reduction to five dimensions} \label{Section7}

We consider the regular case, in which $\xi:=V_1$ generates a free
isometric circle action, and assume that this action lifts to the spin
structure as an $S^1$-action. The quotient is then a smooth spin
manifold:
\begin{equation}\label{eq:global-circle-reduction}
 S^1\longrightarrow M^6\xrightarrow{\ \pi\ }N^5.
\end{equation}
Since $\xi$ has unit norm, $\eta_1=\xi^\flat$ is a connection form,
and $\dd\pi$ identifies $\xi^\perp$ with $\pi^*TN$. The invariant
geometric data determine a metric $g_5$, a function $\phi_5$, a
two-form $\alpha$, and a three-form $H_5$ on $N$, with
\begin{align}
 g&=\pi^*g_5+\eta_1\otimes\eta_1,
 &\phi&=\pi^*\phi_5,
 \label{eq:global-metric-curvature}\\
 \dd\eta_1&=\pi^*\alpha,
 &H&=\pi^*H_5+\eta_1\wedge\pi^*\alpha.
 \label{eq:global-H-decomposition}
\end{align}
Indeed, $\iota_\xi H=\dd\eta_1$, so $H-\eta_1\wedge \dd\eta_1$
is horizontal and invariant and therefore descends to $N$.

\subsection{Reduction of the spinorial equations}

Let us first explain how the spinor descends. In the setting
considered above, $\dd\eta_1$ is a primitive $(1,1)$-form on the
transverse four-dimensional distribution. It therefore belongs to
the Lie algebra of the common $\mathrm{SU}(2)$ stabilizer and
annihilates $\epsilon_\pm$. The vertical gravitino equation
$\nabla^H_\xi\epsilon_\pm=0$, together with the Kosmann formula,
then implies that $\epsilon_\pm$ are invariant under the lifted
circle action.

The unit vector field $\xi$ reduces the spin frame bundle from
$\mathrm{Spin}(6)$ to $\mathrm{Spin}(5)$. Since the negative
half-spin representation restricts to the spin representation in
dimension five, we obtain an equivariant identification
$I_-:S_M^-\longrightarrow\pi^*S_N$.
Thus the invariant unit spinor $\epsilon_-$ determines a unique
unit spinor $\chi\in\Gamma(S_N)$ satisfying
\begin{equation}
 I_-(\epsilon_-)=\pi^*\chi.
\end{equation}
We orient $N$ by
$\operatorname{vol}_6=\pi^*\operatorname{vol}_5\wedge\eta_1$
and fix its Clifford action by
\begin{equation}\label{eq:global-Clifford-convention}
 I_-\bigl(\eta_1\cdot\pi^*u\cdot\psi\bigr)
 =\ii\,\pi^*u\cdot I_-(\psi),
 \qquad u\in\Omega^1(N),\quad\psi\in\Gamma(S_M^-).
\end{equation}
The product on the left preserves chirality. On the right, Clifford
multiplication is identified with the one on $S_N$, with
$u\cdot u=|u|_{g_5}^2\operatorname{id}$ and
$\operatorname{vol}_5\cdot=\operatorname{id}$.
This identification also intertwines the actions of horizontal
two-forms.

Let us now reduce the gravitino equation. For $X\in\Gamma(TN)$,
write $X^h$ for its horizontal lift. In the positive Clifford
convention,
\begin{equation}
 \nabla^H_{X^h}
 =\nabla^g_{X^h}-\frac14(\iota_{X^h}H)\cdot.
\end{equation}
The O'Neill formulas \cite{ON} give the curvature contribution
$\frac14\bigl(\pi^*(\iota_X\alpha)\wedge\eta_1\bigr)\cdot$
to the Levi--Civita spin connection. Since
\begin{equation}
 \iota_{X^h}H
 =\pi^*(\iota_XH_5)+\pi^*(\iota_X\alpha)\wedge\eta_1,
\end{equation}
this contribution cancels with the vertical part of the torsion.
The horizontal and vertical gravitino equations therefore give,
respectively,
\begin{align}
 \nabla^{g_5}_X\chi-\frac14(\iota_XH_5)\cdot\chi&=0,
 \qquad X\in\Gamma(TN),
 \label{eq:global-5d-gravitino}\\
 \alpha\cdot\chi&=0.
 \label{eq:global-alpha-annihilates}
\end{align}

To reduce the dilatino equation, we need the algebraic relation
between $\epsilon_+$ and $\epsilon_-$. Recall that the complex
one-form defining $V$ is
$\eta(X)=\langle\epsilon_-,X^\flat\cdot\epsilon_+\rangle$.
Since $V_1=\operatorname{Re}V$ and the real and imaginary parts of
$V$ are orthonormal, we have
\begin{equation}
 \langle\epsilon_-,\eta_1\cdot\epsilon_+\rangle
 =\eta(V_1)=1.
\end{equation}
Both $\epsilon_-$ and $\eta_1\cdot\epsilon_+$ have unit norm.
Equality in the Cauchy--Schwarz inequality therefore gives
$\eta_1\cdot\epsilon_+=\epsilon_-$, and hence
\begin{equation}\label{eq:global-chiral-relation}
 \epsilon_+=\eta_1\cdot\epsilon_-.
\end{equation}
This is a consequence of the normalized spinorial bilinear on $M$.

The positive-chirality component of the dilatino equation is
\begin{equation}
 \left(\dd\phi-\frac12H\right)\cdot\epsilon_-
 +\kappa\epsilon_+=0.
\end{equation}
The term $\eta_1\wedge\pi^*\alpha$ acts trivially on
$\epsilon_-$ because $\pi^*\alpha\cdot\epsilon_-=0$.
Multiplying the remaining equation by $\eta_1$, using
\eqref{eq:global-chiral-relation}, and applying $I_-$ gives
\begin{equation}\label{eq:global-5d-dilatino}
 \left(\dd\phi_5-\frac12H_5\right)\cdot\chi-\ii\kappa\chi=0.
\end{equation}
Here we used that \eqref{eq:global-Clifford-convention} also holds
for horizontal three-forms. The negative-chirality component gives
the same equation, since horizontal odd forms anticommute with
$\eta_1$. Thus \eqref{eq:global-5d-gravitino},
\eqref{eq:global-alpha-annihilates}, and
\eqref{eq:global-5d-dilatino} give the reduced spinorial system.

Finally, the closedness of the six-dimensional torsion gives
\begin{equation}\label{eq:global-reduced-Bianchi}
 \dd\alpha=0,
 \qquad \dd H_5=-\alpha\wedge\alpha,
\end{equation}
by differentiating \eqref{eq:global-H-decomposition}.

\subsection{Relation with Gauntlett--Kim geometry}

We now rewrite the reduced system in variables adapted to the
comparison with Gauntlett--Kim geometry. Set
\begin{equation}\label{eq:global-GK-fields}
 B:=2\phi_5,
 \qquad F_{\mathrm{GK}}:=\ee^{-B}\star_5H_5,
 \qquad H_5=\ee^B\star_5F_{\mathrm{GK}}.
\end{equation}
With the orientation and Clifford convention fixed above, every
two-form $F$ on $N$ satisfies
\begin{equation}
 (\star_5F)\cdot\chi=-F\cdot\chi,
 \qquad
 (\iota_X\star_5F)\cdot\chi
 =-(X^\flat\wedge F)\cdot\chi.
\end{equation}
Consequently, the reduced gravitino and dilatino equations become
\begin{subequations}\label{eq:global-5d-GK-like-system}
\begin{align}
 \nabla^{g_5}_X\chi
 +\frac14\ee^B(X^\flat\wedge F_{\mathrm{GK}})\cdot\chi&=0,
 \label{eq:global-5d-GK-like-gravitino}\\
 (\dd B)\cdot\chi+\ee^BF_{\mathrm{GK}}\cdot\chi-2\ii\kappa\chi&=0.
 \label{eq:global-5d-GK-like-dilatino}
\end{align}
\end{subequations}
These equations are supplemented by
$\alpha\cdot\chi=0$ and the reduced Bianchi identity
\begin{equation}\label{eq:global-alpha-deformed-GK-equation}
 \dd\bigl(\ee^B\star_5F_{\mathrm{GK}}\bigr)
 =-\alpha\wedge\alpha.
\end{equation}

Moreover, by \eqref{eq:global-H-decomposition} and $\vol_6=\pi^*\vol_5\wedge \eta_1$, we have
\begin{equation}
    \star_6 H=\star_5 H_5\wedge\eta_1-\star_5 \alpha.
\end{equation}
Then, \eqref{eq:steady-gradient-generalized-ricci-soliton} gives
\begin{align}
    \dd F_{\mathrm{GK}} &=0, \\
    \dd(\ee^{-B}\star_5\alpha) &=F_{\mathrm{GK}}\wedge\alpha.
\end{align}
The flux couplings in \eqref{eq:global-5d-GK-like-system} have the same form as those
in the Gauntlett--Kim system
\cite[Equation~(3.6)]{Gauntlett:2007ts}. 

To compare with~\cite{Gauntlett:2007ts}, first observe that their system in (3.6) is not defined for $n=2$, which is the case at hand. Their system can be extended to $n=2$, after a change of normalization. In particular, put $c=\frac{n-2}{2}$ and define
\begin{equation}
    \tilde g =c^{-2} g\,, \quad\ee^{\tilde B}=c^{-2}\ee^{B}\,, \quad \tilde F_{\mathrm{GK}}=c F_{\mathrm{GK}}\,, 
\end{equation}
with the frame $\tilde e_a=c e_a$ that is orthonormal for $\tilde g$. In these new variables, (3.6) of~\cite{Gauntlett:2007ts} becomes
\begin{subequations} \label{eq:GK-KSE-rescaled}
    \begin{align}
        \nabla^{\tilde{g}}_X\chi+\ii\frac{n-2}{4}X^\flat\cdot\chi+\frac{1}{4}\ee^{\tilde B}(X^\flat\wedge \tilde{F}_{\mathrm{GK}})\cdot \chi &=0\,, \label{eq:GK-KSE-rescaled-gravitino}\\
          (\dd \tilde{B})\cdot\chi+\ee^{\tilde{B}}\tilde{F}_{\mathrm{GK}}\cdot\chi+\ii(n-1)\chi&=0 \label{eq:GK-KSE-rescaled-dilatino}\,.
    \end{align}
\end{subequations}
while their (3.7) and the Bianchi identity $\dd F_{\mathrm{GK}}=0$ become
\begin{equation}\label{eq:GK-flux-EOM-Bianchi-rescaled}
    \dd \tilde{F}_{\mathrm{GK}}=0\,, \quad \dd\bigg(\ee^{(3-n)\tilde{B}}\star_{\tilde{g}}\tilde{F}_{\mathrm{GK}}\bigg)=0\,.
\end{equation}
These equations are now well-defined for $n=2$. There is another way to derive the form of the GK equations for $n=2$, they arise from a torus reduction of the seven-dimensional GK geometry (that is $n=3$) on $Y_7=Y_5\times T^2$. For ease of notation, we will drop the tildes for the rescaled structure and denote $F_{\mathrm{GK}}$ with $F$ without subscripts. In the rescaled structure, we have by~\cite[(3.14)]{Gauntlett:2007ts}
\begin{equation}\label{GK-T2-Y5}
   g_7=g_5+\ee^{B}(\dd u_1^2+\dd u_2^2)\,,\quad  F_7=F_5-\dd u_1\wedge \dd u_2\,,
\end{equation}
where $B$ and $F_5$ are the pull-backs on $Y_5$; we also define
\begin{equation}
    \theta^s\coloneqq \ee^{B/2}\dd u_s\,, \text{ where $s=1,2,$} \quad \mathrm{vol}_{T^2}\coloneqq \theta^1\wedge \theta^2\,.
\end{equation}

Decompose the spinor bundle of $Y_7$ as $\tau^* S(Y_5)\otimes \C^2$. A $p$-form $\beta$ on $Y_5$ acts on this spinor bundle as $(\beta\cdot)\otimes \sigma_3^p$ and $\theta^s$ acts as $\mathbf{1}\otimes \sigma_s$. A spinor, constant on $T^2$, decomposes on $\tau^* S(Y_5)\otimes \C^2$ as
\begin{equation}
    \eta=\chi_1\otimes \begin{pmatrix} 1 \\ 0 \end{pmatrix}+\chi_2\otimes \begin{pmatrix} 0 \\ 1 \end{pmatrix}\,,
\end{equation}
where $\chi_1,\chi_2\in\Gamma(S(Y_5))$. Since $\ee^B F_7\cdot=\ee^B F_5\cdot-\vol_{T^2}\cdot$, therefore~\eqref{eq:GK-KSE-rescaled-dilatino} with $n=3$ becomes
\begin{equation} \label{eq:GK7-torus-reduction-dilatino}
    \dd B\cdot\chi_1+\ee^{B}F_5\cdot\chi_1+\ii\chi_1=0 \,, \quad -\dd B\cdot\chi_2+\ee^{B}F_5\cdot\chi_2+3\ii\chi_2=0\,.
\end{equation}
We now move to~\eqref{eq:GK-KSE-rescaled-gravitino} with $n=3$. First, we observe that the spin connection receives a contribution from the warp factor over the torus
\begin{subequations}
\begin{align}
        \nabla^{g_7}_X\eta=\frac{1}{4}X\cdot\dd B\cdot\eta\,, \text{ if } X\in\Gamma(TT^2) \\
     \nabla^{g_7}_X\eta=\nabla^{g_5}_X\chi_1\otimes \begin{pmatrix} 1 \\ 0 \end{pmatrix}+\nabla^{g_5}_X\chi_2\otimes \begin{pmatrix} 0 \\ 1 \end{pmatrix} \,, \text{ if } X\in\Gamma(TY_5)
\end{align}
\end{subequations}
Let us start with the components along $T^2$. Given that $\theta^s\wedge F_7=\theta^s\wedge F_5$, the equation becomes $\frac{1}{4}\theta^s\cdot(\dd B\cdot \eta+\ii\eta+\ee^B F_5\cdot\eta)=0$, which becomes
\begin{equation} \label{eq:GK7-torus-reduction-gravitino-along-torus}
    \dd B\cdot\chi_1+\ee^{B}F_5\cdot\chi_1+\ii\chi_1=0 \,, \quad -\dd B\cdot\chi_2+\ee^{B}F_5\cdot\chi_2+\ii\chi_2=0\,.
\end{equation}
Subtracting the equations for $\chi_2$ in~\eqref{eq:GK7-torus-reduction-dilatino} and~\eqref{eq:GK7-torus-reduction-gravitino-along-torus} yields $\chi_2=0$. For $\chi_1$, the equations coincide and are exactly~\eqref{eq:GK-KSE-rescaled-dilatino} with $n=2$. When $X\in\Gamma(TY_5)$, then $(X^\flat\wedge\vol_{T^2})\cdot=X^\flat\cdot\vol_{T^2}\cdot=\ii (X\cdot)\otimes \mathbf{1}$, so the remaining~\eqref{eq:GK-KSE-rescaled-gravitino} with $n=3$ becomes
\begin{equation}
    \nabla^{g_5}_X\chi_1+\frac{\ii}{4}X\cdot \chi_1-\frac{\ii}{4}X\cdot\chi_1+\frac{1}{4}\ee^B(X^\flat\wedge F_5)\cdot\chi_1= \nabla^{g_5}_X\chi_1+\frac{1}{4}\ee^{B}(X^\flat\wedge F_5)\cdot \chi_1 =0\,,
\end{equation}
which is exactly~\eqref{eq:GK-KSE-rescaled-gravitino} for $n=2$. Moreover, $\star_7 F_7=\ee^B\star_5 F_5\wedge \dd u_1\wedge \dd u_2-\ee^{-B}\vol_5$, therefore
\begin{equation}
    \dd\star_7 F_7=\dd(\ee^B\star_5 F_5)\wedge \dd u_1\wedge \dd u_2=0 \iff \dd(\ee^B\star_5 F_5)=0\,,
\end{equation}
which is~\eqref{eq:GK-flux-EOM-Bianchi-rescaled} for $n=2$.

After setting $n=2$ and $\kappa=-\frac{1}{2}$,~\eqref{eq:GK-KSE-rescaled} becomes exactly our~\eqref{eq:global-5d-GK-like-system}. 
For $\alpha=0$, \eqref{eq:global-alpha-deformed-GK-equation}
coincides with the GK flux equation~\eqref{eq:GK-flux-EOM-Bianchi-rescaled}. Thus, a solution of~\eqref{eq:global-5d-GK-like-system} on $Y_5$ with $\alpha=0$ and $\kappa=-\frac{1}{2}$ can be lifted to a solution of the Gauntlett--Kim system on $Y_5\times T^2$, which by~\cite{Gauntlett:2007ts} is a solution of Type IIB supergravity on $\AdS_3\times Y_5\times T^2$.

For non-zero circle curvature, the quadratic source
$-\alpha\wedge\alpha$ comes from the contribution
$\eta_1\wedge\pi^*\alpha$ to the closed six-dimensional torsion.
This is the transgression mechanism underlying the deformation,
analogous to that considered in \cite{Donos:2008ug}.

\section{Outlook}

The main result of this paper is a spinorial geometry formulation of  BHE geometry.
One of the advantages of our approach is that one does not have to require compactness and as a result, both compact and non-compact geometries  are described in a unified framework. We have shown that, from the spinorial point of view, this  system  is a very simple  deformation of the Type I supergravity supersymmetry equations. In this framework,  the $\mathbb{T}^2$-symmetry
 reduction and the holonomy reduction of the 
Bismut connection occur naturally. In the process, we  uncovered a handful of novel BHE structures,  that do not admit compact solutions, for which complete examples appear harder to construct.

 The problems discussed here  illustrate to
 a hopefully broad audience the versatility of spinorial geometry, widespread in the physics community, in the context of differential geometry. Thus, we hope that this paper can be useful in bridging the theoretical physics and differential geometry communities.

Our constructions bring under the umbrella of spinorial geometry BHE and GK structures, thus illuminating the relationship already noted in \cite{KS-26,ALS,ALL-26}, and allowing one to import techniques developed for one type of geometry, to the other type. For example, the 
extremization problem  associated to GK geometry \cite{Couzens:2018wnk,Gauntlett:2018dpc,Hosseini:2019ddy,Gauntlett:2019roi,Couzens:2022agr}  warrants a promising research avenue in the context of BHE.

We believe that there are several directions worth pursuing, both in the mathematics as well as in the physics sides.
For example, a more thorough exploration of the examples in the various subclasses of BHE geometries that we unearthed would be desirable, especially in the non-compact settings.
In this paper we focussed on BHE geometries  in six dimensions, but it should be straightforward to consider general $2n$-dimensional manifolds, that should be related to the 
GK Killing spinor   equations in dimension $2n-1$ \cite{Gauntlett:2002sc}. Another obvious extension of our analysis is to consider the generalization of our spinorial system  in various dimensions and with different 
structure groups, following the footsteps of \cite{Gauntlett:2003cy}.

The explicit embeddings of our GKSEs in different ten-dimensional supergravity 
theories that we provided (see Appendix \ref{sugraembeddings}), shed light onto the seemingly accidental  analogies between supersymmetric solutions that had been noticed before in the literature.

In particular, we have shown that the $\kappa \in \mathbb{R}$ equations arise in the context of AdS$_3$ solutions of massive Type IIA supergravity, and as a result we have identified a class of AdS$_3 \times I \times$ \BHER supersymmetric solutions, previously discussed in  \cite{Couzens:2022agr}. Alternatively, the same type of \BHER structure can be embedded in Type I supergravity, giving rise to AdS$_3 \times S^1 \times$ \BHER supersymmetric solutions, whose appearance was recently observed in \cite{Papadopoulos:2026xdj}. These two classes are in fact related by a concrete 1-1 map, that is a combination of T-duality and S-duality. 

Furthermore, we have pointed out that  the $\kappa \in \ii\mathbb{R}$ equations arise in the context of generalized Type I supergravity \cite{Arutyunov:2015mqj,Wulff:2016tju}, thus establishing that $\mathbb{R}^{1,3} \times$ \BHEI are supersymmetric solutions of this theory.
In the physics literature, generalized supergravity has mainly been studied in the context of integrable deformations of AdS string sigma models. Our results suggest that its supersymmetric solutions may give rise to a rich class of novel geometric structures. Its equations ensure scale invariance of the string worldsheet without imposing the usual Weyl-invariance conditions. Indeed, we found that \BHEI geometries must be non-compact, in agreement with the results of \cite{Papadopoulos:2024uvi,Papadopoulos:2024tgs}. 
It would be interesting to explore the implications of generalized supergravity, including the Type IIB version, for holography, and look for possible connections with recent work on scale covariance and generalized conformal symmetry \cite{Conti:2026xem,Bobev:2026kld}.

The various embeddings that we provided can be used in practice to construct new supersymmetric solutions in one set-up, from solutions obtained in a different set-up, or using examples constructed in the mathematics literature. For instance, the new compact examples of \BHER manifolds presented recently in  \cite{apostolov2026gluingconstructionzextremalkahler} give rise to novel supersymmetric AdS$_3$ solutions of (massive) Type IIA as well as Type I supergravities, that should be holographic dual to classes of two-dimensional supersymmetric conformal field theories.

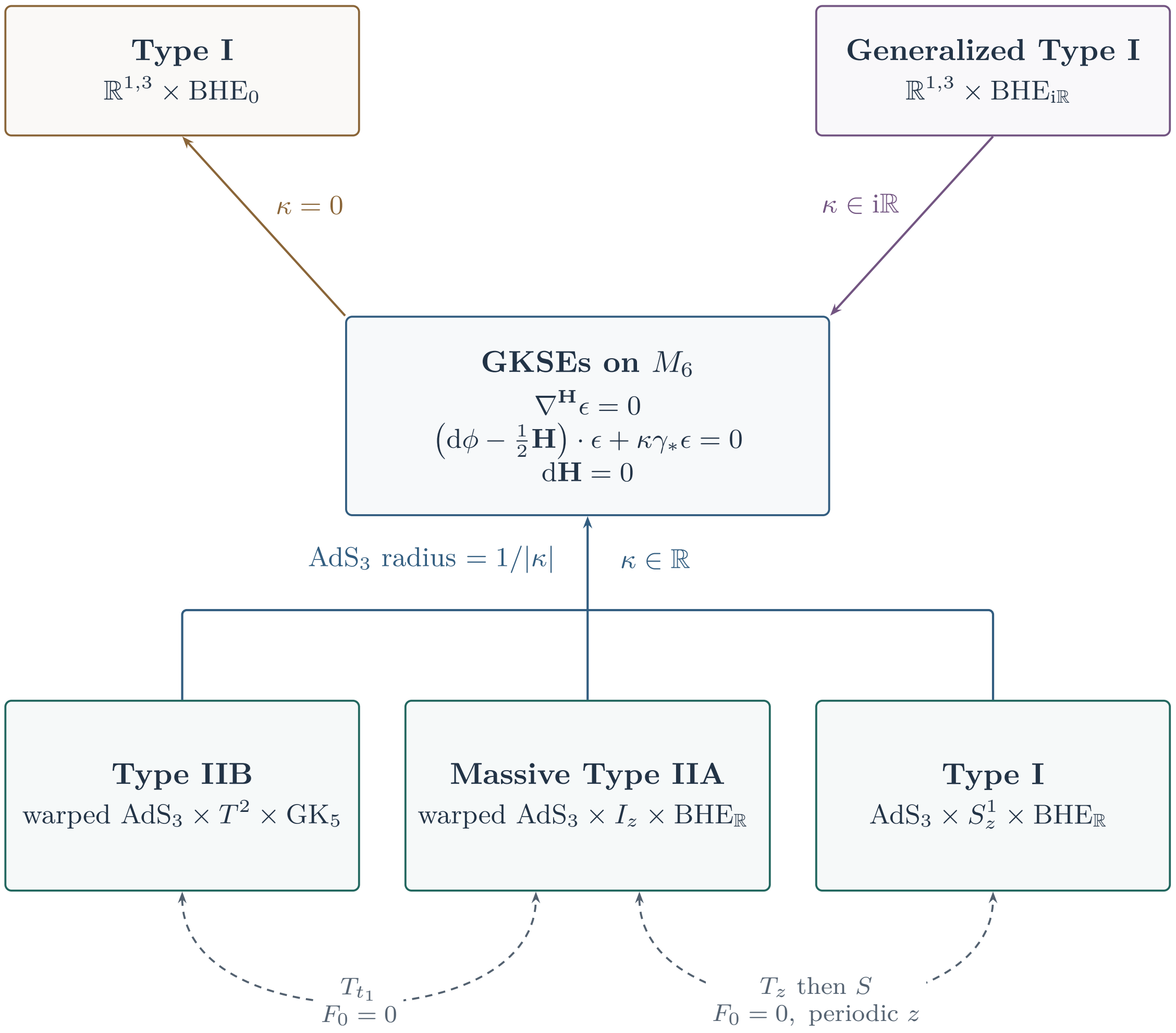
\begin{figure}[ht]
\centering
\captionsetup{font=normalsize}
\begingroup

\definecolor{Ink}{HTML}{223346}
\definecolor{Blue}{HTML}{325D80}
\definecolor{Plum}{HTML}{735582}
\definecolor{Teal}{HTML}{22675F}
\definecolor{Ochre}{HTML}{8A6538}
\definecolor{Grey}{HTML}{536170}

\tikzset{
  box/.style={rectangle, rounded corners=2pt, draw=#1, fill=#1!4!white, text=Ink,
    line width=.65pt, inner xsep=4pt, inner ysep=5pt, align=center, font=\footnotesize},
  upper/.style={minimum width=41mm, minimum height=15mm},
  theory/.style={box=Teal, minimum width=41mm, minimum height=22mm},
  arrow/.style={-{Stealth[length=4pt,width=3pt]}},
  embed/.style={draw=Blue, line width=.72pt, rounded corners=1.5pt},
  dual/.style={draw=Grey, line width=.65pt, dash pattern=on 2.4pt off 2pt,
    {Stealth[length=3.8pt,width=2.8pt]}-{Stealth[length=3.8pt,width=2.8pt]}},
  embed label/.style={font=\footnotesize, text=Blue, inner sep=2pt},
  dual label/.style={font=\scriptsize, text=Grey, fill=white, inner sep=2pt, align=center},
}

\newcommand{\head}[1]{{\small\bfseries #1}}
\newcommand{\sub}[1]{{\scriptsize\itshape #1}}
\newcommand{\GKAdS}{\mathrm{AdS}_3}

\noindent\makebox[\linewidth][c]{%
\resizebox{0.80\linewidth}{!}{%
\begin{tikzpicture}

\node[box=Ochre, upper] (flat) at (-4.7,4)
  {\head{Type I}\\[2pt] $\mathbb R^{1,3}\times\text{\BHEO}$};
\node[box=Plum, upper] (gen) at (4.7,4)
  {\head{Generalized Type I}\\[2pt] $\mathbb R^{1,3}\times\text{\BHEI}$};

\node[box=Blue, minimum width=56mm, minimum height=23mm] (GKSE) at (0,0)
  {\head{GKSEs on \(M_6\)}\\[3pt]
   $\nabla^\Ht\epsilon=0$\\
   $\big(\mathrm d\phi-\tfrac12\Ht\big)\cdot\epsilon+\kappa\gamma_*\epsilon=0$\\
   $\mathrm d\Ht=0$};

\node[theory] (IIB) at (-4.7,-4.4)
  {\head{Type IIB}\\ [2pt]
   $\text{warped }\GKAdS\times T^2\times\mathrm{GK}_5$};
\node[theory] (IIA) at (0,-4.4)
  {\head{Massive Type IIA}\\ [2pt]
   $\text{warped }\GKAdS\times I_z\times\text{\BHER}$};
\node[theory] (I) at (4.7,-4.4)
  {\head{Type I}\\ [2pt]
   $\GKAdS\times S^1_z\times\text{\BHER}$};

\draw[embed, arrow, draw=Ochre] (GKSE.north west) --
  node[midway, above right, font=\footnotesize, text=Ochre] {$\kappa=0$} (flat.south);
\draw[embed, arrow, draw=Plum] (gen.south) --
  node[midway, above left, font=\footnotesize, text=Plum] {$\kappa\in\ii\mathbb R$} (GKSE.north east);

\coordinate (join) at (0,-2.25);
\draw[embed] (IIB.north) |- (join);
\draw[embed] (I.north)   |- (join);
\draw[embed, arrow] (IIA.north) -- (GKSE.south);

\node[embed label, anchor=east] at ([xshift=-3mm,yshift=-5mm]GKSE.south)
  {$\GKAdS$ radius $=1/|\kappa|$};
\node[embed label, anchor=west] at ([xshift=3mm,yshift=-5mm]GKSE.south)
  {$\kappa\in\mathbb R$};

\draw[dual] (IIB.south) .. controls +(down:17mm) and +(down:17mm) ..
  node[midway, dual label] {$T_{t_1}$\\[-1pt] $F_0=0$}
  ([xshift=-6mm]IIA.south);
\draw[dual] ([xshift=6mm]IIA.south) .. controls +(down:17mm) and +(down:17mm) ..
  node[midway, dual label] {$T_z\ \text{then}\ S$\\[-1pt] $F_0=0,\ \text{periodic }z$}
  (I.south);

\end{tikzpicture}%
}}

\endgroup

\caption{The embedding from generalized Type I is discussed in
Section~\ref{gensugemb}.
The embedding from massive Type IIA is discussed in
Section~\ref{mass2aemb}.
In the massless case ($F_0=0$), after compactifying the $I_z$
interval to a circle $S^1_z$, a T-duality along $S^1_z$
followed by an S-duality leads to the Type I embedding in
the bottom right corner \cite{Papadopoulos:2026xdj}.
Instead, a T-duality along the $S^1_{t_1}$-circle
(that exists canonically in the \BHER{} geometry) leads to
the Type IIB setting in the bottom left corner.
Its relation to the six-dimensional GKSEs is understood
through this T-duality to the massless Type IIA embedding.
The five-dimensional ``transgressive'' GK manifold
\cite{Donos:2008ug} is obtained as
$\mathrm{GK}_5=M_6/S^1_{t_1}$,
$\pi\colon M_6\to\mathrm{GK}_5$ and
$\mathrm d\eta_1=\pi^*\alpha$.
When $\alpha=0$, $\mathrm{GK}_5\times T^2$ reduces to an
ordinary seven-dimensional GK-geometry.
The $\kappa=0$ arrow leads back to the \BHEO{} geometry
of Type I.}
\label{fig:web-of-embeddings}
\end{figure}

The web of embeddings and dualities is summarized in Figure \ref{fig:web-of-embeddings}, and it would be interesting to explore the dual field theory interpretation of these relationships.

\medskip

 \appendix

\section{Non Hermitian $\SU(3)$ structure} \label{sec:appendix}

In this section we consider non-generic solutions with $\kappa\neq0$ and
$\epsilon_-=\alpha\epsilon_+^*$, for some constant
$\alpha\in\mathbb C$. Since $\phi$ and $H$ are real and Clifford
multiplication by an odd form anticommutes with spinor conjugation,
conjugating the chiral dilatino equation for $\epsilon_-$ and using
$\epsilon_-=\alpha\epsilon_+^*$ gives
\begin{equation}
     \bar\alpha\left(\dd\phi-\frac12H\right)\cdot\epsilon_+
 =\bar\kappa\epsilon_+^*.
\end{equation}
The equation for $\epsilon_+$ therefore implies
$\kappa|\alpha|^2=\bar\kappa$.
Taking absolute values, and using $\kappa\neq0$, we obtain
$|\alpha|=1$ and hence $\kappa\in\mathbb R$. We may therefore write
$\alpha=\ee^{\ii\varphi}$, with $\varphi\in\mathbb R$. The gravitino
equation imposes no additional condition on this constant phase.
Define $
 \tilde\epsilon_\pm=\ee^{-\ii\varphi/2}\epsilon_\pm.$
Then $
 \tilde\epsilon_-
 =\ee^{-\ii\varphi/2}\ee^{\ii\varphi}\epsilon_+^*
 =\ee^{\ii\varphi/2}\epsilon_+^*
 =\tilde\epsilon_+^*.$
Dropping the tildes, the spinorial equations reduce to
\begin{align}
 \nabla^H\epsilon_+&=0,
 \label{eq:grav:appendix}\\
 \left(\dd\phi-\frac12H\right)\cdot\epsilon_+
 -\kappa\epsilon_+^*&=0,
 \label{eq:dil:appendix}
\end{align}
with $\kappa\in\mathbb R\setminus\{0\}$.
The spinor $\epsilon_+$ defines an $\SU(3)$ structure, which we denote
by $(g,J,\omega,\Omega)$\footnote{The $\SU(3)$ structure associated
to $\epsilon_-=\epsilon_+^*$ is $(g,-J,-\omega,-\bar\Omega)$.}.
Since its norm is constant, we may normalize $|\epsilon_+|^2=1$.

The mixed spinorial bilinear $\eta$ can still be defined without
assuming genericity. However, for every vector field $X$,
\begin{equation}
 \eta(X)
 :=\langle\epsilon_-,X^\flat\cdot\epsilon_+\rangle
 =\langle\epsilon_+^*,X^\flat\cdot\epsilon_+\rangle
 =\mathsf B(\epsilon_+,X^\flat\cdot\epsilon_+)=0,
\end{equation}
by \eqref{eq:spin-B-symmetry-kform} and the symmetry of $\mathsf B$.
Thus its metric-dual vector field $V$ also vanishes. The two spinors
have common stabilizer $\SU(3)$, so they do not induce an
$\SU(2)$ structure. Moreover,
$\mathsf B(\epsilon_+,\epsilon_-)=|\epsilon_+|^2=1$, and hence the
solution is not static. The integrability argument in
Proposition~\ref{prop:integrability-induced-complex-structures}
therefore does not apply. We will show below that $J$ is in fact
non-integrable whenever $\kappa\neq0$.

Theorem~\ref{thm:grks-generalized-ricci-soliton-kappa1} and
Proposition~\ref{prop:scalar-dilaton-equation} apply verbatim.
Thus $(g,H,f)$, with $f=2\phi$, satisfies the soliton equations
\eqref{eq:steady-gradient-generalized-ricci-soliton} and the scalar
identity \eqref{eq:scalar-dilaton-equation}, namely
\begin{equation}
 2|\operatorname{grad}\phi|^2-\Delta\phi
 -\frac1{12}|H|^2+2\kappa^2=0.
\end{equation}
As in Remark~\ref{rmk:gotheorem}, compact solutions are therefore
not excluded.

Theorem~\ref{wellknowtheorem} also applies, since it relies only on
the gravitino equation. In particular, $\nabla^H$ preserves
$(g,J,\omega,\Omega)$ and
\begin{equation}\label{eq:nongeneric-characteristic}
 N\in\Omega^3(M),
 \qquad H=\dd^{\mathrm{c}}\omega-N,
 \qquad \rho^H=0.
\end{equation}
Since we also require $\dd H=0$, the almost Hermitian structure is
both \emph{almost SKT} and \emph{almost CYT}, in the terminology
of Fino--Fowdar~\cite[Definitions~2.2 and~2.3]{Fino_2026}.
We will return below to their examples, in particular the
non-integrable structure on $S^3\times S^3$ in
\cite[Example~5.8]{Fino_2026}.

Let us first determine the Lee form $\theta$.
Since $\nabla^H J=0$, it is still given by the torsion trace
\[
 \theta(X)=-\frac12\sum_{a=1}^6 H(JX,e_a,Je_a),
\]
where $\{e_a\}$ is a local orthonormal frame. This identity does
not require the integrability of $J$ (see for instance \cite{Ivanov:2025zkd}). The  proof of Lemma~\ref{lem:Lee-forms} still applies and gives
\begin{equation}\label{eq:nongeneric-Lee}
 \theta=2\dd\phi,
 \qquad
 \dd\bigl(\ee^{-2\phi}\omega\wedge\omega\bigr)=0.
\end{equation}

Let us now compute $N$ using that $N=-4H^{(3,0)+(0,3)}$.
To determine this component of $H$, pair
\eqref{eq:dil:appendix} with $\epsilon_+$ using $\mathsf B$.
Since $
 \mathsf B(\epsilon_+,\dd\phi\cdot\epsilon_+)=0,$ and $
 \mathsf B(\epsilon_+,\epsilon_+^*)=1,$
we find $
 -\frac12\mathsf B(\epsilon_+,H\cdot\epsilon_+)=\kappa.$
As in the proof of
Proposition~\ref{prop:integrability-induced-complex-structures},
the definition of $\Omega$ then yields
\begin{equation}
     g_{\mathbb C}(H,\Omega)
 =-6\mathsf B(\epsilon_+,H\cdot\epsilon_+)=12\kappa.
\end{equation}
Since $g_{\mathbb C}(\Omega,\bar\Omega)=48$,
$g_{\mathbb C}(\Omega,\Omega)=0$ and $H$ is real, 
\begin{equation}\label{eq:nongeneric-H-pure}
 H^{(3,0)+(0,3)}
 =\frac{\kappa}{4}(\Omega+\bar\Omega)
 =\frac{\kappa}{2}\operatorname{Re}\Omega.
\end{equation}
and therefore
\begin{equation}\label{eq:nongeneric-N}
 N=-2\kappa\operatorname{Re}\Omega.
\end{equation}
In particular, $J$ is in general not integrable. Substituting into
\eqref{eq:nongeneric-characteristic} gives
\begin{equation}\label{nonintegrableJ}
 H=\dd^{\mathrm{c}}\omega+2\kappa\operatorname{Re}\Omega.
\end{equation}
This agrees with Fino--Fowdar~\cite[Proposition~4.1]{Fino_2026},
after translating the conventions for the structure forms.

We next compute $\dd\Omega$. Since $\Omega$ is $\nabla^H$-parallel,
\begin{equation}\label{eq:nongeneric-dOmega-torsion}
 \dd\Omega=\sum_{a=1}^6(\iota_{e_a}H)\wedge(\iota_{e_a}\Omega).
\end{equation}
The $\SU(3)$ identities give that $
 \dd\Omega=\theta\wedge\Omega+\kappa\,\omega\wedge\omega.$
Using $\theta=2\dd\phi$, we conclude that
\begin{equation}\label{eq:nongeneric-weighted-Omega}
 \dd\bigl(\ee^{-2\phi}\Omega\bigr)
 =\kappa\, \ee^{-2\phi}\omega\wedge\omega.
\end{equation}

In conclusion, the $\SU(3)$ structure satisfies
\begin{subequations}\label{eq:gmw:modified}
\begin{align}
 \dd\bigl(\ee^{-2\phi}\Omega\bigr)
 &=\kappa\, \ee^{-2\phi}\omega\wedge\omega,
 \label{eq:gmw-a-mod}\\
 \dd\bigl(\ee^{-2\phi}\omega\wedge\omega\bigr)&=0,
 \label{eq:gmw-b-mod}\\
 H&=\dd^{\mathrm{c}}\omega+2\kappa\operatorname{Re}\Omega.
 \label{eq:gmw-c-mod}
\end{align}
\end{subequations}
These equations give a non-Hermitian deformation of the Type~I
conditions in \cite{str,Hull,Gauntlett:2003cy}. It can be shown that \cite[Example~5.8]{Fino_2026} solves the equations \eqref{eq:gmw:modified}, together with the Bianchi identity $dH=0$ with $\kappa=2$ and $\theta=0$, that is, the dilaton is constant.
Conversely, they recover the spinorial system. 

Indeed, let $(g,J,\omega,\Omega)$ be an $\SU(3)$ structure satisfying $\dd \dd^{\mathrm{c}} \omega+2\kappa \, \dd\operatorname{Re}\Omega=0$ and 
\eqref{eq:gmw:modified}, for some smooth function $\phi$ and
constant $\kappa\in\mathbb R\setminus\{0\}$ and let $\epsilon_+$ be its associated unit spinor.  We show below that it satisfies
\eqref{eq:grav:appendix}--\eqref{eq:dil:appendix} with $dH=0$.

By \eqref{eq:gmw-a-mod}, the intrinsic torsion component
obstructing the skew-symmetry of the Nijenhuis tensor
vanishes; see \cite[Proposition~4.1]{Fino_2026}.
Thus there exists a unique Hermitian connection with
skew torsion $\hat H=\dd^{\mathrm{c}}\omega-N$.
Moreover, \eqref{eq:gmw-b-mod} gives $\theta=2\dd\phi$.

As in Section~\ref{sec:sufficiency}, write
$\nabla^{\hat H}\Omega=\ii\alpha\otimes\Omega$, with $\alpha$
a real one-form. The same torsion calculation, now
including the $(3,0)+(0,3)$-component of $\hat H$, gives
\begin{equation}
  \dd\Omega
 =(\ii\alpha+\theta)\wedge\Omega
 +\frac1{12}g_{\mathbb C}(\hat H,\Omega)\,\omega\wedge\omega.   
\end{equation}
Comparing types with \eqref{eq:gmw-a-mod}, we obtain $
 \alpha\wedge\Omega=0,$ and $
 g_{\mathbb C}(\hat H,\Omega)=12\kappa.$
Since $\alpha$ is real, the first identity implies
$\alpha=0$, exactly as in Section~\ref{sec:sufficiency}.
The second identity, together with the fact that $\hat H$ is real,
gives
\[
 \hat H^{(3,0)+(0,3)}
 =\frac{\kappa}{2}\operatorname{Re}\Omega,
 \qquad
 N=-4\hat H^{(3,0)+(0,3)}
 =-2\kappa\operatorname{Re}\Omega.
\]
Consequently,
$\hat H=\dd^{\mathrm{c}}\omega+2\kappa\operatorname{Re}\Omega=H$,
and $\nabla^H$ preserves the defining spinor $\epsilon_+$.

Finally, the spinorial decomposition used in
Section~\ref{sec:sufficiency} gives
\begin{equation}
    H\cdot\epsilon_+
 =-\frac16g_{\mathbb C}(H,\Omega)\epsilon_+^*
   +\theta\cdot\epsilon_+
 =-2\kappa\epsilon_+^*+2\dd\phi\cdot\epsilon_+.
\end{equation}

Indeed, the pure-type component contributes only to the
coefficient of $\epsilon_+^*$, while the computation of
the real one-form component is unchanged.
This proves \eqref{eq:dil:appendix}.

\medskip
We conclude this appendix by showing that non-generic solutions can be constructed from the static solutions of the GKSE with real $\kappa$ considered in the main body of the paper. Let $\epsilon=\epsilon_++\epsilon_-$ be a static solution of the GKSE, 
with $\kappa\in\mathbb R\setminus\{0\}$ and
$|\epsilon_\pm|=1$, and let
$(\eta,\omega^T,\Omega^T)$ be its associated
$\SU(2)$ structure, with $\eta=\eta_1+\ii\eta_2$.
For every $z\in U(1)$, consider the spinor $\epsilon_z:=\epsilon+z\epsilon^*$. We show that it
is a non-generic solution of the GKSE with the same data $(g,H,\phi,\kappa)$.

Since $\kappa$ is real, $\epsilon^*$ satisfies the same
spinorial equations as $\epsilon$, as shown in the proof
of Proposition~\ref{prop:beyond-static}.
By complex linearity, $\epsilon_z$ is therefore a solution.
Since $|z|=1$, its chiral components satisfy $
 (\epsilon_z)_-
 =\epsilon_-+z\epsilon_+^*
 =z\bigl((\epsilon_z)_+\bigr)^*.$
Moreover, staticity implies
$\langle \epsilon_+,\epsilon_-^*\rangle=0$, so
$|(\epsilon_z)_+|^2=2$.
Thus $\epsilon_z$ is non-zero and non-generic. Writing $z=e^{\ii\vartheta}$, we replace $\epsilon_z$ by
$2^{-1/2}e^{-\ii\vartheta/2}\epsilon_z$,
so that its chiral components have unit norm and satisfy
$(\epsilon_z)_-=((\epsilon_z)_+)^*$. The $\SU(3)$ structure associated with the spinor  $\epsilon_z$ is different from the original  $\SU(3)$ structures, while $(g,H,\phi,\kappa)$ and the connection $\nabla^H$ are independent of $z$ and remain
fixed. In particular, in terms of the $\SU(2)$ structure
$(\eta,\omega^T,\Omega^T)$ 
this reads 
\begin{equation}
  \omega_z = \eta_1 \wedge \eta_2 - \operatorname{Im}\!\left( \bar{z}\, \Omega^T \right),
  \qquad
  \Omega_z = -\eta \wedge \left( \operatorname{Re}\!\left( \bar{z}\, \Omega^T \right) + i\, \omega^T \right),
\end{equation}
 and  it satisfies
\eqref{eq:gmw:modified}, with
$\theta_z=2\dd\phi$ and
$N_z=-2\kappa\operatorname{Re}\Omega_z\neq0$.
Its characteristic connection has holonomy contained
in $\SU(2)$, and is flat precisely when the connection
of the original static solution is flat.

For example, the cohomogeneity-one family on
$S^3\times S^3$ constructed\footnote{See also \cite{podesta2026newinhomogeneousbismuthermitianeinstein}.} by \cite{Donos:2008hd} (based on \cite[Appendix C]{Donos:2008ug}), which depends on two parameters $(a,\varepsilon)$
 admits a static $\SU(2)$ structure that solves the GKSE, as recalled in Section \ref{section:examples}. Therefore, non-generic partner solutions, depending on $z\in U(1)$ exist. 
For $\varepsilon\neq0$, taking $z=1$ gives an almost CYT
solution with non-constant dilaton and non-flat
characteristic connection.
For $\varepsilon=0$, the dilaton becomes constant and the
background becomes homogeneous and Bismut-flat.
In this limit, for $a=\pm1$ and after a homothety fixing $\kappa=2$,
the non-generic construction recovers the structure of~\cite[Example~5.8]{Fino_2026}.

\section{Geometry with unequal chiral norms}
\label{app:beyond-static}

In Section~\ref{sec:su2su3}, we assumed that the static spinor
$\epsilon=\epsilon_++\epsilon_-$ satisfies
$|\epsilon_+|=|\epsilon_-|=1$. In this appendix we explore
the case when the structure is static but the norms of
$\epsilon_\pm$ are not equal. Throughout, we assume
$\kappa\neq0$. By Proposition~\ref{prop:beyond-static},
every generic solution admits a static representative
with the same data $(g,H,\phi,\kappa)$.

The purpose of this appendix is to show that unequal chiral
norms still lead to BHE structures. However, for real $\kappa$,
they may prevent the transverse Hermitian structure from
being locally conformally K\"ahler.
\medskip

Since $\epsilon_\pm$ have constant norm by
Proposition~\ref{prop:beyond-static}, we explore the geometry
induced by the normalized spinors $
 \widehat\epsilon_\pm
 :=\frac{\epsilon_\pm}{|\epsilon_\pm|}.$
Although both spinors now have unit norm, this separate
normalization introduces two different coefficients in
the dilatino equations. Set $
 \kappa_\pm
 :=\frac{\kappa|\epsilon_\mp|}{|\epsilon_\pm|}.
$
The normalized spinors satisfy
\begin{align}
 \nabla^H\widehat\epsilon_\pm&=0,
 \label{eq:grav-normalized-unequal}\\
 \left(\dd\phi-\frac12H\right)\cdot\widehat\epsilon_+
 &=\kappa_+\widehat\epsilon_-,
 \label{eq:dil-normalized-unequal-plus}\\
 \left(\dd\phi-\frac12H\right)\cdot\widehat\epsilon_-
 &=-\kappa_-\widehat\epsilon_+.
 \label{eq:dil-normalized-unequal-minus}
\end{align}
In particular, $\kappa_+\kappa_-=\kappa^2$, while
$\kappa_+=\kappa_-$ if and only if
$|\epsilon_+|=|\epsilon_-|$.
Thus, when the norms are unequal, the normalized spinor
does not satisfy the original system with a single parameter.

Since $\widehat\epsilon_\pm$ are unitary and
$\nabla^H$-parallel, they define two $\nabla^H$-parallel
$\SU(3)$ structures $(g,J^\pm,\omega^\pm,\Omega^\pm)$,
which are almost CYT by Theorem~\ref{wellknowtheorem}.
The static condition is preserved by normalization.
The same pairing argument as in
Proposition~\ref{prop:integrability-induced-complex-structures}
therefore gives $
 \mathsf B(\widehat\epsilon_\pm,
           H\cdot\widehat\epsilon_\pm)=0.$
Indeed, the terms involving $\dd\phi$ vanish by
\eqref{eq:spin-B-symmetry-kform}, and those involving
$\kappa_\pm$ vanish by the static condition.
Consequently, both $J^\pm$ are integrable and the two
induced structures are BHE.

The normalized spinor
$\widehat\epsilon=\widehat\epsilon_++\widehat\epsilon_-$
also induces an $\SU(2)$ structure. Let $
 \eta(X)
 :=\langle\widehat\epsilon_-,X^\flat\cdot
          \widehat\epsilon_+\rangle,
$ and $
 \eta=\eta_1+\ii\eta_2,$
and denote its metric dual by $V=V_1+\ii V_2$, where
$V_i=\eta_i^\sharp$.
The fields $V_1,V_2$ are orthonormal and
$\nabla^H$-parallel. Hence they are Killing and $
 \dd\eta_i=\iota_{V_i}H.$

The computation of Lemma~\ref{lem:Lee-forms} applies to
\eqref{eq:dil-normalized-unequal-plus} and
\eqref{eq:dil-normalized-unequal-minus}, giving
\begin{equation}\label{eq:Lee-forms-unequal}
 \theta^+=2\dd\phi-2\operatorname{Re}(\kappa_+\bar\eta),
 \qquad
 \theta^-=2\dd\phi+2\operatorname{Re}(\kappa_-\eta).
\end{equation}
Let $\mathcal H=\ker\eta_1\cap\ker\eta_2$.
The difference $\theta^+-\theta^-$ is a linear combination
of $\eta_1,\eta_2$, and therefore vanishes on $\mathcal H$.

With our convention of summing over all indices, the
torsion trace identity reads $
 \frac12\operatorname{tr}_{\omega^\pm}H=J^\pm\theta^\pm
$.
Using $\omega^\pm=\omega^T\pm\eta_1\wedge\eta_2$,
we obtain, for every horizontal vector field $X$,
\begin{equation}
   2H(X,V_1,V_2)
 =\frac12
   \bigl(\operatorname{tr}_{\omega^+}H
        -\operatorname{tr}_{\omega^-}H\bigr)(X)
 =\bigl(J^+\theta^+-J^-\theta^-\bigr)(X)
 =0.  
\end{equation}
Here we used that $J^+$ and $J^-$ agree on $\mathcal H$
and preserve this distribution.
By skew-symmetry, $H(V_1,V_2,\cdot)$ also vanishes on
the vertical distribution. Thus $
 H(V_1,V_2,\cdot)=0.$
It follows that both $\dd\eta_i$ are horizontal and $
 [V_1,V_2]=-H(V_1,V_2,\cdot)^\sharp=0.$
Consequently, the torsion decomposes as
\begin{equation}\label{eqn:H2}
 H=\eta_1\wedge \dd\eta_1+\eta_2\wedge \dd\eta_2+H^T,
\end{equation}
where $H^T$ is the torsion of the transverse Hermitian structure.

Equation \eqref{eqn:H2}, together with $H=\dd^{\mathrm{c}}_\pm \omega_\pm$, shows that $\dd\eta_i$ have type $(1,1)$.
Together with $\nabla^H J^\pm=0$ and
$\nabla^H V_i=0$, this gives
$\mathcal L_{V_i}J^\pm=0$.
Thus the transverse metric and complex structure are basic
and induce a Hermitian structure on each local leaf space.
Moreover, Cartan's formula gives
$\mathcal L_{V_j}\dd\eta_i=0$, so the curvature forms
$\dd\eta_i$ are also basic.
\medskip

The computation of
Lemma~\ref{lem:dilaton-vertical-derivatives}, applied to the
normalized equations, gives
\begin{equation}\label{eq:vertical-dilaton-unequal}
 2V_1(\phi)=\operatorname{Re}(\kappa_+-\kappa_-),
 \qquad
 2V_2(\phi)=\operatorname{Im}(\kappa_++\kappa_-).
\end{equation}
In particular, these derivatives are constant.
The transverse Lee form is
\begin{equation}\label{eq:Lee-transverse-unequal}
 \begin{aligned}
 \theta^T
 &=2\dd\phi-2V_1(\phi)\eta_1-2V_2(\phi)\eta_2\\
 &=2\dd\phi
   -\operatorname{Re}(\kappa_+-\kappa_-)\eta_1
   -\operatorname{Im}(\kappa_++\kappa_-)\eta_2,
 \end{aligned}
\end{equation}
and satisfies $\dd\omega^T=\theta^T\wedge\omega^T$.
Indeed, this follows by averaging the identities for
$\dd\omega^\pm$, exactly as in the unit-norm case.
The expression in \eqref{eq:Lee-transverse-unequal}
is horizontal by \eqref{eq:vertical-dilaton-unequal}
and invariant under $V_1,V_2$, hence basic.

The argument used in
\eqref{eq:weighted-holomorphic-forms}--%
\eqref{eq:transverse-holomorphic-form} similarly gives
\begin{equation}\label{eq:omegaT}
 \dd\Omega^T=
 \Bigl(
 \theta^T
 -\ii\operatorname{Im}(\kappa_+-\kappa_-)\eta_1
 +\ii\operatorname{Re}(\kappa_++\kappa_-)\eta_2
 \Bigr)\wedge\Omega^T.
\end{equation}
Here $\Omega^T$ is horizontal, but need not be basic.

Since the transverse geometry has real dimension four,
$H^T=-\star_T\theta^T$.
The torsion decomposition \eqref{eqn:H2} therefore shows
that the Bianchi identity is equivalent to
\begin{equation}\label{eq:bianchi}
 -\dd\star_T\theta^T
 +\dd\eta_1\wedge \dd\eta_1
 +\dd\eta_2\wedge \dd\eta_2=0.
\end{equation}
Comparing the vertical components of
$\frac12\operatorname{tr}_{\omega^\pm}H
 =J^\pm\theta^\pm$, using
\eqref{eq:Lee-forms-unequal} and
\eqref{eq:Lee-transverse-unequal}, gives
\begin{equation}\label{eq:traces}
 \frac12\operatorname{tr}_{\omega^T}(\dd\eta_1)
 =-\operatorname{Im}(\kappa_+-\kappa_-),
 \qquad
 \frac12\operatorname{tr}_{\omega^T}(\dd\eta_2)
 =\operatorname{Re}(\kappa_++\kappa_-).
\end{equation}
Taking the transverse Hodge star of \eqref{eq:bianchi}
and using \eqref{eq:traces}, we obtain
\begin{equation}\label{eq:box-unequal-norms}
 \dd_T^*\theta^T
 +\bigl|\kappa_++\overline{\kappa_-}\bigr|^2
 -\frac12|\dd\eta_1|_{g^T}^2
 -\frac12|\dd\eta_2|_{g^T}^2=0.
\end{equation}
This is the transverse scalar equation in the general
Hermitian setting. Note that applying \cite{GGP}, the Bismut Ricci form of the transverse metric is 
\begin{equation}
    \rho^{B,T}=-\operatorname{Im}(\kappa_+-\kappa_-) \dd \eta_1 + \operatorname{Re}(\kappa_++\kappa_-) \dd \eta_2.
\end{equation}
This is useful for the reverse construction analogue to Theorems \ref{thm:real-kappa-local-model} and \ref{thm:imaginary-kappa-local-model}.

\medskip

If $M$ is compact, $\phi$ attains a maximum, where
$V_1(\phi)=V_2(\phi)=0$.
Since these derivatives are constant by
\eqref{eq:vertical-dilaton-unequal}, they vanish identically.
For real $\kappa$, this implies $\kappa_+=\kappa_-$,
and hence $|\epsilon_+|=|\epsilon_-|$.
For purely imaginary $\kappa\neq0$, the identity
$2V_2(\phi)=\operatorname{Im}(\kappa_++\kappa_-)$
gives a contradiction, since $\kappa_\pm$ are positive
real multiples of $\kappa$.
Thus compact solutions necessarily have real $\kappa$
and equal chiral norms.

\medskip

Suppose now that $\kappa\in\mathbb R$.
Then \eqref{eq:vertical-dilaton-unequal} gives $
 V_2(\phi)=0,
 $ and $
 V_1(\phi)=\frac12(\kappa_+-\kappa_-),$
and \eqref{eq:Lee-transverse-unequal} reduces to
\begin{equation}\label{eq:Lee-unequal}
 \theta^T=2\dd\phi-2V_1(\phi)\eta_1.
\end{equation}
Since $V_1(\phi)$ is constant, $
 \dd\theta^T=-2V_1(\phi)\dd\eta_1.$
In particular, when the norms of $\epsilon_\pm$ are unequal,
$V_1(\phi)\neq0$, and the transverse structure is locally
conformally K\"ahler if and only if $\dd\eta_1=0$.

\medskip

We now discuss sufficiency. The starting point is analogous
to that of Section~\ref{sec:sufficiency}, with
\eqref{eq:modified-Hull-Strominger} unchanged and
$\dd H=0$, while
\eqref{eq:sufficiency-Leep-form} and
\eqref{eq:sufficiency-Leem-form} are replaced by
\begin{subequations}
\label{eq:sufficiency-Lee-forms-general-case}
\begin{align}
 \leep &=2\dd\phi-2\mathrm{Re}(\kappa_+\vformbar),\\
 \leem &=2\dd\phi+2\mathrm{Re}(\kappa_-\vform),
\end{align}
\end{subequations}
where
$\kappa_\pm=\kappa a_\mp/a_\pm$
for positive constants $a_\pm$.

Let $\widehat\epsilon_\pm$ be the compatible unit spinors
defining the $\SU(2)$ structure.
The dilatino argument proceeds as in
Section~\ref{sec:sufficiency} up to
\eqref{eq:gravitino-converse-lee-form}.
Substituting
\eqref{eq:sufficiency-Lee-forms-general-case} then gives
\begin{equation}
 \left(\dd\phi-\frac12H\right)\cdot\widehat\epsilon_\pm
 \mp\kappa_\pm\widehat\epsilon_\mp=0.
\end{equation}
The gravitino argument is unchanged and yields
\begin{equation}
 \nabla^H\widehat\epsilon_\pm=0.
\end{equation}
Finally, set $\epsilon_\pm=a_\pm\widehat\epsilon_\pm$.
Since $a_\pm$ are constant and
$a_\pm\kappa_\pm=\kappa a_\mp$, we obtain
\begin{equation}
 \nabla^H\epsilon_\pm=0,
 \qquad
 \left(\dd\phi-\frac12H\right)\cdot\epsilon_\pm
 \mp\kappa\epsilon_\mp=0.
\end{equation}
Thus $\epsilon=\epsilon_++\epsilon_-$ solves the original
GKSE, with $|\epsilon_\pm|=a_\pm$.

\medskip

\medskip

Appendices A and B conclude our analysis of the possible solutions of GKSE. We summarize their geometry in the following table.


\begin{table}[ht]
\centering
\begingroup
\renewcommand{\arraystretch}{1.8}
\setlength{\tabcolsep}{10pt}

\resizebox{0.9\textwidth}{!}{%
\begin{tabular}{|l|c|c|}
\hline
 & $\kappa\in\mathbb R\setminus\{0\}$
 & $\kappa=\ii\lambda\in\ii\mathbb R\setminus\{0\}$ \\
\hline

\multicolumn{3}{|c|}{\textbf{Non-generic}}\\
\hline

$\SU(3)$ structure
&
$\dd(\ee^{-2\phi}\Omega)
 =\kappa\ee^{-2\phi}\omega\wedge\omega$
&
--- \\

&
$\dd(\ee^{-2\phi}\omega\wedge\omega)=0$
&
--- \\

Torsion
&
$H=\dd^c\omega+2\kappa\operatorname{Re}\Omega$
&
--- \\

$\dd H=0$
&
$\dd\dd^c\omega=-2\kappa\,\dd\operatorname{Re}\Omega$
&
--- \\

Compact $M_6$
& allowed & excluded \\

Non-compact $M_6$
& allowed & excluded \\

\hline
\multicolumn{3}{|c|}{\textbf{Generic, unequal chiral norms}}\\
\hline

Transverse $\omega^T$
&
$\dd\omega^T=\theta^T\wedge\omega^T$
&
$\dd\omega^T=\theta^T\wedge\omega^T$
\\

Transverse $\Omega^T$
&
$\dd\Omega^T=
 \bigl(\theta^T+\ii(\kappa_++\kappa_-)\eta_2\bigr)
 \wedge\Omega^T$
&
$\dd\Omega^T=
 \bigl(\theta^T-\ii(\lambda_+-\lambda_-)\eta_1\bigr)
 \wedge\Omega^T$
\\

Transverse Lee form
&
$\theta^T=2\dd\phi-(\kappa_+-\kappa_-)\eta_1$
&
$\theta^T=2\dd\phi-(\lambda_++\lambda_-)\eta_2$
\\

Curvature type
&
$\dd\eta_i\in\Omega^{1,1}$
&
$\dd\eta_i\in\Omega^{1,1}$
\\

Curvature traces
&
$\displaystyle
 \frac12\operatorname{tr}_{\omega^T}(\dd\eta_1)=0$
&
$\displaystyle
 \frac12\operatorname{tr}_{\omega^T}(\dd\eta_1)
 =-(\lambda_+-\lambda_-)$
\\

&
$\displaystyle
 \frac12\operatorname{tr}_{\omega^T}(\dd\eta_2)
 =\kappa_++\kappa_-$
&
$\displaystyle
 \frac12\operatorname{tr}_{\omega^T}(\dd\eta_2)=0$
\\

$\dd H=0$
&
$\dd_T^*\theta^T+(\kappa_++\kappa_-)^2
-\frac12|\dd\eta_1|_{g^T}^2
-\frac12|\dd\eta_2|_{g^T}^2=0$
&
$\dd_T^*\theta^T+(\lambda_+-\lambda_-)^2
-\frac12|\dd\eta_1|_{g^T}^2
-\frac12|\dd\eta_2|_{g^T}^2=0$
\\

Compact $M_6$
& excluded & excluded \\

Non-compact $M_6$
& allowed & allowed \\

\hline
\end{tabular}%
}
\endgroup

\caption{Summary of the non-generic and unequal-norm cases
of the GKSE with $\kappa\neq0$.
The generic case with equal chiral norms is summarized
in Table~1.
By Proposition~\ref{prop:beyond-static}, every generic
solution admits a static representative with the same
$(g,H,\phi,\kappa)$; the chiral norms in this table refer
to this representative.
In the non-generic case, the induced almost complex
structure is not integrable.
Here $\kappa_\pm=\kappa|\epsilon_\mp|/|\epsilon_\pm|$ and,
in the imaginary case, $\kappa_\pm=\ii\lambda_\pm$.
The curvature forms have transverse type $(1,1)$.}
\label{tab:spinorial-cases}
\end{table}

\section{A guide for physicists}

Despite being standard in differential geometry, the notations we have used throughout the paper are less common in the theoretical/mathematical physics literature. The aim of this appendix is to provide a handy dictionary to convert the results in this paper to a notation more familiar to physicists, and to summarize the main results in a way that is simpler to reference in the physics literature. 

\subsection{Conventions}
We use the basis on $\mathrm{Cliff}(6)$ where the generators are purely imaginary and antisymmetric. As a matter of notation, we will use this dictionary (\cf\eqref{eq:spin-complex-bilinear-form})
\begin{subequations}
\begin{align}
	\langle\eps,\gamma_{m_1\ldots m_k}\psi\rangle &=\eps^\dagger\gamma_{m_1\ldots m_k}\psi, \\
	\mathsf B(\eps,\gamma_{m_1\ldots m_k}\psi) &= (\eps^c)^\dagger\gamma_{m_1\ldots m_k}\psi=\eps^T\gamma_{m_1\ldots m_k}\psi. 
\end{align}
\end{subequations}
From \eqref{eq:spin-chirality-local}, we recall
\begin{equation}
    \gamma_*=-\ii \gamma^{123456}.
\end{equation}
Moreover, in this appendix, and in the next, we adopt the following convention on the indices, which is common in the supergravity literature: latin letters $m,n,p,\ldots$ denote ``curved'' indices, while $a,b,c,\ldots$ denote ``flat'' indices. Upper-case is reserved for the ten-dimensional fields, while lower-case for the six-dimensional internal geometries. The ``curved'' and ``flat'' indices are related by multiplication via the vielbein. Namely, writing the orthonormal frame in the coordinate basis $e^a=e^a_m\dd x^m$, we call $e^a_m$ the vielbein, and we find the ``curved object'' via multiplication with the vielbein. For example,
\begin{equation}
\gamma_m=e^a_m \gamma_a,
\end{equation}
so that 
\begin{equation}
    \{\gamma_a,\gamma_b\}=2\delta_{ab}\mathbf{1}\iff \{\gamma_m,\gamma_n\}=2g_{mn}\mathbf{1}.
\end{equation}

\subsection{Main results} We now summarize the main results of the paper in the conventions and notation just introduced.

\subsubsection{Deformed KSE system} The deformed Killing spinor equations are (\cf \eqref{eq:grav}, \eqref{eq:dil})
\begin{subequations}\label{eq:app-KSE-modified}
	\begin{align}
		&(\mathcal{A}+\kappa\gamma_*)\eps\coloneqq\left(\gamma^m\nabla_m \phi-\frac{1}{12}H_{mnp}\gamma^{mnp}+\kappa\gamma_*\right)\eps=0\,, \label{eq:app-KSE-modified-1} \\  
		&\mathcal{D}\eps\coloneqq\left(\nabla_m
		-\frac{1}{8}H_{mnp}\gamma^{np}\right)\eps=0\,. \label{eq:app-KSE-modified-2}
	\end{align}
\end{subequations} 
Decomposing the Killing spinor into its chiral components $\eps=\epsp+\epsm$, the system decouples into (\cf \eqref{eq:KSE-chiral-gravitino}, \eqref{eq:KSE-chiral-dilatino})
\begin{subequations}\label{eq:app-KSE-modified-decoupled}
	\begin{align}
		\mathcal{A}\epspm\mp\kappa\epsmp&=(\gamma^a\nabla_a \phi)\epspm-\frac{1}{12}H_{abc}\gamma^{abc}\epspm\mp \kappa\epsmp=0\,,\label{eq:app-KSE-modified-decoupled-1} \\  
		\mathcal{D}\epspm&=\left(\nabla_m
		-\frac{1}{8}H_{mab}\gamma^{ab}\right)\epspm=0\,. \label{eq:app-KSE-modified-decoupled-2}
	\end{align}
\end{subequations} 
\subsubsection{Integrability conditions} From the integrability conditions on the KSE~\eqref{eq:app-KSE-modified}, we obtain (\cf\eqref{eq:steady-gradient-generalized-ricci-soliton}, \eqref{eq:scalar-dilaton-equation})
\begin{subequations} \label{eq:app-EOM}
	\begin{align}
		&R_{mn}-\frac{1}{4}H_{mab}H_n{}^{ab}+2\nabla_m\nabla_n\phi=0\,, \\
		&-\nabla^2\phi+2(\nabla\phi)^2-\frac{1}{12}H_{abc}H^{abc}+2\kappa^2=0\,,
	\end{align}
\end{subequations}
and the equation of motion for $H$
\begin{equation} \label{eq:app-H-flux-EOM}
	\nabla_aH^a{}_{mn}-2\nabla_a\phi H^a{}_{mn}=0 \quad (\text{or equivalently } \dd (\ee^{-2\phi} \star H) =0) \,.
\end{equation}
In the limit $\kappa\to 0$,~\eqref{eq:app-EOM} and~\eqref{eq:app-H-flux-EOM} become (A.4) of~\cite{Gauntlett:2003cy}\footnote{Note that the conventions in~\cite{Gauntlett:2003cy} are $H_\mathrm{there}=-H_\mathrm{here}$.} when $F=0$.
\subsubsection{$\SU(3)$ and $\SU(2)$ structures} With a static $\SU(2)$ structure $(\epsp, \epsm)$ we introduce the bilinears (\cf\eqref{eq:spin-eta-V-local}, \eqref{eq:spin-omegaT-local}, \eqref{eq:spin-OmegaT-local})
\begin{subequations}
    \begin{align}
 \eta_a&=\epsilon_-^\dagger\gamma_a\epsilon_+,
 \quad
 V^a=\epsilon_-^\dagger\gamma^a\epsilon_+,\\
 \omega^T_{ab}
 &=-\frac{\ii}{2}\left(
 \epsilon_+^\dagger\gamma_{ab}\epsilon_+
 +\epsilon_-^\dagger\gamma_{ab}\epsilon_-
 \right),\\
 \Omega^T_{ab}
 &=\epsilon_+^T\gamma_{ab}\epsilon_-.
    \end{align}
\end{subequations}
The $\SU(3)^\pm$-structures can be written in terms of the transverse $\SU(2)$ structure with (\cf \eqref{eq:spin-SU2-to-SU3-omega-global},~\eqref{eq:spin-SU2-to-SU3-Omega-global})
\begin{align}
 \omega^\pm
 &=\omega^T\mp\frac{\ii}{2}\overline\eta\wedge\eta,\\
 \Omega^+
 &=-\eta\wedge\Omega^T,
 &
 \Omega^-
 &=\overline\eta\wedge\Omega^T.
\end{align}
From the gravitino equation~\eqref{eq:app-KSE-modified-decoupled-2}, we find (\cf\eqref{eq:static-SU2-parallel-vector-fields})
\begin{equation}
    \nabla^H_mV_n=0,
\end{equation}
which gives (\cf \eqref{eq:static-SU2-deta})
\begin{equation}
    \nabla_{(m}\eta_{n)}=0, \quad  \nabla_{[m}\eta_{n]}=\frac{1}{2}V^aH_{amn}.
\end{equation}

One finds $\mathcal{L}_V H=0$ and (\cf\eqref{eq:dilaton-vertical-derivatives})
\begin{equation}
    \mathcal L_{\vvre}\phi=0,\quad \mathcal{L}_{\vvim}\phi=\Imp\,\kappa.
\end{equation}
By suitably contracting the dilatino-like equations with gamma matrices and (the conjugate of) $\epspm$, one finds that the $\SU(3)^\pm$-structures satisfy (\cf\eqref{eqn:symplectic}, \eqref{eq:weighted-holomorphic-forms} and Lemma~\ref{lem:Lee-forms})
\begin{subequations}\label{eq:app-su3}
	\begin{align}
		\dd (\ee^{-2\phi} \fundp) &= \ee^{-2\phi} \star H - 2\ee^{-2\phi} \Rep(\kappa\vformbar)\wedge\omega^+\,,  \label{eq:app-symplectic-6dp}\\
		\dd (\ee^{-2\phi} \fundm) &=  -\ee^{-2\phi} \star H   +   2\ee^{-2\phi} \Rep(\kappa
		\vform)\wedge\omega^-\,, \label{eq:app-symplectic-6dm} \\
		\dd (\ee^{-2\phi} \fundp\wedge\fundp) &= - 2\ee^{-2\phi} \Rep(\kappa\vformbar)\wedge\fundp\wedge\fundp\,, \label{eq:app-generalised-conformally-balanced-6dp}\\
		\dd (\ee^{-2\phi}\fundm\wedge\fundm ) &= +2\ee^{-2\phi} \Rep(\kappa\vform)\wedge\fundm\wedge\fundm\,, \label{eq:app-generalised-conformally-balanced-6dm}\\
		\dd (\ee^{-2\phi}\Omega^+) &= - \ee^{-2\phi}\kappa \vformbar\wedge\Omega^+\,, \label{eq:app-holomorphic-6dp} \\
		\dd (\ee^{-2\phi}\Omega^-)&=   \ee^{-2\phi}\kappa  \vform\wedge\Omega^-\,.\label{eq:app-holomorphic-6dm}
	\end{align}
\end{subequations}

With analogous techniques, one finds that the transverse $\SU(2)$ structure obeys (\cf\eqref{eq:transverse-holomorphic-form})
\begin{subequations}
	\begin{align}
	\dd (\ee^{-2\phi}\holot)&=  2\ii \kappa\ee^{-2\phi} \vformim \wedge \holot\,, \label{eq:app-holomorphic-4d} \\
		\Rep(\kappa)\,  \dd\left(\ee^{-2\phi} \vformre \wedge\fundt\right)&=0 \,, \label{eq:app-transverse-condition} \\
		\Rep(\kappa)\, \left[  \dd  \vformre +2  \Imp(\kappa)  \vformre \wedge \vformim\right]\wedge \fundt &= 0 \, . \label{eq:app-kappareorim}
		\end{align}
\end{subequations}

With the conditions found above, one can describe the geometry of \BHER and \BHEI\!\!. For example, when $\kappa\in\R$, one can use \eqref{eq:app-kappareorim} in \eqref{eq:app-transverse-condition} to give, together with \eqref{eq:app-holomorphic-4d}, that the transverse metric is conformally K\"ahler. The geometry of \BHER and \BHEI (in the equal chiral norms cases) is summarized in Table \ref{tab:real-vs-imaginary}.

\section{A web of embeddings}
\label{sugraembeddings}

In this appendix we show how our \BHER geometry on $M_6$ emerge in massive IIA compactifications on $\AdS_3\times I\times M_6$ and Type I compactifications on $\AdS_3\times S^1\times M_6$; moreover, we write down explicitly a T-duality from the massless limit of these massive IIA backgrounds to the Type IIB $\AdS_3\times Y_5\times T^2$ solutions of \cite{Gauntlett:2007ts}, in agreement with what was shown in \cite{Couzens:2022agr} for the solution of \cite{Donos:2008ug}. A different T-duality, relates these backgrounds to the Type I solutions on $\AdS_3\times S^1\times M_6$, showing that \BHER geometry emerges on $\AdS_3$ backgrounds of Type I, as anticipated in \cite{Papadopoulos:2026xdj}. We conclude by showing that \BHEI is the geometry of the internal space $M_6$ of an $\mathcal{N}=2$ compactification of the Type I generalized supergravity introduced in \cite{Wulff:2016tju}.

\subsection{Embedding in massive Type IIA}\label{mass2aemb}  We now show that, given a solution of our GKSE with $\kappa\in\R$, we can uplift the solution to a supersymmetric background of massive IIA on $\AdS_3\times M_6\times I_z$. To avoid confusion with the NSNS flux $H$, in this section we denote with $\Ht$ the torsion of the BHE structure. We take the Ansatz (in the string frame)
\begin{equation} \label{eq:app-mIIA-Ansatz}
\begin{aligned}
    g_{10} &=A^2\bigg(g(\AdS_3)+g_6\bigg)+\ee^{2B}\dd z^2, \quad A= \ee^{-\phi/2}h^{-1/8}, \quad \ee^B=\ee^{\phi/2}h^{-3/8}, \quad  h=c_0+2 F_0 z. \\
    \end{aligned}
\end{equation}
with the ten-dimensional dilaton and fluxes
\begin{equation} \label{eq:app-mIIA-Ansatz-fluxes}
   \ee^{-2\Phi}=\ee^\phi h^{5/4}, \quad F_2 =0,\quad H_{10}=0, \quad F_4=(2\kappa\,\vol(\AdS_3)+\Ht)\wedge \dd z.
\end{equation}
The Killing spinor equations then read\footnote{Here we are using the convention $\slashed{\alpha}_k=\frac{1}{k!}\alpha_{n_1\ldots n_k}\Gamma^{n_1\ldots n_k}$.}
\begin{align}
    \bigg(\nabla_M+\frac{\ee^\Phi}{8}(F_0+\slashed{F}_4)\Gamma_M\bigg)\varepsilon&=0,\label{eq:app-gravitino-mIIA}\\
    \bigg(\slashed{\nabla}\Phi+\frac{\ee^\Phi}{4}(5F_0+\slashed{F}_4)\bigg)\varepsilon&=0, \label{eq:app-dilatino-mIIA}
\end{align}
where we have introduced the non-chiral spinor of Type IIA $\varepsilon=\varepsilon_++\varepsilon_-$.
The frame is $e^\alpha=A\hat e^\alpha$, $e^a=A\tilde{e}^a$, $e^{\underline{z}}=\ee^B\dd z$, and $\mathrm{Cliff}(1,9)$ decomposes as\footnote{We take the generators $\hat\gamma^\alpha$ of $\mathrm{Cliff}(1,2)$ to be in the Majorana basis, and the basis of $\mathrm{Cliff}(6)$ we have already used above, where $\gamma^a$ are purely imaginary and antisymmetric; then our $\Gamma^A$ are a real Majorana basis for $\mathrm{Cliff(1,9)}$.}
\begin{align}
        \Gamma^\alpha&=\hat{\gamma}^\alpha\otimes\mathbf{1}\otimes \sigma_1, \\
    \Gamma^{\underline{z}}&=\mathbf{1}\otimes\gamma_*\otimes \sigma_2, \\
    \Gamma^a &= \mathbf{1}\otimes\gamma^a\otimes \sigma_2.
\end{align}
We find
\begin{equation}
    \ee^\Phi\slashed{F}_4=A^{-1}(2\kappa\Gamma^{012{\underline{z}}}+\slashed{\Ht}\Gamma^{\underline{z}}), \quad \slashed{\nabla}\Phi=-\frac{1}{2}A^{-1}\slashed\nabla\phi-\frac{5}{4}\ee^{\Phi}F_0\Gamma^{\underline{z}}.
\end{equation}
We take the Ansatz for the spinor $\varepsilon$
\begin{equation}\label{eq:app-mIIA-embedding-spinor-Ansatz}
    \varepsilon=A^{1/2}\zeta\otimes \eta + \mathrm{c. c.},
\end{equation}
where $\zeta$ is a Killing spinor of $\AdS_3$ and $\eta$ a doublet of complex spinors on $M_6\times I_z$, and we assume the projection
\begin{equation}\label{eq:app-mIIA-projection}
     (\gamma_*\otimes \sigma_2)\eta=\eta \implies  \Gamma^{\underline{z}}\varepsilon=\varepsilon.
\end{equation}
We define $\sigma_2u_\pm=\pm u_\pm$, so that the projection above can be written as
\begin{equation}
    \eta=\epsp\otimes u_++\ii\epsm\otimes u_-,
\end{equation}
where $\gamma_*\epspm=\pm\epspm$. Note that $\sigma_1 u_\pm=\pm\ii u_\mp$.

Since the Killing spinor equations are real, we will carry on the computation for the first term of \eqref{eq:app-mIIA-embedding-spinor-Ansatz}; the corresponding conditions for the complex conjugate are automatically satisfied. Let us start with~\eqref{eq:app-dilatino-mIIA}. The Romans mass terms combine into $\frac{5}{4}\ee^\Phi F_0(1-\Gamma^{\underline{z}})$, which vanishes on $\varepsilon$ due to~\eqref{eq:app-mIIA-projection}. The rest amounts to
\begin{equation}
    -\frac{1}{2}A^{-1}\left(\Gamma^a\nabla_a\phi-\frac{1}{2}\slashed{\Ht}-\kappa \Gamma^{012}\right)\varepsilon=-\frac{1}{2}A^{-1/2}\zeta\otimes(\mathcal{A}\otimes\sigma_2-\kappa\sigma_1)\eta=0.
\end{equation}
Noting that
\begin{equation} \label{eq:app-dilatino-from-IIAm-embedding}
    (\mathcal{A}\otimes\sigma_2-\kappa\sigma_1)\eta=(\mathcal{A}\epsp-\kappa\epsm)\otimes u_+-\ii (\mathcal{A\epsm+\kappa\epsp})\otimes u_-,
\end{equation}
this implies our~\eqref{eq:app-KSE-modified-decoupled-1}.

Moving to~\eqref{eq:app-gravitino-mIIA}, we have to take into account the contribution of the warp factor in the spin connection, namely
\begin{align}
    \nabla_\alpha &=A^{-1}\hat\nabla_\alpha+\frac{1}{2}\Gamma_\alpha(\Gamma^a A^{-1}\partial_a\log A+\Gamma^{\underline{z}}\ee^{-B}\partial_z\log A), \\
    \nabla_{\underline{z}} &= \ee^{-B}\partial_z+\frac{1}{2}\Gamma_{\underline{z}}\Gamma^a A^{-1}\partial_a B, \\
    \nabla_a &= A^{-1}(\tilde{\nabla}_a+\frac{1}{2}\Gamma_a{}^b\partial_b\log A)+\frac{1}{2}\Gamma_a\Gamma^{\underline{z}}\ee^{-B}\partial_z\log A,
\end{align}
with $\partial_a\log A=-\frac{1}{2}\partial_a\phi$, $\partial_a B=\frac{1}{2}\partial_a\phi$ and $\ee^{-B}\partial_z \log A=-\frac{1}{4}\ee^\Phi F_0$. In every component, except along $z$, the Romans mass terms combine into $\frac{1}{8}\ee^\Phi F_0\Gamma_A(1-\Gamma^{\underline{z}})$, which vanishes on $\varepsilon$. More precisely,
\begin{equation}
\begin{aligned}
    \bigg(\nabla_\alpha+\frac{\mathrm{e}^\Phi}{8}(F_0+\slashed F_4)\Gamma_\alpha\bigg)\varepsilon
&=\frac{1}{8}\ee^\Phi F_0\,\Gamma_\alpha(1-\Gamma^{\underline{z}})\varepsilon
-\frac{1}{4}A^{-1}\Gamma_\alpha\bigg(\slashed\nabla\phi-\frac{1}{2}\slashed\Ht-\kappa\Gamma^{012}\bigg)\varepsilon \\
&\quad+A^{-1/2}\bigg(\hat\nabla_\alpha\zeta-\frac{\kappa}{2}\hat\gamma_\alpha\zeta\bigg)\otimes\eta ,
\end{aligned}
\end{equation}
\begin{equation}
\begin{aligned}
\bigg(\nabla_{\underline{z}}+\frac{\ee^\Phi}{8}(F_0+\slashed F_4)\Gamma_{\underline{z}}\bigg)\varepsilon
&=-\frac{1}{8}\ee^\Phi F_0\,(1-\Gamma^{\underline{z}})\varepsilon
-\frac{1}{4}A^{-1}\bigg(\slashed\nabla\phi-\frac{1}{2}\slashed \Ht-\kappa\Gamma^{012}\bigg)\varepsilon\\
&\quad+A^{1/2}\ee^{-B}\,\zeta\otimes\partial_z\eta ,
\end{aligned}
\end{equation}
\begin{equation}
\begin{aligned}
\bigg(\nabla_a+\frac{\ee^\Phi}{8}(F_0+\slashed F_4)\Gamma_a\bigg)\varepsilon
&=\frac{1}{8}\ee^\Phi F_0\,\Gamma_a(1-\Gamma^{\underline{z}})\varepsilon
-\frac{1}{4}A^{-1}\Gamma_a\bigg(\slashed\nabla\phi-\frac{1}{2}\slashed \Ht-\kappa\Gamma^{012}\bigg)\varepsilon\\
&\quad+A^{-1/2}\,\zeta\otimes(\mathcal D_a\otimes\mathbf 1)\eta .
\end{aligned}
\end{equation}
Thus, we are left with
\begin{equation} \label{eq:app-gravitino-from-IIAm-embedding}
    \bigg(\nabla_m-\frac{1}{8}\Ht_{mnp}\gamma^{np}\bigg)\epspm=0, \quad \partial_z\eta=0, \quad \hat\nabla_\mu\zeta=\frac{\kappa}{2}\hat\gamma_\mu\zeta.
\end{equation}
The last one is the Killing spinor equation on $\AdS_3$, from which we read $\ell=1/|\kappa|$.

Therefore, we have shown in \eqref{eq:app-dilatino-from-IIAm-embedding} and \eqref{eq:app-gravitino-from-IIAm-embedding} that our GKSE, when $\kappa$ is real, naturally emerge as the KSE on the internal $M_6$ of a massive Type IIA compactification on $\AdS_3\times M_6\times I_z$.

Let us now set $F_0=0$ and $c_0=1$, then~\eqref{eq:app-mIIA-Ansatz} becomes
\begin{equation}\label{eq:app-IIAm-embedding-massless}
    g_{10}=\ee^{-\phi}\bigg(g(\AdS_3)+g_6\bigg)+\ee^{\phi}\dd z^2, \quad \Phi=-\frac{\phi}{2}, \quad F_4=\bigg(2\kappa\,\vol(\AdS_3)+\Ht\bigg)\wedge\dd z.
\end{equation}
From Theorem \ref{thm:real-kappa-local-model}, we have
\begin{equation}
     g_6=\frac{\Rk}{8\kappa^2}g_K+\eta_1^2+\eta_2^2, \quad \ee^{2\phi}=\frac{\Rk}{8\kappa^2},\quad \dd\vformre=\alpha,\quad \dd\vformim=\frac{\rhok}{2\kappa}.
\end{equation}
This is the massless limit of the massive GK solution of ~\cite[(5.33)]{Couzens:2022agr},
with the identifications
$m=\kappa,$
$\eta_2=\frac{1}{2m}(d\psi+P_K),$
$\eta_1=d\varphi+t_1\Sigma,$
where $\varphi$, $t_1$ and $\Sigma$ are those of \cite{Couzens:2022agr}.

We would now like to T-dualize along $\vvre$.
First, recall from \eqref{eq:local-static-fibration-forms} that $\vformre=\dd t_1+\sigma_1$ and $\vformim=\dd t_2+\sigma_2$, and decompose as in Section \ref{Section7}
\begin{equation}
    g_6=g_5+(\dd t_1+\sigma_1)^2, \quad g_5=\ee^{2\phi}g_K+\eta_2^2, \quad \Ht=\Ht_5+\vformre\wedge\alpha.
\end{equation}
The T-duality along $V_1$ yields the Type IIB background
\begin{align}\label{massive-GK-after-T-duality}
    g_{10} &=\ee^{-\phi}\bigg(g(\AdS_3)+g_5\bigg)+\ee^{\phi}(\dd \tilde{x}^2+\dd z^2), \\
    F_5 &=(1+\star)[(2\kappa\,\vol(\AdS_3)+\Ht_5)\wedge \dd z\wedge \dd \tilde x], \\
    F_3 &=-\alpha\wedge\dd z, \\
    H &=-\alpha\wedge\dd\tilde{x}, \\
    \Phi &=0,
\end{align}
where $\tilde x$ is the T-dual coordinate. Note that
\begin{equation}
    \star[(2\kappa\,\vol(\AdS_3)+\Ht_5)\wedge \dd z\wedge \dd \tilde x] = -2\kappa\,\ee^{-2\phi}\vol_5-\vol(\AdS_3)\wedge F_{\mathrm{GK}}, \quad F_{\mathrm{GK}}=\ee^{-2\phi}\star_5 \Ht_5.
\end{equation}
From \eqref{eq:H-real}, we have
\begin{equation}
    \Ht_5=\dd^{\mathrm{c}}(\ee^{2\phi})\wedge\omega_K+\eta_2\wedge \dd\eta_2. 
\end{equation}
Writing $\rho_K=\frac{R_K}{4}\omega_K+\rho_0$, with $\rho_0$ its anti-self-dual part, one can see $\star_4\rho_K=\frac{R_K}{2}\omega_K-\rho_K$. Therefore,
\begin{equation}
   F_{\mathrm{GK}} =\ee^{-2\phi}\star_5 H_5=\ee^{-2\phi}(-\dd \eta_2+2\dd\phi\wedge\eta_2+2\kappa\,\ee^{2\phi}\omega_K)=2\kappa \omega_K-\dd(\ee^{-2\phi}\eta_2).
\end{equation}
The fluxes in \eqref{massive-GK-after-T-duality} then become
\begin{equation}
\begin{aligned}
        F_5 &=[2\kappa\vol(\AdS_3)+\dd^{\mathrm{c}}(\ee^{2\phi})\wedge\omega_K+\eta_2\wedge\dd\eta_2]\wedge \dd z\wedge \dd\tilde x+ \\ &\hphantom{=}\vol(\AdS_3)\wedge[\dd(\ee^{-2\phi}\eta_2)-2\kappa\omega_K]-\kappa\ee^{2\phi}\omega_K\wedge\omega_K\wedge\eta_2, \\
        F_3 &=-\dd\eta_1\wedge \dd z, \\
        H &=-\dd\eta_1\wedge \dd \tilde x.
\end{aligned}    
\end{equation}
The Type IIB Bianchi identity $\dd F_5=H\wedge F_3$ then gives the ``Box equation''
\begin{equation}
    \dd\bigg(\dd ^{\mathrm{c}}(\ee^{2\phi})\wedge\omega_K+\eta_2\wedge \dd\eta_2\bigg)=-\dd\eta_1\wedge\dd\eta_1.
\end{equation}
In particular, for $\alpha=0$, we have recovered the $\AdS_3\times Y_7$ Type IIB GK solution, with $Y_7=Y_5\times T^2$, \cf \eqref{GK-T2-Y5}. When $\alpha\neq 0$, the Bianchi identity gives $\dd \Ht_5=-\alpha\wedge \alpha$, as in~\eqref{eq:global-alpha-deformed-GK-equation}; this is the transgression mechanism introduced in GK geometry in \cite{Donos:2008ug}.

\subsection{Embedding in Type I}\label{typIsugemb}  We now take \eqref{eq:app-mIIA-Ansatz} and impose that the Romans mass vanishes $F_0=0$ and $c_0=1$, obtaining again \eqref{eq:app-IIAm-embedding-massless}. If we also require that $z$ is periodic, then~\eqref{eq:app-mIIA-Ansatz} is the T-dual of 
\begin{equation}
    g_{10}=\ee^{-\phi}\bigg(g(\AdS_3)+g(M_6)+\dd y^2\bigg), \quad \Phi=-\phi, \quad F_3=\Ht+2\kappa\,\vol(\AdS_3),
\end{equation}
where $y$ is the T-dual coordinate.
But this is, in turn, the S-dual of 
\begin{equation}
    g_{10}=g(\AdS_3)+g(M_6)+\dd y^2, \quad \Phi=\phi, \quad H=\Ht+2\kappa\,\vol(\AdS_3)
\end{equation}
Therefore, the massless sector of the family~\eqref{eq:app-mIIA-Ansatz} is equivalent via T- and S-duality to the Type I $\AdS_3\times M_6\times S_1$ embedding.

\subsection{Embedding in generalized supergravity} \label{gensugemb} Here we show that our GKSE naturally arise on the geometry of the internal space of $M_6$ of a compactification of the generalized supergravity introduced in~\cite{Wulff:2016tju} on $\R^{1,3}\times M_6$. The generalized Type I Killing spinor equations (in our conventions\footnote{$H_{\mathrm{there}}=-H_{\mathrm{here}}$.}) read
\begin{subequations}
	\begin{align}
		\bigg(\Gamma^M X_M-\frac{1}{12}H_{MNP}\Gamma^{MNP}\bigg)\varepsilon &=0, \label{eq:app-generalised-KSE-dilatino} \\
		\bigg(\nabla_M-\frac{1}{8}H_{MNP}\Gamma^{NP}\bigg)\varepsilon&=0, \label{eq:app-generalised-KSE-gravitino}
	\end{align}
\end{subequations} 
where $\varepsilon$ is a Majorana--Weyl spinor. The integrability conditions of these KSE read
\begin{subequations}
	\begin{align}
		R_{MN}+2\nabla_{(M}X_{N)}-\frac{1}{4}H_{MPQ}H_{N}{}^{PQ}=0, \\
		\nabla^SH_{SMN}-2X^SH_{SMN}+4\nabla_{[M}X_{N]}=0, \\
		\nabla^M X_M-2 X^M X_M+\frac{1}{12}H_{MNP}H^{MNP}=0,
	\end{align}
\end{subequations}
with the Bianchi identity
\begin{equation}
	\dd H=0.
\end{equation}

When $X_M=\nabla_M\phi$, these are the standard Type I equations, while if $\dd X\neq0$, no dilaton $X=\dd\Phi$ exists even locally, and the solution is genuinely generalized.

Now we take the direct product Ansatz for the metric
\begin{equation}
	g_{10}=g(\R^{1,3})+g(M_6),
\end{equation}
$H$ to be purely internal on $M_6$ and we define
\begin{equation}
	X^{\pm}=\frac{1}{2}\theta^{\pm}=\dd\phi+W^{\pm}, \text{ where } W^\pm= \mp \Rep(\kappa) \eta_1-\Imp(\kappa) \eta_2.
\end{equation}
The Clifford algebra decomposes as\footnote{As we did before, here $\hat\gamma_a$ are taken to be the generators of $\mathrm{Cliff}(1,3)$ in the Majorana basis, $\gamma_a$ are purely imaginary, therefore $\Gamma_a$ are purely real.}
\begin{equation}
	\Gamma_\mu=\hat\gamma_\mu\otimes \mathbf{1},\quad \Gamma_m=\hat\gamma_*\otimes \gamma_m,
\end{equation}
so $\Gamma_*=\hat\gamma_*\otimes\gamma_*$, $\Gamma^{np}=\mathbf{1}\otimes \gamma^{np}$ and $\Gamma^{mnp}=\hat\gamma_*\otimes \gamma^{mnp}$.
When $\kappa=\ii\lambda$, with $\lambda\in\R$, we have
\begin{equation}\label{eq:app-generalised-sugra-X-imaginary-case}
    X^+=X^-=\dd\phi-\lambda\eta_2\coloneqq X,
\end{equation}
and we take the $\mathcal{N}=2$ spinor Ansatz
\begin{align}
	\varepsilon=\zeta_1\otimes \epsilon_+ +\zeta_2\otimes (\epsm)^*+ \mathrm{c.c.}, 
\end{align}
where $\zeta_1$ and $\zeta_2$ are constant spinors on $\R^{1,3}$ of positive chirality, namely $\hat\gamma_*\zeta_{1,2}=\zeta_{1,2}.$
Using this Ansatz as input in~\eqref{eq:app-generalised-KSE-dilatino}, it gives
\begin{equation}
	\left(\gamma^m X_m-\frac{1}{12}H_{mnp}\gamma^{mnp}\right)\epspm=0,
\end{equation}
which, using~\eqref{eq:static-vector-Clifford-action}, becomes exactly our~\eqref{eq:app-KSE-modified-1}. 

Analogously, from~\eqref{eq:app-generalised-KSE-gravitino} we find~\eqref{eq:app-KSE-modified-2}. Therefore, there is a natural embedding of our BHE in the $\kappa\in\ii\R$ as $\mathcal{N}=2$ solutions of generalized supergravity on $\R^{1,3}\times M_6$.

\subsubsection{Equations of motion} The equations of motion reduce to the ones on $M_6$
\begin{subequations}
	\begin{align}
		R_{mn}+2\nabla_{(m}X_{n)}-\frac{1}{4}H_{mpq}H_{n}{}^{pq}=0, \\
		\nabla^sH_{smn}-2X^sH_{smn}+4\nabla_{[m}X_{n]}=0, \\
		\nabla^m X_m-2 X^m X_m+\frac{1}{12}H_{mnp}H^{mnp}=0,
	\end{align}
\end{subequations}
One can check that substituting~\eqref{eq:app-generalised-sugra-X-imaginary-case} in these equations of motion yields our integrability conditions~\eqref{eq:app-EOM} and~\eqref{eq:app-H-flux-EOM} with $\kappa=\ii\lambda$.

\medskip

\noindent\textbf{Acknowledgments.}  
Anna Fino is partially supported by GNSAGA (INdAM).

\bibliographystyle{plain}
\bibliography{BCFM-bib}
\end{document}